# Rotated reference frames in radiative transport theory

Manabu Machida*
Department of Informatics, Faculty of Engineering, Kindai University, Higashi-Hiroshima, Japan
*Corresponding author: email address: machida@hiro.kindai.ac.jp

**Abstract**
Rotated reference frames offer fast algorithms for the radiative transport equation (RTE). We review the singular-eigenfunction approach and related numerical methods for the multi-dimensional RTE with rotated reference frames.

**Contents**

## 1. Introduction

In this review, we consider the radiative transport equation (RTE). The equation governs transport phenomena in which particles are scattered and absorbed in a medium but the scattering among the particles can be neglected. The RTE is a type of the linear Boltzmann equation. In the RTE, the speed of particles is quite often assumed to be a constant. Light propagation in random media is governed by the RTE since the light propagation can be viewed as moving particles which undergo scattering and absorption (Ishimaru, 1978). See Carminati and Schotland, 2021 for the relation between Maxwell's equations and the RTE. Random media (i.e., media in which scattering and absorption occur randomly) include fog, cloud, biological tissue, and the interstellar medium. Moreover, neutrons in a reactor obey the RTE (Duderstadt and Martin, 1979). In most parts of this review, we assume that the absorption coefficient $\mu_a$ and scattering coefficient $\mu_s$ are constants.

The singular-eigenfunction approach was explored as early as 1945 by Davison, 1945. After some further efforts such as Van Kampen, 1955, Davison, 1957, and Wigner, 1961, Case established a way of finding solutions with separation of variables (Case, 1959; Case, 1960). In 1960, Case considered the time-independent one-dimensional radiative transport equation with isotropic scattering and solved the equation with separation of variables by finding singular eigenfunctions (Case, 1960). The method was soon extended to the case of anisotropic scattering without (McCormick and Kuscer, 1965; Mika, 1961) and with (McCormick and Kuscer, 1966) azimuthal dependence. In this singular-eigenfunction approach, solutions to the one-dimensional radiative transport equation are given by a superposition of singular eigenfunctions. The existence and uniqueness of such solutions were proved (Larsen, 1974; Larsen, 1975; Larsen and Habetler, 1973; Larsen, Sancaktar, and Zweifel, 1975). Case's approach was revisited recently and the RTE for isotropic scattering was considered in half spaces of arbitrary dimension (d'Eon and McCormick, 2019).

In 1964, Dede used rotated reference frames to solve the three-dimensional radiative transport equation with the spherical-harmonic approximation (cf., the $P_N$ method) (Dede, 1964). Dede pointed out that equations in three dimensions reduce to one-dimensional equations if reference frames are rotated in the direction of the Fourier vector. Kobayashi developed Dede's calculation and computed coefficients in the spherical-harmonic expansion by solving a three-term recurrence relation recursively starting with the initial term (Kobayashi, 1977). Numerical calculation was performed for a two-dimensional rectangular domain in Kobayashi, 1977.

In 2004, Markel independently arrived at rotated reference frames and obtained the coefficients in terms of eigenvalues and eigenvectors of the tridiagonal matrix originating from the three-term recurrence relation (Markel, 2004). By computing the specific intensity for a three-dimensional infinite medium, in Markel, 2004, it was showed that the algorithm is efficient (Markel, 2004). Moreover it was found that complex unit vectors can be used to rotate reference frames (Panasyuk, Schotland, and Markel, 2006). The new formulation can be viewed as separation of variables in which the eigenvalues are separation constants (Schotland and Markel, 2007). This generalization made it possible to solve boundary value problems in the form of plane-wave decomposition (Machida, Panasyuk, Schotland, and Markel, 2010). The series of numerical algorithms started with Markel, 2004 was named the method of rotated reference frames.

Until this point, the notion of rotated reference frames had always come together with the spherical-harmonic expansion of the specific intensity. Then it was realized that the technique of rotated reference frames is not necessarily linked to spherical harmonics, and singular eigenfunctions in the one-dimensional transport theory were extended to three dimensions (Machida, 2014). In particular, the fundamental solution with singular eigenfunctions of the one-dimensional RTE for anisotropic scattering could be extended to three dimensions. Indeed, the method of rotated reference frames is a three-dimensional extension of the spherical-harmonic expansion (Barichello, Garcia, and Siewert, 1998; Sanchez and McCormick, 1982) in Caseology (the singular-eigenfunction approach started by Case). Using singular eigenfunctions, the fundamental solution of the RTE was obtained in flatland (i.e., two-dimensional space) (Machida, 2016b).

## 2. Radiative transport equation

Let $\mathbf{r}$ be the position vector in three dimensions and $\hat{\mathbf{s}}$ be a unit vector which specifies the direction of light that propagates in a random medium. We used the hat symbol $\widehat{\ }$ to indicate that $\hat{\mathbf{s}}$ is a unit vector ($\hat{\mathbf{s}}\cdot\hat{\mathbf{s}}=1$). In this article, a vector with the hat symbol is a unit vector. The vector $\mathbf{r}$ is written as

$$\mathbf{r}=\begin{pmatrix}x\\y\\z\end{pmatrix},\quad -\infty<x<\infty,\quad -\infty<y<\infty,\quad -\infty<z<\infty\,. \tag{2.1}$$

The vector $\hat{\mathbf{s}}$ is given by

$$\hat{\mathbf{s}}=\begin{pmatrix}\boldsymbol{\omega}\\\cos\theta\end{pmatrix},\quad \boldsymbol{\omega}=\begin{pmatrix}\sin\theta\cos\varphi\\\sin\theta\sin\varphi\end{pmatrix},\quad 0\le\theta\le\pi\,,\quad 0\le\varphi<2\pi\,. \tag{2.2}$$

In addition to the set $\mathbb{R}$ of real numbers and the sets $\mathbb{R}^d$ of $d$-dimensional real vectors ($d=2,3$), we will use notation:

$$\mathbb{S}^2=\{(\theta,\varphi);\ 0\le\theta\le\pi\,,\ 0\le\varphi<2\pi\}\,. \tag{2.3}$$

The RTE is written as

$$\frac{1}{c}\frac{\partial}{\partial t}I(\mathbf{r},\hat{\mathbf{s}},t)+\hat{\mathbf{s}}\cdot\nabla I(\mathbf{r},\hat{\mathbf{s}},t)+(\mu_a+\mu_s)I(\mathbf{r},\hat{\mathbf{s}},t)=\mu_s\int_{\mathbb{S}^2}p(\hat{\mathbf{s}},\hat{\mathbf{s}}')I(\mathbf{r},\hat{\mathbf{s}}',t)d\hat{\mathbf{s}}'+S(\mathbf{r},\hat{\mathbf{s}},t)\,, \tag{2.4}$$

where $t$ is time and $I(\mathbf{r},\hat{\mathbf{s}},t)$ is the specific intensity at position $\mathbf{r}$ in direction $\hat{\mathbf{s}}$ at time $t$. The speed of light in the medium, which is assumed to be a positive constant, is denoted by $c$. In the second term on the left-hand side, the nabla $\nabla$ means

$$\nabla=\begin{pmatrix}\partial/\partial x\\\partial/\partial y\\\partial/\partial z\end{pmatrix}. \tag{2.5}$$

In general, the absorption coefficient $\mu_a$ and scattering coefficient $\mu_s$ are functions of $\mathbf{r},\hat{\mathbf{s}}$. As described above, we will assume that $\mu_a,\mu_s$ are constants unless otherwise stated. The sum of these optical parameters is called the total attenuation $\mu_t=\mu_a+\mu_s$. The scattering phase function $p(\hat{\mathbf{s}},\hat{\mathbf{s}}')$, described below, can depend on $\mathbf{r}$ but this dependence is neglected. On the right-hand side, $S(\mathbf{r},\hat{\mathbf{s}},t)$ is the source term. The RTE (2.4) satisfies suitable boundary conditions and an initial condition. In the time-independent case, the RTE is written as

$$\hat{\mathbf{s}}\cdot\nabla I(\mathbf{r},\hat{\mathbf{s}})+\mu_t I(\mathbf{r},\hat{\mathbf{s}})=\mu_s\int_{\mathbb{S}^2}p(\hat{\mathbf{s}},\hat{\mathbf{s}}')I(\mathbf{r},\hat{\mathbf{s}}')d\hat{\mathbf{s}}'+S(\mathbf{r},\hat{\mathbf{s}})\,. \tag{2.6}$$

The stationary propagation is governed by (2.6). The time-dependent RTE (2.4) is reduced to the RTE (2.6) by the Fourier transform or Laplace transform. In this sense, it is crucial to solve (2.6). The main target of this review is the RTE in (2.6). Boundary conditions are imposed when (2.6) is solved.

The scattering phase function $p(\hat{\mathbf{s}}, \hat{\mathbf{s}}')$ describes the probability that light propagating in direction $\hat{\mathbf{s}}'$ changes its direction to $\hat{\mathbf{s}}$ by scattering. We have

$$\int_{\mathbb{S}^2} p(\hat{\mathbf{s}}, \hat{\mathbf{s}}') d\hat{\mathbf{s}} = \int_0^{2\pi} \int_0^{\pi} p(\hat{\mathbf{s}}, \hat{\mathbf{s}}') \sin\theta \; d\theta \, d\varphi = 1 \,. \tag{2.7}$$

The anisotropic factor or scattering asymmetry parameter $\mathrm{g} \in [-1,1]$ is introduced as the average of $\hat{\mathbf{s}} \cdot \hat{\mathbf{s}}'$:

$$\mathrm{g} = \int_{\mathbb{S}^2} \hat{\mathbf{s}} \cdot \hat{\mathbf{s}}' p(\hat{\mathbf{s}}, \hat{\mathbf{s}}') d\hat{\mathbf{s}}' \,. \tag{2.8}$$

In the case of isotropic scattering, where a particle is scattered in all directions with equal probability, $p(\hat{\mathbf{s}}, \hat{\mathbf{s}}') = 1/(4\pi)$.

For a unit vector $\hat{\mathbf{k}}$ (i.e., $\hat{\mathbf{k}} \cdot \hat{\mathbf{k}} = \hat{k}_x^2 + \hat{k}_y^2 + \hat{k}_z^2 = 1$), we introduce an operator $\mathcal{R}_{\hat{\mathbf{k}}}$ that rotates the reference frame so that the $z$-axis lies in the direction of $\hat{\mathbf{k}}$. For example, we can write

$$\mathcal{R}_{\hat{\mathbf{k}}} \cos\theta = \mathcal{R}_{\hat{\mathbf{k}}} \hat{\mathbf{s}} \cdot \hat{\mathbf{z}} = \hat{\mathbf{s}} \cdot \hat{\mathbf{k}} \,, \tag{2.9}$$

where $\hat{\mathbf{z}}$ is a unit vector in the direction of the positive $z$-axis. Later, we will define the operator $\mathcal{R}_{\hat{\mathbf{k}}}$ using Wigner's $D$-matrices. But for now, let us understand $\mathcal{R}_{\hat{\mathbf{k}}}$ rather intuitively.

# 3. The first example of rotated reference frames for the RTE

Let us begin with the following homogeneous RTE.

$$\hat{\mathbf{s}} \cdot \nabla I(\mathbf{r}, \hat{\mathbf{s}}) + \mu_t I(\mathbf{r}, \hat{\mathbf{s}}) = \frac{\mu_s}{4\pi} \int_{\mathbb{S}^2} I(\mathbf{r}, \hat{\mathbf{s}}) \, d\hat{\mathbf{s}} \,. \tag{3.1}$$

We assume that $\mu_t, \mu_s$ are positive constants. To find the solution to (3.1), we assume the separated form as (Machida, 2016c)

$$I(\mathbf{r}, \hat{\mathbf{s}}) = e^{-\mu_t \hat{\mathbf{k}} \cdot \mathbf{r}/\nu} \mathcal{R}_{\hat{\mathbf{k}}} \Phi_\nu(\hat{\mathbf{s}}), \tag{3.2}$$

where $\hat{\mathbf{k}}$ is a unit vector and $\nu$ is a separation constant. Recall that $\mathcal{R}_{\hat{\mathbf{k}}}$ was introduced in the end of Sec. 2. Here, $\Phi_\nu(\hat{\mathbf{s}})$ is some function which depends on $\hat{\mathbf{s}}$ but is independent of $\mathbf{r}$. We impose the normalization condition as

$$\int_{\mathbb{S}^2} \Phi_\nu(\hat{\mathbf{s}}) d\hat{\mathbf{s}} = 2\pi \,. \tag{3.3}$$

Since the integral over all angles does not change by the choice of reference frames, we have

$$\int_{\mathbb{S}^2} \mathcal{R}_{\hat{\mathbf{k}}} \Phi_\nu(\hat{\mathbf{s}}) d\hat{\mathbf{s}} = 2\pi \,. \tag{3.4}$$

By substituting the separated form (3.2) for $I$ in the homogeneous equation (3.1), we obtain

$$\mu_t \left(1 - \frac{\hat{\mathbf{s}} \cdot \hat{\mathbf{k}}}{\nu}\right) \mathcal{R}_{\hat{\mathbf{k}}} \Phi_\nu(\hat{\mathbf{s}}) = \frac{\mu_s}{2} \,. \tag{3.5}$$

The left-hand side of (3.5) can be rewritten as

$$\mu_t \left(1 - \frac{\hat{\mathbf{s}} \cdot \hat{\mathbf{k}}}{\nu}\right) \mathcal{R}_{\hat{\mathbf{k}}} \Phi_\nu(\hat{\mathbf{s}}) = \mu_t \left(1 - \frac{\mathcal{R}_{\hat{\mathbf{k}}} \cos\theta}{\nu}\right) \mathcal{R}_{\hat{\mathbf{k}}} \Phi_\nu(\hat{\mathbf{s}}) = \mathcal{R}_{\hat{\mathbf{k}}} \, \mu_t \left(1 - \frac{\cos\theta}{\nu}\right) \Phi_\nu(\hat{\mathbf{s}}) \,, \tag{3.6}$$

where we used (2.9). Hence by the inverse rotation,

$$\mu_t \left(1 - \frac{\cos\theta}{\nu}\right) \Phi_\nu(\hat{\mathbf{s}}) = \frac{\mu_s}{2} \,. \tag{3.7}$$

Let us introduce the notation

$$\mu = \cos\theta \,, \qquad \varpi = \frac{\mu_s}{\mu_t} \,. \tag{3.8}$$

Then we arrive at

$$\left(1 - \frac{\mu}{\nu}\right) \Phi_\nu(\hat{\mathbf{s}}) = \frac{\varpi}{2} \,. \tag{3.9}$$

This equation (3.9) has been intensively studied in the one-dimensional transport theory and $\Phi_\nu(\hat{\mathbf{s}})$ is called the singular eigenfunction (see Appendix A). The separation constant $\nu$ is known to be continuous in $(-1,1)$ and has two discrete values $\pm\nu_0$ ($\nu_0 > 1$), where they are solutions to the following transcendental equation (Case and Zweifel, 1967).

$$1 - \varpi\nu_0 \tanh^{-1}\frac{1}{\nu_0} = 0 \,. \tag{3.10}$$

Thus, for each pair of $(\hat{\mathbf{k}}, \nu)$, we obtain the solution $I(\mathbf{r}, \hat{\mathbf{s}})$ in (3.2) of the three-dimensional RTE in (3.1) using the knowledge of the one-dimensional transport theory. The choice of the unit vector $\hat{\mathbf{k}}$ is arbitrary at this moment. Later, when we obtain the fundamental solution or solve boundary-value problems, suitable $\hat{\mathbf{k}}$ will be chosen. The unit vector $\hat{\mathbf{k}}$ is not necessarily a real vector. Indeed, as described in Sec. 4, we will consider complex vectors $\hat{\mathbf{k}} \in \mathbb{C}^3$ ($\mathbb{C}$ denotes the set of complex numbers).

# 4. Fundamental solution to the three-dimensional RTE

## 4.1. Separated solution

Let us consider the case of anisotropic scattering, i.e., $p(\hat{\mathbf{s}}, \hat{\mathbf{s}}')$ is not a constant. We assume that scatterers are spherically symmetric and model $p(\hat{\mathbf{s}}, \hat{\mathbf{s}}')$ as

$$p(\hat{\mathbf{s}}, \hat{\mathbf{s}}') = \frac{1}{4\pi} \sum_{l=0}^{l_{\max}} \beta_l P_l(\hat{\mathbf{s}} \cdot \hat{\mathbf{s}}') = \sum_{l=0}^{l_{\max}} \sum_{m=-l}^{l} \frac{\beta_l}{2l+1} Y_{lm}(\hat{\mathbf{s}}) Y_{lm}^*(\hat{\mathbf{s}}') \,, \tag{4.1}$$

where $\beta_l$ are coefficients. Here, * means complex conjugate. We note the orthogonality relation for Legendre polynomials:

$$\int_{-1}^{1} P_l(\mu) P_{l'}(\mu) \, d\mu = \frac{2}{2l+1} \delta_{ll'} \,, \tag{4.2}$$

where $\delta_{ll'}$ is the Kronecker delta. We see $\beta_0 = 1$ from the calculation,

$$1 = \int_{\mathbb{S}^2} p(\hat{\mathbf{s}}, \hat{\mathbf{s}}') \, d\hat{\mathbf{s}}' = \frac{1}{4\pi} \sum_{l=0}^{l_{\max}} \beta_l \int_{\mathbb{S}^2} P_l(\hat{\mathbf{s}} \cdot \hat{\mathbf{s}}') d\hat{\mathbf{s}}' = \beta_0 \,. \tag{4.3}$$

Moreover since $p(\hat{\mathbf{s}}, \hat{\mathbf{s}}')$ is a probability density and has to be nonnegative,

$$1 = \int_{\mathbb{S}^2} p(\hat{\mathbf{s}}, \hat{\mathbf{s}}') \, d\hat{\mathbf{s}}' \geq \left| \int_{\mathbb{S}^2} p(\hat{\mathbf{s}}, \hat{\mathbf{s}}') P_l(\hat{\mathbf{s}} \cdot \hat{\mathbf{s}}') d\hat{\mathbf{s}}' \right| = \frac{|\beta_l|}{2l+1} \,, \quad 1 \leq l \leq l_{\max} \,. \tag{4.4}$$

That is,

$$|\beta_l| \leq 2l + 1 \,. \tag{4.5}$$

We can write $\beta_1 = 3\mathrm{g}$, where g was given in (2.8). In the case of isotropic scattering, $\beta_0 = 1$ and $\beta_l = 0$ for $l \geq 1$. Spherical harmonics $Y_{lm}(\hat{\mathbf{s}})$ are given by

$$Y_{lm}(\hat{\mathbf{s}}) = \sqrt{\frac{2l+1}{4\pi}\frac{(l-m)!}{(l+m)!}} P_l^m(\mu) e^{im\varphi}\,, \tag{4.6}$$

where $P_l^m(\mu)$ are associated Legendre polynomials. They satisfy

$$\mu P_l^m(\mu) = \frac{l-m+1}{2l+1} P_{l+1}^m(\mu) + \frac{l+m}{2l+1} P_{l-1}^m(\mu)\,, \tag{4.7}$$

$$\int_{-1}^{1} P_l^m(\mu) P_{l'}^m(\mu)\, d\mu = \frac{2(l+m)!}{(2l+1)(l-m)!}\delta_{ll'}\,. \tag{4.8}$$

Moreover,

$$P_m^m(\mu) = (-1)^m (2m-1)!!\,(1-\mu^2)^{m/2}\,,\quad P_{m+1}^m(\mu) = (2m+1)\mu P_m^m(\mu)\,,\quad m \geq 0\,. \tag{4.9}$$

The Henyey-Greenstein phase function $p_{\mathrm{H-G}}$ is given by (Henyey and Greenstein, 1941)

$$\beta_l = (2l+1)\mathrm{g}^l\,,\quad l_{\max} = \infty\,. \tag{4.10}$$

That is,

$$p_{\mathrm{H-G}}(\hat{\mathbf{s}}, \hat{\mathbf{s}}') = \sum_{l=0}^{\infty}\sum_{m=-l}^{l} \mathrm{g}^l Y_{lm}(\hat{\mathbf{s}}) Y_{lm}^*(\hat{\mathbf{s}}') = \frac{1}{4\pi}\frac{1-\mathrm{g}^2}{(1+\mathrm{g}^2-2\mathrm{g}\hat{\mathbf{s}}\cdot\hat{\mathbf{s}}')^{3/2}}\,. \tag{4.11}$$

Let us consider the homogeneous RTE:

$$\hat{\mathbf{s}}\cdot\nabla I(\mathbf{r},\hat{\mathbf{s}}) + \mu_t I(\mathbf{r},\hat{\mathbf{s}}) = \mu_s \int_{\mathbb{S}^2} p(\hat{\mathbf{s}},\hat{\mathbf{s}}') I(\mathbf{r},\hat{\mathbf{s}}')\, d\hat{\mathbf{s}}'\,. \tag{4.12}$$

We assume the following separated solution.

$$I(\mathbf{r},\hat{\mathbf{s}}) = e^{-\mu_t \hat{\mathbf{k}}\cdot\mathbf{r}/\nu} \mathcal{R}_{\hat{\mathbf{k}}} \Phi_\nu^m(\hat{\mathbf{s}}). \tag{4.13}$$

Here,

$$\Phi_\nu^m(\hat{\mathbf{s}}) = \phi^m(\nu,\mu)(1-\mu^2)^{|m|/2} e^{im\varphi}\,, \tag{4.14}$$

where $\phi^m(\nu,\mu)$ satisfies

$$\int_{-1}^{1} \phi^m(\nu,\mu)(1-\mu^2)^{|m|}\, d\mu = 1\,. \tag{4.15}$$

The assumed form (4.13) has a similar structure to the separated solution in the one-dimensional transport theory by (Case, 1960; McCormick and Kuščer, 1966; Mika, 1961). For the singular eigenfunctions, not only orthogonality relations but also biorthogonality relations are known (McCormick and Kuščer, 1966). We refer the reader to McCormick and Kuščer, 1973 and Kuščer and McCormick, 1991 (the Busbridge polynomials were used).

By substituting the separated solution (4.13) for $I(\mathbf{r},\hat{\mathbf{s}})$ in the homogeneous equation (4.12), we obtain

$$\left(1 - \frac{\mathcal{R}_{\hat{\mathbf{k}}}\mu}{\nu}\right) \mathcal{R}_{\hat{\mathbf{k}}} \Phi_\nu^m(\hat{\mathbf{s}}) = \varpi \int_{\mathbb{S}^2} p(\hat{\mathbf{s}},\hat{\mathbf{s}}') \mathcal{R}_{\hat{\mathbf{k}}} \Phi_\nu^m(\hat{\mathbf{s}}')\, d\hat{\mathbf{s}}'\,. \tag{4.16}$$

Noting the fact that $p(\hat{\mathbf{s}},\hat{\mathbf{s}}')$ depends only on the dot product $\hat{\mathbf{s}}\cdot\hat{\mathbf{s}}'$ and $p(\hat{\mathbf{s}},\hat{\mathbf{s}}') = p(\mathcal{R}_{\hat{\mathbf{k}}}\hat{\mathbf{s}}, \mathcal{R}_{\hat{\mathbf{k}}}\hat{\mathbf{s}}')$ holds, (4.16) reduces to

$$\left(1 - \frac{\mu}{\nu}\right) \Phi_\nu^m(\hat{\mathbf{s}}) = \varpi \int_{\mathbb{S}^2} p(\hat{\mathbf{s}},\hat{\mathbf{s}}') \Phi_\nu^m(\hat{\mathbf{s}}')\, d\hat{\mathbf{s}}'\,. \tag{4.17}$$

Equation (4.17) has been intensively studied in the singular-eigenfunction approach for the one-dimensional transport equation.

In what follows, we will introduce $g_l^m(\nu)$. We multiply (4.17) by $Y_{lm}^*(\hat{\mathbf{s}})$ and integrate over $\hat{\mathbf{s}}$:

$$\int_{\mathbb{S}^2}\left(1-\frac{\mu}{\nu}\right)Y_{lm}^*(\hat{\mathbf{s}})\Phi_\nu^m(\hat{\mathbf{s}})d\hat{\mathbf{s}}=\varpi\int_{\mathbb{S}^2}Y_{lm}^*(\hat{\mathbf{s}})\int_{\mathbb{S}^2}p(\hat{\mathbf{s}},\hat{\mathbf{s}}')\Phi_\nu^m(\hat{\mathbf{s}}')d\hat{\mathbf{s}}'\,d\hat{\mathbf{s}}\,. \tag{4.18}$$

Thus, noting (4.7),

$$\sqrt{(l+1)^2-m^2}g_{l+1}^m(\nu)+\sqrt{l^2-m^2}g_{l-1}^m(\nu)=\nu h_l g_l^m(\nu)\,, \tag{4.19}$$

where

$$h_l=2l+1-\varpi\beta_l\Theta(l_{\max}-l)\,, \tag{4.20}$$

$$g_l^m(\nu)=(-1)^m\sqrt{\frac{(l-m)!}{(l+m)!}}\int_{-1}^1 P_l^m(\mu)\phi^m(\nu,\mu)(1-\mu^2)^{|m|/2}\,d\mu\,. \tag{4.21}$$

Here, $\Theta$ is the step function defined as $\Theta(x)=1$ for $x\geq 0$ and $\Theta(x)=0$ for $x<0$. By using (4.9),

$$g_{m+1}^m(\nu)=\frac{\nu h_m}{\sqrt{2l+1}}g_m^m(\nu)\,,\quad m\geq 0. \tag{4.22}$$

Moreover, by (4.9), (4.15), (4.21), we have

$$g_m^m(\nu)=\frac{(2m-1)!!}{\sqrt{(2m)!}}=\frac{\sqrt{(2m)!}}{2^m m!}\,,\quad m\geq 0\,. \tag{4.23}$$

The polynomials which are defined by the three-term recurrence relation (4.19) are called the normalized Chandrasekhar polynomials (Garcia and Siewert, 1989; Garcia and Siewert, 1990). See Chandrasekhar, 1950 for the Chandrasekhar polynomials. We note that

$$g_l^{-m}(\nu)=(-1)^m g_l^m(\nu)\,,\quad g_l^m(-\nu)=(-1)^{l+m}g_l^m(\nu)\,. \tag{4.24}$$

See Ganapol, 2014 and Garcia and Siewert, 1990 for the numerical calculation of Chandrasekhar's polynomials.

Singular eigenfunctions $\phi^m(\nu,\mu)$ are obtained as (see Appendix A)

$$\phi^m(\nu,\mu)=\frac{\varpi\nu}{2}\mathcal{P}\frac{g^m(\nu,\mu)}{\nu-\mu}+\lambda^m(\nu)(1-\mu^2)^{-|m|}\delta(\nu-\mu)\,, \tag{4.25}$$

where $\mathcal{P}$ denotes Cauchy's principal value and $\delta(\cdot)$ denotes Dirac's delta function. Here we introduced

$$g^m(\nu,\mu)=\sum_{l=|m|}^{l_{\max}}\beta_l p_l^m(\mu)g_l^m(\nu)\,, \tag{4.26}$$

where using associated Legendre polynomials $P_l^m(\mu)$, polynomials $p_l^m(\mu)$ are defined as

$$p_l^m(\mu)=(-1)^m\sqrt{\frac{(l-m)!}{(l+m)!}}P_l^m(\mu)(1-\mu^2)^{-|m|/2}\,. \tag{4.27}$$

The polynomials $p_l^m(\mu)$ satisfy the following three-term recurrence relation.

$$\sqrt{(l+1)^2-m^2}p_{l+1}^m(\mu)+\sqrt{l^2-m^2}p_{l-1}^m(\mu)=(2l+1)\mu p_l^m(\mu) \tag{4.28}$$

with initial terms,

$$p_m^m(\mu)=\frac{(2m-1)!!}{\sqrt{(2m)!}}=\frac{\sqrt{(2m)!}}{2^m m!}\,,\quad p_{m+1}^m(\mu)=\sqrt{2l+1}\mu p_m^m(\mu)\,,\quad m\geq 0. \tag{4.29}$$

The relation $p_l^{-m}(\mu)=(-1)^m p_l^m(\mu)$ can be used to derive formulae for $p_l^m(\mu), m<0$.

Let us call the separated solutions (4.13) eigenmodes.

### 4.2. Complex unit vectors

So far, $\hat{\mathbf{k}}$ has been an arbitrary unit vector. To proceed further, we specify the unit vector $\hat{\mathbf{k}}$. We write

$$\boldsymbol{\rho} = \begin{pmatrix} x \\ y \end{pmatrix} , \tag{4.30}$$

and introduce

$$\mathbf{q} = \begin{pmatrix} q_x \\ q_y \end{pmatrix} = q \begin{pmatrix} \cos \varphi_{\mathbf{q}} \\ \sin \varphi_{\mathbf{q}} \end{pmatrix} \tag{4.31}$$

for $-\infty < q_x < \infty$, $-\infty < q_y < \infty$, $q \geq 0$, and $0 \leq \varphi_{\mathbf{q}} < 2\pi$. We set

$$\hat{\mathbf{k}} = \hat{\mathbf{k}}(\nu, \mathbf{q}) = \begin{pmatrix} -i\nu\mathbf{q} \\ \hat{k}_z \end{pmatrix} . \tag{4.32}$$

By this, eigenmodes have the form of plane-wave decomposition:

$$I(\mathbf{r}, \hat{\mathbf{s}}) = e^{i\mu_t \mathbf{q}\cdot\boldsymbol{\rho}} e^{-\mu_t \hat{k}_z z/\nu} \mathcal{R}_{\hat{\mathbf{k}}(\nu,\mathbf{q})} \Phi_\nu^m(\hat{\mathbf{s}}). \tag{4.33}$$

By the relation $\hat{\mathbf{k}} \cdot \hat{\mathbf{k}} = 1$, we obtain

$$\hat{k}_z = \hat{k}_z(\nu q) = \sqrt{1 + (\nu q)^2} . \tag{4.34}$$

That is, $\hat{k}_z$ is real. The eigenmodes decay in the $z$ direction as $e^{-\mu_t \hat{k}_z z/\nu}$ if $\nu$ is positive and they grow if $\nu$ is negative.

### 4.3. Orthogonality relations

Rotated singular eigenfunctions $\mathcal{R}_{\hat{\mathbf{k}}} \Phi_\nu^m(\hat{\mathbf{s}})$ satisfy the following orthogonality relations.

$$\int_{\mathbb{S}^2} \mu \left( \mathcal{R}_{\hat{\mathbf{k}}_1} \Phi_{\nu_1}^{m_1}(\hat{\mathbf{s}}) \right) \left( \mathcal{R}_{\hat{\mathbf{k}}_2} \Phi_{\nu_2}^{m_2 *}(\hat{\mathbf{s}}) \right) d\hat{\mathbf{s}} = 2\pi \hat{k}_z(\nu_1 q) \mathcal{N}^m(\nu_1) \delta_{\nu_1 \nu_2} \delta_{m_1 m_2} . \tag{4.35}$$

Here, $\mathcal{N}^m(\nu)$ is the normalization factor from the one-dimensional transport theory (see Appendix A):

$$\mathcal{N}^m(\nu) = \int_{-1}^{1} \mu \, [\phi^m(\nu, \mu)]^2 (1 - \mu^2)^{|m|} d\mu . \tag{4.36}$$

The orthogonality relation (4.35) can be derived as follows. First, we can show by direct calculation,

$$\mathcal{R}_{\hat{\mathbf{k}}(\nu,\mathbf{q})} \mu = \hat{k}_z(\nu q)\mu - i\nu q \sqrt{1 - \mu^2} \cos(\varphi - \varphi_{\mathbf{q}}) , \tag{4.37}$$

and

$$\mathcal{R}_{\hat{\mathbf{k}}(\nu,\mathbf{q})}^{-1} \mu = \hat{k}_z(\nu q)\mu - i|\nu q| \sqrt{1 - \mu^2} \cos \varphi . \tag{4.38}$$

We note that the homogeneous equation (4.16) can be expressed as

$$\left(1 - \frac{\mathcal{R}_{\hat{\mathbf{k}}} \mu}{\nu}\right) \mathcal{R}_{\hat{\mathbf{k}}} \Phi_\nu^m(\hat{\mathbf{s}}) = \varpi \sum_{l=0}^{l_{\max}} \sum_{m=-l}^{l} \frac{\beta_l}{2l+1} Y_{lm}(\hat{\mathbf{s}}) \int_{\mathbb{S}^2} Y_{lm}^*(\hat{\mathbf{s}}') \mathcal{R}_{\hat{\mathbf{k}}} \Phi_\nu^m(\hat{\mathbf{s}}') d\hat{\mathbf{s}}' . \tag{4.39}$$

Using the above equation we can write two equations:

$$\left( \mathcal{R}_{\hat{\mathbf{k}}_2} \Phi_{\nu_2}^{m_2}(\hat{\mathbf{s}}) \right) \left( 1 - \frac{\mathcal{R}_{\hat{\mathbf{k}}_1} \mu}{\nu_1} \right) \mathcal{R}_{\hat{\mathbf{k}}_1} \Phi_{\nu_1}^{m_1}(\hat{\mathbf{s}})$$

$$= \varpi \sum_{l=0}^{l_{max}} \sum_{m=-l}^{l} \frac{\beta_l}{2l+1} Y_{lm}(\hat{\mathbf{s}}) \left(\mathcal{R}_{\hat{\mathbf{k}}_2} \Phi_{\nu_2}^{m_2}(\hat{\mathbf{s}})\right) \int_{\mathbb{S}^2} Y_{lm}^*(\hat{\mathbf{s}}') \mathcal{R}_{\hat{\mathbf{k}}_1} \Phi_{\nu_1}^{m_1}(\hat{\mathbf{s}}') d\hat{\mathbf{s}}' \,, \tag{4.40}$$

$$\left(\mathcal{R}_{\hat{\mathbf{k}}_1} \Phi_{\nu_1}^{m_1}(\hat{\mathbf{s}})\right) \left(1 - \frac{\mathcal{R}_{\hat{\mathbf{k}}_2} \mu}{\nu_2}\right) \mathcal{R}_{\hat{\mathbf{k}}_2} \Phi_{\nu_2}^{m_2}(\hat{\mathbf{s}})$$

$$= \varpi \sum_{l=0}^{l_{max}} \sum_{m=-l}^{l} \frac{\beta_l}{2l+1} Y_{lm}(\hat{\mathbf{s}}) \left(\mathcal{R}_{\hat{\mathbf{k}}_1} \Phi_{\nu_1}^{m_1}(\hat{\mathbf{s}})\right) \int_{\mathbb{S}^2} Y_{lm}^*(\hat{\mathbf{s}}') \mathcal{R}_{\hat{\mathbf{k}}_2} \Phi_{\nu_2}^{m_2}(\hat{\mathbf{s}}') d\hat{\mathbf{s}}' \,. \tag{4.41}$$

By subtraction and integration over $\hat{\mathbf{s}}$, we obtain

$$\int_{\mathbb{S}^2} \left(\frac{\mathcal{R}_{\hat{\mathbf{k}}_2} \mu}{\nu_2} - \frac{\mathcal{R}_{\hat{\mathbf{k}}_1} \mu}{\nu_1}\right) \left(\mathcal{R}_{\hat{\mathbf{k}}_1} \Phi_{\nu_1}^{m_1}(\hat{\mathbf{s}})\right) \left(\mathcal{R}_{\hat{\mathbf{k}}_2} \Phi_{\nu_2}^{m_2}(\hat{\mathbf{s}})\right) \, d\hat{\mathbf{s}}$$

$$= \left(\frac{\hat{k}_z(\nu_2 q)}{\nu_2} - \frac{\hat{k}_z(\nu_1 q)}{\nu_1}\right) \int_{\mathbb{S}^2} \mu \left(\mathcal{R}_{\hat{\mathbf{k}}_1} \Phi_{\nu_1}^{m_1}(\hat{\mathbf{s}})\right) \left(\mathcal{R}_{\hat{\mathbf{k}}_2} \Phi_{\nu_2}^{m_2}(\hat{\mathbf{s}})\right) \, d\hat{\mathbf{s}} = 0 \,. \tag{4.42}$$

This means

$$\int_{\mathbb{S}^2} \mu \left(\mathcal{R}_{\hat{\mathbf{k}}_1} \Phi_{\nu_1}^{m_1}(\hat{\mathbf{s}})\right) \left(\mathcal{R}_{\hat{\mathbf{k}}_2} \Phi_{\nu_2}^{m_2}(\hat{\mathbf{s}})\right) \, d\hat{\mathbf{s}} = 0 \quad \text{for } \nu_1 \neq \nu_2 \,. \tag{4.43}$$

Suppose that $\nu = \nu_1 = \nu_2$, $\hat{\mathbf{k}} = \hat{\mathbf{k}}_1 = \hat{\mathbf{k}}_2$, $m_1 \neq m_2$. In this case, we have

$$\int_{\mathbb{S}^2} \mu \left(\mathcal{R}_{\hat{\mathbf{k}}} \Phi_{\nu}^{m_1}(\hat{\mathbf{s}})\right) \left(\mathcal{R}_{\hat{\mathbf{k}}} \Phi_{\nu}^{m_2}(\hat{\mathbf{s}})\right) \, d\hat{\mathbf{s}} = \int_{\mathbb{S}^2} \left(\mathcal{R}_{\hat{\mathbf{k}}}^{-1} \mu\right) \Phi_{\nu}^{m_1}(\hat{\mathbf{s}}) \Phi_{\nu}^{m_2}(\hat{\mathbf{s}}) d\hat{\mathbf{s}}$$

$$= \hat{k}_z(\nu q) \int_{\mathbb{S}^2} \mu \Phi_{\nu}^{m_1}(\hat{\mathbf{s}}) \Phi_{\nu}^{m_2}(\hat{\mathbf{s}}) d\hat{\mathbf{s}}$$

$$= 2\pi \hat{k}_z(\nu q) \delta_{m_1,-m_2} \int_{-1}^{1} \mu \phi^{m_1}(\nu,\mu) \phi^{-m_1}(\nu,\mu) (1-\mu^2)^{|m_1|} \, d\mu \,. \tag{4.44}$$

We note that ((4.14) and Appendix A)

$$\phi^{-m}(\nu,\mu) = \phi^{m}(\nu,\mu), \quad \Phi_{\nu}^{-m}(\hat{\mathbf{s}}) = \Phi_{\nu}^{m*}(\hat{\mathbf{s}}). \tag{4.45}$$

Thus we have shown (4.35).

### 4.4. Separation constant

From the one-dimensional transport theory, the separation constant $\nu$ is known to be real, and an eigenvalue ($\nu \notin [-1,1]$) or in the continuous spectrum ($\nu \in (-1,1)$). If $\nu$ is an eigenvalue, then $-\nu$ is also an eigenvalue. The number $2M^m$ of eigenvalues ($-l_{\max} \leq m \leq l_{\max}$) is determined by $\varpi$ and $p(\hat{\mathbf{s}}, \hat{\mathbf{s}}')$. For example, when $p(\hat{\mathbf{s}}, \hat{\mathbf{s}}') = 1/(4\pi)$ for isotropic scattering ($l_{\max} = 0$), $M^0 = 1$ and there are two eigenvalues $\pm\nu_0$ ($\nu_0 > 1$), which are the solutions to (Case and Zweifel, 1967)

$$1 - \varpi \nu \tanh^{-1} \frac{1}{\nu} = 0 \,. \tag{4.46}$$

In general, an eigenvalue is a solution to

$$\Lambda^m(\nu) = 0 \,, \tag{4.47}$$

where

$$\Lambda^m(\nu) = 1 - \frac{\varpi \nu}{2} \int_{-1}^{1} \frac{g^m(\nu,\mu)}{\nu - \mu} (1-\mu^2)^{|m|} \, d\mu \,, \quad \nu \notin [-1,1] \,. \tag{4.48}$$

We will denote eigenvalues as

$$\nu = \pm\nu_j^m \,, \;\; \nu_j^m > 1 \,, \;\; j = 0,1,\cdots, M^m - 1 \,. \tag{4.49}$$

We can numerically evaluate eigenvalues $\pm\nu_j^m$ as eigenvalues of a tridiagonal matrix $A(m)$ $(-l_{\max} \le m \le l_{\max})$ which is given by

$$A(m) = \begin{pmatrix} 0 & a_{|m|+1} & 0 & & & \\ a_{|m|+1} & 0 & a_{|m|+2} & & & \\ 0 & a_{|m|+2} & 0 & \ddots & & \\ & & \ddots & \ddots & & \\ & & & & 0 & a_{l_B} \\ & & & & a_{l_B} & 0 \end{pmatrix} , \tag{4.50}$$

where $l_B$ $(\ge l_{\max})$ is sufficiently large and

$$a_l = a_l(m) = \sqrt{\frac{l^2 - m^2}{h_l h_{l-1}}} \,. \tag{4.51}$$

The matrix $A(m)$ has $(l_B - |m| + 1)/2$ or $(l_B - |m|)/2$ positive eigenvalues depending on whether $l_B - |m| + 1$ is even or odd. To see how $A(m)$ is obtained, we first show that eigenvalues are zeros of $g_l^m$ as $l \to \infty$ (Garcia and Siewert, 1982).

Let us define

$$q_l^m(w) = \frac{1}{2}\int_{-1}^{1} \frac{p_l^m(\mu)}{w - \mu}(1-\mu^2)^{|m|}\, d\mu \;, \;\; w \notin [-1,1] \,. \tag{4.52}$$

For $\nu \notin [-1,1]$, the three-term recurrence relation of $p_l^m$ implies

$$\sqrt{(l+1)^2 - m^2}\, q_{l+1}^m(\nu) = (2l+1)\nu q_l^m(\nu) - \sqrt{l^2 - m^2}\, q_{l-1}^m(\nu) \\ -\big(\mathrm{sgn}(m)\big)^m \frac{\sqrt{(2|m|)!}}{(2|m|-1)!!}\delta_{l|m|} \,, \tag{4.53}$$

where $\mathrm{sgn}(m) = 1$ for $m \ge 0$ and $\mathrm{sgn}(m) = -1$ for $m < 0$. We multiply $q_l^m(\nu)$ on both sides of (4.19) and multiply $g_l^m(\nu)$ on both sides of (4.53). By subtracting the former from the latter, we obtain

$$\sqrt{(l+1)^2 - m^2}\big(q_{l+1}^m(\nu) g_l^m(\nu) - q_l^m(\nu) g_{l+1}^m(\nu)\big) = (2l + 1 - h_l)\nu q_l^m(\nu) g_l^m(\nu) \\ -\sqrt{l^2 - m^2}\big(q_{l-1}^m(\nu) g_l^m(\nu) - q_l^m(\nu) g_{l-1}^m(\nu)\big) - \delta_{l|m|} \,. \tag{4.54}$$

By taking the sum over $l$ from $|m|$ to $l_B$, we obtain

$$\sqrt{(l_B+1)^2 - m^2}\big(q_{l_B+1}^m(\nu) g_{l_B}^m(\nu) - q_{l_B}^m(\nu) g_{l_B+1}^m(\nu)\big) \\ = \sum_{l=|m|}^{l_B} (2l + 1 - h_l)\nu q_l^m(\nu) g_l^m(\nu) - 1 \,. \tag{4.55}$$

Noting that $\Lambda^m(\nu) = 1 - \varpi\nu \sum_{l=|m|}^{l_{\max}} \beta_l g_l^m(\nu) q_l^m(\nu)$, we obtain (the Christoffel-Darboux formula)

$$\Lambda^m(\nu) = \sqrt{(l_B+1)^2 - m^2}\left(q_{l_B}^m(\nu) g_{l_B+1}^m(\nu) - q_{l_B+1}^m(\nu) g_{l_B}^m(\nu)\right) . \tag{4.56}$$

Next we multiply $p_l^m(\nu)$ on both sides of (4.53) and multiply $q_l^m(\nu)$ on both sides of (4.28). By subtracting the former from the latter and by summing the resulting expression over $l$ from $|m|$ to $l_B$, we obtain

$$1 = \sqrt{(l_B+1)^2 - m^2}\left(p_{l_B+1}^m(\mu) q_{l_B}^m(\nu) - p_{l_B}^m(\nu) q_{l_B+1}^m(\nu)\right) . \tag{4.57}$$

Similarly, we multiply $p_l^m(\nu)$ on both sides of (4.19) and $g_l^m(\nu)$ on both sides of (4.28), subtract the former from the latter, and take the sum over $l$ from $|m|$ to $l_B$. We obtain

$$\varpi \nu g^m(\nu,\nu) = \sqrt{(l_B+1)^2 - m^2}\left(p_{l_B+1}^m(\mu) g_{l_B}^m(\nu) \ - p_{l_B}^m(\nu) g_{l_B+1}^m(\nu)\right) . \tag{4.58}$$

We obtain using (4.56), (4.57), and (4.58),

$$p_{l_B+1}^m(\nu)\Lambda^m(\nu) = g_{l_B+1}^m(\nu) - \varpi \nu g^m(\nu,\nu) q_{l_B+1}^m(\nu) . \tag{4.59}$$

We note that

$$\lim_{l\to\infty} \frac{q_l^m(w)}{p_l^m(w)} = \lim_{l\to\infty} \frac{Q_l^m(w)}{P_l^m(w)} = 0 , \quad w \notin [-1,1] ,$$

where $Q_l^m$ is the associated Legendre polynomial of the second kind. Therefore we obtain (Garcia and Siewert, 1982)

$$\Lambda^m(\nu) = \lim_{l_B\to\infty} \frac{g_{l_B+1}^m(\nu)}{p_{l_B+1}^m(\nu)} . \tag{4.60}$$

Hence, eigenvalues are zeros of $g_l^m$ as $l \to \infty$.

Let us rewrite (4.19) as

$$a_l(m)\sqrt{h_{l-1}} g_{l-1}^m(\nu) + a_{l+1}(m)\sqrt{h_{l+1}} g_{l+1}^m(\nu) = \nu\sqrt{h_l} g_l^m(\nu) . \tag{4.61}$$

The recurrence relation (4.61) can be further rewritten in a matrix-vector form if $g_{l_B+1}^m = 0$. That is, if $\nu$ is an eigenvalue of the matrix $A(m)$, it means $g_{l_B+1}^m(\nu) = 0$, i.e., $\nu$ is a zero of $g_{l_B+1}^m$. Hence, $\nu_j^m$ can be numerically obtained as eigenvalues of $A(m)$ for sufficiently large $l_B$. In Garcia and Siewert, 1989, numerical techniques to accurately compute $\nu_j^m$ were developed.

### 4.5. Fundamental solution

Let us consider the fundamental solution $G(\mathbf{r},\hat{\mathbf{s}};\mathbf{r}_0,\hat{\mathbf{s}}_0)$ or the Green's function for the whole space, which satisfies

$$\hat{\mathbf{s}}\cdot\nabla G(\mathbf{r},\hat{\mathbf{s}};\mathbf{r}_0,\hat{\mathbf{s}}_0) + \mu_t G(\mathbf{r},\hat{\mathbf{s}};\mathbf{r}_0,\hat{\mathbf{s}}_0) = \mu_s \int_{\mathbb{S}^2} p(\hat{\mathbf{s}},\hat{\mathbf{s}}') G(\mathbf{r},\hat{\mathbf{s}}';\mathbf{r}_0,\hat{\mathbf{s}}_0) d\hat{\mathbf{s}}' + \delta(\mathbf{r}-\mathbf{r}_0)\delta(\hat{\mathbf{s}}-\hat{\mathbf{s}}_0) \tag{4.62}$$

for $\mathbf{r},\mathbf{r}_0 \in \mathbb{R}^3$ and $\hat{\mathbf{s}},\hat{\mathbf{s}}_0 \in \mathbb{S}^2$. Note that $G(\mathbf{r},\hat{\mathbf{s}};\mathbf{r}_0,\hat{\mathbf{s}}_0) \to 0$ as $|\mathbf{r}| \to \infty$. Here, multi-dimensional Dirac's delta functions were introduced as

$$\delta(\mathbf{r}-\mathbf{r}_0) = \delta(x-x_0)\delta(y-y_0)\delta(z-z_0) , \quad \delta(\hat{\mathbf{s}}-\hat{\mathbf{s}}_0) = \delta(\mu-\mu_0)\delta(\varphi-\varphi_0) . \tag{4.63}$$

We express $G(\mathbf{r},\hat{\mathbf{s}};\mathbf{r}_0,\hat{\mathbf{s}}_0)$ as

$$G(\mathbf{r},\hat{\mathbf{s}};\mathbf{r}_0,\hat{\mathbf{s}}_0) = \sum_{m=-l_{\max}}^{l_{\max}} \times \begin{cases} \int_{\mathbb{R}^2}\int_+ A^m(\nu,\mathbf{q}) e^{i\mu_t \mathbf{q}\cdot\boldsymbol{\rho}} e^{-\mu_t \hat{k}_z(\nu q)(z-z_0)/\nu} \mathcal{R}_{\hat{\mathbf{k}}(\nu,\mathbf{q})} \Phi_\nu^m(\hat{\mathbf{s}}) d\nu\, d\mathbf{q} , & z > z_0 , \\ \int_{\mathbb{R}^2}\int_+ B^m(\nu,\mathbf{q}) e^{i\mu_t \mathbf{q}\cdot\boldsymbol{\rho}} e^{\mu_t \hat{k}_z(\nu q)(z-z_0)/\nu} \mathcal{R}_{\hat{\mathbf{k}}(-\nu,\mathbf{q})} \Phi_{-\nu}^m(\hat{\mathbf{s}}) d\nu\, d\mathbf{q} , & z < z_0 . \end{cases} \tag{4.64}$$

Here, for a function $f^m(\nu)$, the notation $\int_+ f^m(\nu)\, d\nu$ means,

$$\int_+ f^m(\nu)\, d\nu = \sum_{j=0}^{M^m-1} f(\nu_j^m) + \int_0^1 f^m(\nu)\, d\nu . \tag{4.65}$$

In (4.64), $A^m(\nu,\mathbf{q}), B^m(\nu,\mathbf{q})$ will be determined so that $G(\mathbf{r},\hat{\mathbf{s}};\mathbf{r}_0,\hat{\mathbf{s}}_0)$ satisfies the RTE.

Let us consider the Fourier transform:

$$\begin{cases} \tilde{G}(\mathbf{q};z,\hat{\mathbf{s}};\mathbf{r}_0,\hat{\mathbf{s}}_0) = \int_{\mathbb{R}^2} e^{-i\mu_t\mathbf{q}\cdot\boldsymbol{\rho}} G(\mathbf{r},\hat{\mathbf{s}};\mathbf{r}_0,\hat{\mathbf{s}}_0) d\boldsymbol{\rho} \ , \\ G(\mathbf{r},\hat{\mathbf{s}};\mathbf{r}_0,\hat{\mathbf{s}}_0) = \frac{{\mu_t}^2}{(2\pi)^2}\int_{\mathbb{R}^2} e^{i\mu_t\mathbf{q}\cdot\boldsymbol{\rho}} \tilde{G}(\mathbf{q};\ z,\hat{\mathbf{s}};\mathbf{r}_0,\hat{\mathbf{s}}_0)\, d\mathbf{q} \ . \end{cases} \tag{4.66}$$

Here the Fourier transform $\tilde{G}(\mathbf{q};z,\hat{\mathbf{s}};\mathbf{r}_0,\hat{\mathbf{s}}_0)$ satisfies

$$\mu\frac{\partial}{\partial z}\tilde{G}(\mathbf{q};z,\hat{\mathbf{s}};\mathbf{r}_0,\hat{\mathbf{s}}_0) + (1+i\boldsymbol{\omega}\cdot\boldsymbol{q})\mu_t\tilde{G}(\mathbf{q};z,\hat{\mathbf{s}};\mathbf{r}_0,\hat{\mathbf{s}}_0)$$

$$= \mu_s\int_{\mathbb{S}^2} p(\hat{\mathbf{s}},\hat{\mathbf{s}}')\tilde{G}(\mathbf{q};z,\hat{\mathbf{s}}';\mathbf{r}_0,\hat{\mathbf{s}}_0)d\hat{\mathbf{s}}' + e^{-i\mu_t\mathbf{q}\cdot\boldsymbol{\rho}_0}\delta(z-z_0)\delta(\hat{\mathbf{s}}-\hat{\mathbf{s}}_0) \ . \tag{4.67}$$

If we integrate (4.67) from $z=z_0-0$ to $z=z_0+0$, we obtain the jump condition:

$$\mu\big(\tilde{G}(\mathbf{q};z_0+0,\hat{\mathbf{s}};\mathbf{r}_0,\hat{\mathbf{s}}_0) - \tilde{G}(\mathbf{q};z_0-0,\hat{\mathbf{s}};\mathbf{r}_0,\hat{\mathbf{s}}_0)\big) = e^{-i\mu_t\mathbf{q}\cdot\boldsymbol{\rho}_0}\delta(\hat{\mathbf{s}}-\hat{\mathbf{s}}_0) \ . \tag{4.68}$$

Using the orthogonality relation (4.35), we obtain

$$A^m(\nu,\mathbf{q}) = \frac{{\mu_t}^2 e^{-i\mu_t\mathbf{q}\cdot\boldsymbol{\rho}_0}}{2\pi\hat{k}_z(\nu q)\mathcal{N}^m(\nu)}\mathcal{R}_{\hat{\mathbf{k}}(\nu,\mathbf{q})}\Phi_\nu^{m*}(\hat{\mathbf{s}}_0), \tag{4.69}$$

$$B^m(\nu,\mathbf{q}) = \frac{{\mu_t}^2 e^{-i\mu_t\mathbf{q}\cdot\boldsymbol{\rho}_0}}{2\pi\hat{k}_z(\nu q)\mathcal{N}^m(\nu)}\mathcal{R}_{\hat{\mathbf{k}}(-\nu,\mathbf{q})}\Phi_{-\nu}^{m*}(\hat{\mathbf{s}}_0). \tag{4.70}$$

Hence,

$$\tilde{G}(\mathbf{q};\ z,\hat{\mathbf{s}};\mathbf{r}_0,\hat{\mathbf{s}}_0) = e^{-i\mu_t\mathbf{q}\cdot\boldsymbol{\rho}_0}\sum_{m=-l_{\max}}^{l_{\max}}\int_+ \frac{e^{-\mu_t\hat{k}_z(\nu q)|z-z_0|/\nu}}{2\pi\hat{k}_z(\nu q)\mathcal{N}^m(\nu)}\mathcal{R}_{\hat{\mathbf{k}}(\pm\nu,\mathbf{q})}\Phi_{\pm\nu}^m(\hat{\mathbf{s}})\Phi_{\pm\nu}^{m*}(\hat{\mathbf{s}}_0)d\nu \ . \tag{4.71}$$

Here, upper signs are chosen for $z>z_0$ and lower signs are chosen for $z<z_0$.

Since singular eigenfunctions $\phi^m(\nu,\mu)$ contain generalized functions and direct computation of $G(\mathbf{r},\hat{\mathbf{s}};\mathbf{r}_0,\hat{\mathbf{s}}_0)$ is complicated, we will rather try to find numerical solutions to the RTE making use of rotated reference frames.

From the fundamental solution in (4.66), (4.71), we can obtain the energy density $U(\mathbf{r})$ of the specific intensity $I(\mathbf{r},\hat{\mathbf{s}})$ which satisfies the following RTE.

$$\hat{\mathbf{s}}\cdot\nabla I(\mathbf{r},\hat{\mathbf{s}}) + \mu_t I(\mathbf{r},\hat{\mathbf{s}}) = \mu_s\int_{\mathbb{S}^2} p(\hat{\mathbf{s}},\hat{\mathbf{s}}')I(\mathbf{r},\hat{\mathbf{s}}')\, d\hat{\mathbf{s}}' + \delta(\mathbf{r}-\mathbf{r}_0) \ . \tag{4.72}$$

The energy density is given by

$$U(\mathbf{r}) = \frac{1}{c}\int_{\mathbb{S}^2} I(\mathbf{r},\hat{\mathbf{s}})\, d\hat{\mathbf{s}} \ . \tag{4.73}$$

We can calculate $U(\mathbf{r})$ as

$$U(\mathbf{r}) = \frac{1}{c}\int_{\mathbb{S}^2}\int_{\mathbb{S}^2} G(\mathbf{r},\hat{\mathbf{s}};\mathbf{r}_0,\hat{\mathbf{s}}_0)\, d\hat{\mathbf{s}}d\hat{\mathbf{s}}_0 = \frac{{\mu_t}^2}{2\pi c}\int_{\mathbb{R}^2} e^{i\mu_t\mathbf{q}\cdot(\boldsymbol{\rho}-\boldsymbol{\rho}_0)}\int_+ \frac{e^{-\mu_t\hat{k}_z(\nu q)|z-z_0|/\nu}}{\hat{k}_z(\nu q)\mathcal{N}^0(\nu)} d\nu\, d\mathbf{q} \ . \tag{4.74}$$

By noticing the Hansen-Bessel formula,

$$J_0(x) = \frac{1}{2\pi}\int_0^{2\pi} e^{ix\cos\varphi}\, d\varphi \ , \ \ x \geq 0 \ , \tag{4.75}$$

where $J_0$ is the Bessel function of order 0, we obtain

$$U(\mathbf{r}) = \frac{{\mu_t}^2}{c}\int_0^\infty qJ_0(\mu_t q|\boldsymbol{\rho}-\boldsymbol{\rho}_0|)\int_+ \frac{e^{-\mu_t\hat{k}_z(\nu q)|z-z_0|/\nu}}{\hat{k}_z(\nu q)\mathcal{N}^0(\nu)} d\nu\, dq \ . \tag{4.76}$$

Since the source is isotropic, without loss of generality we can set

$$\mathbf{r} = \begin{pmatrix} \boldsymbol{\rho}_0 \\ z \end{pmatrix} . \tag{4.77}$$

By setting $Q = \hat{k}_z(\nu q)$,

$$U(z) = U(\mathbf{r}) = \frac{{\mu_t}^2}{c} \int_{+} \int_1^\infty \frac{e^{-\mu_t Q|z-z_0|/\nu}}{\nu^2 \mathcal{N}^0(\nu)} \, dQ \, d\nu \ . \tag{4.78}$$

Thus we obtain the following expression, which holds for general anisotropic scattering (Avram and Machida, 2013).

$$U(z) = \frac{\mu_t}{c|z-z_0|} \left( \sum_{j=0}^{M^0-1} \frac{e^{-\mu_t|z-z_0|/\nu_j^0}}{\nu_j^0 \mathcal{N}^0(\nu_j^0)} + \int_0^1 \frac{e^{-\mu_t|z-z_0|/\nu}}{\nu \mathcal{N}^0(\nu)} \, d\nu \right) . \tag{4.79}$$

In the case of the linear scattering:

$$p(\hat{\mathbf{s}}, \hat{\mathbf{s}}') = \frac{1}{4\pi} + \frac{3\mathrm{g}}{4\pi} \hat{\mathbf{s}} \cdot \hat{\mathbf{s}}' , \tag{4.80}$$

the positive eigenvalue $\nu_0^0$ can be obtained from the matrix $A(0)$ in (4.50) (Zhang, Avram, and Machida, 2022). In the case of isotropic scattering, the formula (4.79) was obtained by Ganapol and Kornreich, 1995.

We admit that the operator $\mathcal{R}_{\hat{\mathbf{k}}(\nu,\mathbf{q})}$ is not yet fully explained. Before exploring $\mathcal{R}_{\hat{\mathbf{k}}(\nu,\mathbf{q})}$ in Sec. 6, let us consider the fundamental solution in the two-dimensional space (flatland).

# 5. Fundamental solution in flatland

In transport theory, the one-dimensional RTE means the RTE with one spatial variable and two angular variables, whereas the three-dimensional RTE has three spatial variables and two angular variables. In this sense, the two-dimensional RTE is an equation with two spatial variables and two angular variables. Here, we consider light propagation in two-dimensional space. This RTE with two spatial variables and one angular variable is often called the RTE in flatland, which was named after a novel by Edwin A. Abbott (Abbott, 1884).

Let us consider the following RTE in flatland:

$$\left[ \begin{pmatrix} \cos\varphi \\ \sin\varphi \end{pmatrix} \cdot \begin{pmatrix} \partial/\partial x \\ \partial/\partial y \end{pmatrix} + \mu_t \right] G(\boldsymbol{\rho}, \varphi; \varphi_0) = \mu_s \int_0^{2\pi} p(\varphi, \varphi') G(\boldsymbol{\rho}, \varphi'; \varphi_0) \, d\varphi' + \delta(\boldsymbol{\rho}) \delta(\varphi - \varphi_0) , \tag{5.1}$$

where $G(\boldsymbol{\rho}, \varphi; \varphi_0) \to 0$ as $|\boldsymbol{\rho}| \to \infty$. We give the scattering phase function as

$$p(\varphi, \varphi') = \frac{1}{2\pi} \sum_{m=-l_{\max}}^{l_{\max}} \beta_m^{(2D)} e^{im(\varphi-\varphi')} = \frac{1}{2\pi} + \frac{1}{\pi} \sum_{m=-l_{\max}}^{l_{\max}} \beta_m^{(2D)} \cos[m(\varphi - \varphi')] \ . \tag{5.2}$$

Here, $\beta_0^{(2D)} = 1$, $\beta_m^{(2D)} \in [-1,1]$, and $\beta_{-m}^{(2D)} = \beta_m^{(2D)}$. We obtain (Machida, 2016b)

$$G(\boldsymbol{\rho}, \varphi; \varphi_0) = \frac{\mu_t}{2\pi} \int_{-\infty}^{\infty} e^{i\mu_t q y} \left[ \sum_{j=0}^{M^{(2D)}-1} \phi^{(2D)} \left( \pm\nu_j^{(2D)}, \varphi - \varphi_{\hat{\mathbf{k}}\left(\pm\nu_j^{(2D)} q\right)} \right) \right.$$

$$\times \phi^{(2D)} \left( \pm\nu_j^{(2D)}, \varphi_0 - \varphi_{\hat{\mathbf{k}}\left(\pm\nu_j^{(2D)} q\right)} \right) \frac{e^{-\sqrt{1+\left(\nu_j^{(2D)} q\right)^2}|x|/\nu_j^{(2D)}}}{\sqrt{1 + \left(\nu_j^{(2D)} q\right)^2} \mathcal{N}^{(2D)}\left(\nu_j^{(2D)}\right)}$$

$$+\int_0^1 \phi^{(2D)}\left(\pm\nu_j^{(2D)}, \varphi - \varphi_{\hat{\mathbf{k}}\left(\pm\nu_j^{(2D)}q\right)}\right)\phi^{(2D)}\left(\pm\nu_j^{(2D)}, \varphi_0 - \varphi_{\hat{\mathbf{k}}\left(\pm\nu_j^{(2D)}q\right)}\right)$$
$$\left.\times \frac{e^{-\sqrt{1+\left(\nu^{(2D)}q\right)^2}|x|/\nu^{(2D)}}}{\sqrt{1+(\nu^{(2D)}q)^2}\mathcal{N}^{(2D)}(\nu^{(2D)})}\, d\nu\right] dq\,, \tag{5.3}$$

where upper signs are used for $x > 0$ and lower signs are used for $x < 0$. The notations are described below.

The separation constant $\nu^{(2D)}$ is continuous in $(-1,1)$ and it can be written for $\nu^{(2D)} \neq [-1,1]$ as $\nu^{(2D)} = \pm\nu_j^{(2D)}$ ($j = 0, \cdots, M^{(2D)} - 1$). Here, $\nu_j^{(2D)}$ are positive roots of $\Lambda^{(2D)}\left(\nu_j^{(2D)}\right) = 0$, where

$$\Lambda^{(2D)}(z) = 1 - \frac{\varpi z}{2\pi}\int_0^{2\pi} \frac{g^{(2D)}(z,\varphi)}{z - \cos\varphi}\, d\varphi\ , \quad z \in \mathbb{C} \setminus [-1,1]\,. \tag{5.4}$$

We note that $M^{(2D)}$ is the number of positive roots and the function $g^{(2D)}(z,\varphi)$ is given by

$$g^{(2D)}(z,\varphi) = 1 + 2\sum_{m=1}^{l_{\max}} \beta_m^{(2D)}\gamma_m(z)\cos m\varphi\ , \tag{5.5}$$

where the polynomials $\gamma_m(z)$ satisfy the following three-term recurrence relation (Machida, 2016b):

$$2zh_m^{(2D)}\gamma_m(z) - \gamma_{m+1}(z) - \gamma_{m-1}(z) = 0 \tag{5.6}$$

with initial terms,

$$\gamma_0(z) = 1\,,\ \ \gamma_1(z) = (1-\varpi)z\,, \tag{5.7}$$

and

$$h_m^{(2D)} = 1 - \varpi\beta_m^{(2D)}\,. \tag{5.8}$$

In flatland, singular eigenfunctions are given by (Machida, 2016b)

$$\phi^{(2D)}\left(\nu^{(2D)},\varphi\right) = \frac{\varpi\nu^{(2D)}}{2\pi}\mathcal{P}\frac{g^{(2D)}\left(\nu^{(2D)},\varphi\right)}{\nu^{(2D)} - \cos\varphi}$$
$$+\frac{\sqrt{1-\nu^{(2D)2}}}{2}\lambda^{(2D)}\left(\nu^{(2D)}\right)\delta\left(\nu^{(2D)} - \cos\varphi\right), \tag{5.9}$$

where $\mathcal{P}$ denotes Cauchy's principal value and

$$\lambda^{(2D)}\left(\nu^{(2D)}\right) = 1 - \frac{\varpi\nu^{(2D)}}{2\pi}\mathcal{P}\int_0^{2\pi}\frac{g^{(2D)}\left(\nu^{(2D)},\varphi\right)}{\nu^{(2D)} - \cos\varphi}\, d\varphi\ , \quad \nu^{(2D)} \in (-1,1)\,. \tag{5.10}$$

Singular eigenfunctions satisfy the orthogonality relation (Machida, 2016b):

$$\int_0^{2\pi} \phi^{(2D)}(\nu,\varphi)\phi^{(2D)}(\nu',\varphi)\cos\varphi\ d\varphi = \mathcal{N}^{(2D)}(\nu)\delta(\nu - \nu')\,. \tag{5.11}$$

The normalization factor is given by

$$\mathcal{N}^{(2D)}(\nu) = \begin{cases} \dfrac{2\nu}{\sqrt{1-\nu^2}}\left[\left(\dfrac{\varpi\nu}{2}\right)^2 g^{(2D)}(\nu,\varphi_\nu)^2 + \lambda^{(2D)}(\nu)^2\right]\,, & \nu \in (-1,1)\,, \\ \left(\dfrac{\varpi\nu}{2}\right)^2 g^{(2D)}(\nu,\varphi_\nu)\dfrac{d\Lambda^{(2D)}(\nu)}{d\nu}\,, & \nu \notin [-1,1]\,, \end{cases} \tag{5.12}$$

where

$$\varphi_\nu = \begin{cases} \cos^{-1}\nu\ , & \nu \in [-1,1]\ , \\ i\cosh^{-1}\nu\ , & \nu > 1\ , \\ \pi + i\cosh^{-1}|\nu|\ , & \nu < -1\ . \end{cases} \tag{5.13}$$

We note that $0 \le \cos^{-1}\nu \le \pi$ for $\nu \in [-1,1]$ and $\cosh^{-1}|\nu| = \ln\left(|\nu| + \sqrt{\nu^2 - 1}\right)$ for $|\nu| > 1$. In singular eigenfunctions in the expression of the fundamental solution (5.3), $\varphi_{\hat{\mathbf{k}}\left(\pm\nu_j^{(2D)}q\right)}$ appears in the forms of $\cos\varphi_{\hat{\mathbf{k}}\left(\pm\nu_j^{(2D)}q\right)}$ and $\sin\varphi_{\hat{\mathbf{k}}\left(\pm\nu_j^{(2D)}q\right)}$. We obtain

$$\cos\varphi_{\hat{\mathbf{k}}\left(\pm\nu_j^{(2D)}q\right)} = \sqrt{1 + \left(\nu_j^{(2D)}q\right)^2}\ ,\quad \sin\varphi_{\hat{\mathbf{k}}\left(\pm\nu_j^{(2D)}q\right)} = \mp i\nu_j^{(2D)}q\ . \tag{5.14}$$

It should be mentioned that for isotropic scattering ($l_{\max} = 0$) in flatland, the positive eigenvalue is obtained as (Machida, 2016b)

$$\nu_0^{(2D)} = \frac{1}{\sqrt{1 - \varpi^2}}\ , \tag{5.15}$$

whereas the one-dimensional RTE (three spatial dimensions with planar symmetry) with isotropic scattering has the largest eigenvalue $\nu_0$ which is a solution of the transcendental equation (4.46) but cannot be explicitly written down.

Let us calculate the energy density $u^{(2D)}$ (up to a constant) for an isotropic source $\delta(\boldsymbol{\rho})$:

$$u^{(2D)} = \int_0^{2\pi}\int_0^{2\pi} G(\boldsymbol{\rho}, \varphi; \varphi_0)\, d\varphi d\varphi_0\ . \tag{5.16}$$

We can set $y = 0$ without loss of generality. Assume $x > 0$. We obtain

$$u^{(2D)} = \frac{\mu_t}{\pi}\left[\sum_{j=0}^{M^{(2D)}-1} \frac{K_0\left(x/\nu_j^{(2D)}\right)}{\nu_j^{(2D)}\mathcal{N}^{(2D)}\left(\nu_j^{(2D)}\right)} + \int_0^1 \frac{K_0(x/\nu)}{\nu\mathcal{N}^{(2D)}(\nu)}\, d\nu\right]\ , \tag{5.17}$$

where $K_0$ is the modified Bessel function of the second kind of order zero. The formula (5.17) is the general result. In particular, for isotropic scattering ($l_{\max} = 0$), we have

$$u^{(2D)} = \frac{\mu_t}{\pi}\left[\frac{K_0\left(x/\nu_0^{(2D)}\right)}{\nu_0^{(2D)}\mathcal{N}^{(2D)}\left(\nu_0^{(2D)}\right)} + \int_0^1 \frac{K_0(x/\nu)}{\nu\mathcal{N}^{(2D)}(\nu)}\, d\nu\right]\ , \tag{5.18}$$

where $\nu_0^{(2D)}$ is given in (5.15). In the case of isotropic scattering, it is also possible to obtain $u^{(2D)}$ with the Fourier transform. It is expected that the relation between (5.18) and the expression via the Fourier transform can be shown similar to the three-dimensional case (Machida, 2016a). See also Ganapol, 2000 and Ganapol, 2015. In Fig. 5.1, $u^{(2D)}$ in (5.18) is compared to the result from finite-difference method (Fujii, 2017).

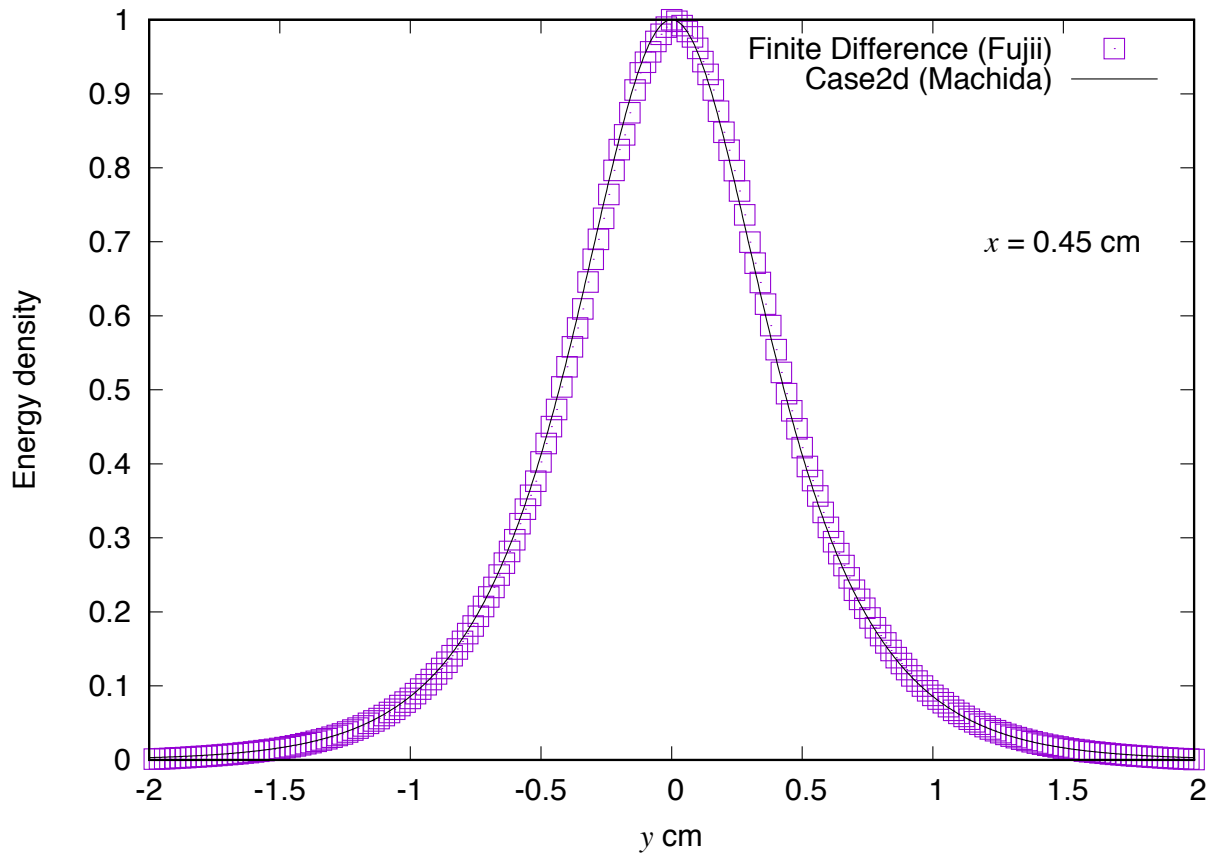

Figure 5.1. The energy density $u^{(2D)}$ in (5.18) is plotted with the result by finite difference method (purple open square) (Fujii, 2017).

# 6. Rotation of reference frames

Let us assume that function $f(\hat{\mathbf{s}}) = f(\theta, \varphi)$ can be expanded by spherical harmonics:

$$f(\hat{\mathbf{s}}) = \sum_{l=0}^{\infty} \sum_{m=-l}^{l} f_{lm} Y_{lm}(\hat{\mathbf{s}}) \ , \tag{6.1}$$

where $f_{lm} \in \mathbb{C}$ are coefficients given by

$$f_{lm} = \int_{\mathbb{S}^2} f(\hat{\mathbf{s}}) Y_{lm}^*(\hat{\mathbf{s}}) d\hat{\mathbf{s}} \ . \tag{6.2}$$

For a unit vector $\hat{\mathbf{k}} \in \mathbb{C}^3$, $\hat{\mathbf{k}} \cdot \hat{\mathbf{k}} = 1$, the operator $\mathcal{R}_{\hat{\mathbf{k}}}$ is defined as

$$\begin{aligned} \mathcal{R}_{\hat{\mathbf{k}}} f(\hat{\mathbf{s}}) &= \sum_{l=0}^{\infty} \sum_{m=-l}^{l} f_{lm} \mathcal{R}_{\hat{\mathbf{k}}} Y_{lm}(\hat{\mathbf{s}}) = \sum_{l=0}^{\infty} \sum_{m=-l}^{l} f_{lm} \sum_{m'=-l}^{l} D_{m'm}^{l}(\varphi_{\hat{\mathbf{k}}}, \theta_{\hat{\mathbf{k}}}, 0) Y_{lm'}(\hat{\mathbf{s}}) \\ &= \sum_{l=0}^{\infty} \sum_{m=-l}^{l} f_{lm} \sum_{m'=-l}^{l} e^{-im'\varphi_{\hat{\mathbf{k}}}} d_{m'm}^{l}(\theta_{\hat{\mathbf{k}}}) Y_{lm'}(\hat{\mathbf{s}}) \ , \end{aligned} \tag{6.3}$$

where $\theta_{\hat{\mathbf{k}}}, \varphi_{\hat{\mathbf{k}}}$ are the polar and azimuthal angles of $\hat{\mathbf{k}}$, respectively. Here, $D_{m'm}^{l}, d_{m'm}^{l}$ are Wigner's *D*- and *d*-matrices, respectively (Varshalovich, Moskalev, and Khersonskii, 1988). The rotation of the reference frame so that the $z$-axis coincides with the direction of $\hat{\mathbf{k}}$ is not unique. In (6.3), we defined the rotation by setting three Euler angles of $D_{m'm}^{l}(\alpha, \beta, \gamma)$ as $\alpha = \varphi_{\hat{\mathbf{k}}}$, $\beta = \theta_{\hat{\mathbf{k}}}$, and $\gamma = 0$ (Markel, 2004).

By direct calculation we can show (Appendix B)

$$\left(\mathcal{R}_{\hat{\mathbf{k}}} Y_{l_1 m_1}(\hat{\mathbf{s}})\right)\left(\mathcal{R}_{\hat{\mathbf{k}}} Y_{l_2 m_2}(\hat{\mathbf{s}})\right) = \mathcal{R}_{\hat{\mathbf{k}}}\left(Y_{l_1 m_1}(\hat{\mathbf{s}}) Y_{l_2 m_2}(\hat{\mathbf{s}})\right) \ . \tag{6.4}$$

This implies the homomorphism:

$$\left(\mathcal{R}_{\hat{\mathbf{k}}} f(\hat{\mathbf{s}})\right)\left(\mathcal{R}_{\hat{\mathbf{k}}} g(\hat{\mathbf{s}})\right) = \mathcal{R}_{\hat{\mathbf{k}}}\left(f(\hat{\mathbf{s}}) g(\hat{\mathbf{s}})\right) . \tag{6.5}$$

The inverse $\mathcal{R}_{\hat{\mathbf{k}}}^{-1}$ is given as follows.

$$\mathcal{R}_{\hat{\mathbf{k}}}^{-1}Y_{lm}(\hat{\mathbf{s}}) = \sum_{m'=-l}^{l} D_{m'm}^{l}(0,-\theta_{\hat{\mathbf{k}}},-\varphi_{\hat{\mathbf{k}}})Y_{lm'}(\hat{\mathbf{s}}) = \sum_{m'=-l}^{l} e^{im\varphi_{\hat{\mathbf{k}}}}d_{mm'}^{l}(\theta_{\hat{\mathbf{k}}})Y_{lm'}(\hat{\mathbf{s}}) \,. \qquad (6.6)$$

Indeed, we have (Machida, 2015)

$$\mathcal{R}_{\hat{\mathbf{k}}}^{-1}\mathcal{R}_{\hat{\mathbf{k}}}Y_{lm}(\hat{\mathbf{s}}) = \sum_{m'=-l}^{l} e^{-im'\varphi_{\hat{\mathbf{k}}}}d_{m'm}^{l}(\theta_{\hat{\mathbf{k}}})\mathcal{R}_{\hat{\mathbf{k}}}^{-1}Y_{lm'}(\hat{\mathbf{s}})$$

$$= \sum_{m'=-l}^{l} d_{m'm}^{l}(\theta_{\hat{\mathbf{k}}}) \sum_{m''=-l}^{l} d_{m'm}^{l}(\theta_{\hat{\mathbf{k}}})Y_{lm''}(\hat{\mathbf{s}}) = Y_{lm}(\hat{\mathbf{s}}), \qquad (6.7)$$

where we used (Varshalovich, Moskalev, and Khersonskii, 1988)

$$\sum_{m'=-l}^{l} d_{m'm}^{l}(\theta_{\hat{\mathbf{k}}})d_{m'm''}^{l}(\theta_{\hat{\mathbf{k}}}) = \delta_{mm''} \,. \qquad (6.8)$$

Let us consider $\theta_{\hat{\mathbf{k}}}, \varphi_{\hat{\mathbf{k}}}$ for $\hat{\mathbf{k}} = \hat{\mathbf{k}}(\nu, \mathbf{q})$, which is a complex vector that satisfies $\hat{\mathbf{k}} \cdot \hat{\mathbf{k}} = 1$. We take the the branch cut of a square root along the positive real axis in the complex plane. Then for $z \in \mathbb{C}$, we have $0 \le \arg(\sqrt{z}) < \pi$. We have for $\theta_{\hat{\mathbf{k}}}$,

$$\cos\theta_{\hat{\mathbf{k}}(\nu,\mathbf{q})} = \hat{\mathbf{z}} \cdot \hat{\mathbf{k}} = \hat{k}_z(\nu q), \quad \sin\theta_{\hat{\mathbf{k}}(\nu,\mathbf{q})} = \sqrt{1 - \left(\cos\theta_{\hat{\mathbf{k}}(\nu,\mathbf{q})}\right)^2} = \sqrt{\hat{k}_x^2 + \hat{k}_y^2} = i|\nu q| \,. \qquad (6.9)$$

We note that $\theta_{\hat{\mathbf{k}}(\nu,\mathbf{q})}$ depends on $\nu, \mathbf{q}$ as $\nu q$. Since $\theta_{\hat{\mathbf{k}}(\nu,\mathbf{q})}$ is complex, we write the complex angle as $\theta_{\hat{\mathbf{k}}(\nu,\mathbf{q})} = i\tau(\nu q)$ (Panasyuk, Schotland, and Markel, 2006). We can write

$$d_{m'm}^{l}\left(\theta_{\hat{\mathbf{k}}(\nu,\mathbf{q})}\right) = d_{m'm}^{l}[i\tau(\nu q)] \,. \qquad (6.10)$$

These $d_{m'm}^{l}$ are analytically continued Wigner's $d$-matrices. We have

$$d_{00}^{0}[i\tau(\nu q)] = 1 \,, \qquad (6.11)$$

$$d_{00}^{1}[i\tau(\nu q)] = \sqrt{1 + (\nu q)^2} \,, \quad d_{01}^{1}[i\tau(\nu q)] = \frac{i}{\sqrt{2}}|\nu q| \,, \qquad (6.12)$$

and

$$d_{1,\pm 1}^{1}[i\tau(\nu q)] = \frac{1 \pm \sqrt{1 + (\nu q)^2}}{2} \,. \qquad (6.13)$$

Note that

$$d_{mm'}^{l} = (-1)^{m+m'}d_{-m,-m'}^{l} = (-1)^{m+m'}d_{m'm}^{l} \,. \qquad (6.14)$$

Wigner $d$-matrices can be computed recursively (Blanco, Flórez, and Bermejo, 1997; Courant and Hilbert, 1953; Edmonds, 1957; Machida, Panasyuk, Schotland, and Markel, 2010). For $\varphi_{\hat{\mathbf{k}}}$, we have

$$\cos\varphi_{\hat{\mathbf{k}}(\nu,\mathbf{q})} = \frac{\hat{k}_x}{\sqrt{\hat{k}_x^2 + \hat{k}_y^2}} = -\frac{\nu}{|\nu|}\cos\varphi_{\mathbf{q}} \,, \quad \sin\varphi_{\hat{\mathbf{k}}(\nu,\mathbf{q})} = \frac{\hat{k}_y}{\sqrt{\hat{k}_x^2 + \hat{k}_y^2}} = -\frac{\nu}{|\nu|}\sin\varphi_{\mathbf{q}} \,. \qquad (6.15)$$

Thus we obtain

$$\varphi_{\hat{\mathbf{k}}(\nu,\mathbf{q})} = \begin{cases} \varphi_{\mathbf{q}} + \pi, & \nu > 0, \\ \varphi_{\mathbf{q}}, & \nu < 0. \end{cases} \qquad (6.16)$$

## 7. Three-dimensional $F_N$ method

The letter F in $F_N$ means facile (Siewert, 1978). The $F_N$ method has been developed in the one-dimensional transport theory. Here, we explore how the $F_N$ method is extended to three dimensions by using rotated reference frames.

In the $F_N$ method, there is no need to evaluate singular functions although the fact that the specific intensity consists of singular eigenfunctions is used. In one dimension, the RTE was solved by the $F_N$ method in the slab geometry for isotropic scattering (Grandjean and Siewert, 1979; Siewert and Benoist, 1979) and anisotropic scattering without (Devaux and Siewert, 1980; Garcia and Siewert, 1985; Siewert, 1978) and with (Garcia and Siewert, 1992; Garcia and Siewert, 1998) azimuthal dependence. The method was also extended to multigroup (Garcia and Siewert, 1981). After finding the specific intensity on the boundary, we can further calculate the specific intensity inside the medium (Garcia and Siewert, 1985). The uniqueness of the solution to the key $F_N$ equation was proved (Larsen, 1982). For isotropic scattering, the three dimensional RTE was solved with the $F_N$ method (Dunn and Siewert, 1985; Siewert and Dunn, 1983) using the pseudo-problem (Williams, 1982), which is based on plane-wave decomposition. See the review article by Garcia (Garcia, 1985).

### 7.1. Half-space in three dimensions

Let us consider light propagation in the half-space which is occupied by a homogeneous random medium (i.e., the absorption and scattering coefficients are constant). The specific intensity $I(\mathbf{r},\hat{\mathbf{s}})$ of light is governed by the following RTE.

$$\begin{cases} \hat{\mathbf{s}}\cdot\nabla I(\mathbf{r},\hat{\mathbf{s}})+\mu_t I(\mathbf{r},\hat{\mathbf{s}})=\mu_s\int_{\mathbb{S}^2}p(\hat{\mathbf{s}},\hat{\mathbf{s}}')I(\mathbf{r},\hat{\mathbf{s}}')d\hat{\mathbf{s}}'+S(\mathbf{r},\hat{\mathbf{s}}), & z>0, \\ I(\mathbf{r},\hat{\mathbf{s}})=0, \quad z=0, \quad 0<\mu\le 1, \\ I(\mathbf{r},\hat{\mathbf{s}})\to 0, \quad z\to\infty. \end{cases} \tag{7.1}$$

We assume that the internal source $S(\mathbf{r},\hat{\mathbf{s}})$ is compact and smooth such that the solution $I(\mathbf{r},\hat{\mathbf{s}})$ of the RTE can be expanded by a finite number of spherical harmonics to a good approximation. By the Placzek lemma (Case, de Hoffmann, and Placzek, 1953), we can consider the following RTE for the whole space.

$$\begin{cases} \hat{\mathbf{s}}\cdot\nabla\psi(\mathbf{r},\hat{\mathbf{s}})+\mu_t\psi(\mathbf{r},\hat{\mathbf{s}})=\mu_s\int_{\mathbb{S}^2}p(\hat{\mathbf{s}},\hat{\mathbf{s}}')\psi(\mathbf{r},\hat{\mathbf{s}}')d\hat{\mathbf{s}}'+\Theta(z-0^+)S(\mathbf{r},\hat{\mathbf{s}})+\mu I(\mathbf{r},\hat{\mathbf{s}})\delta(z), & z\in\mathbb{R}, \\ \psi(\mathbf{r},\hat{\mathbf{s}})\to 0, \quad |z|\to\infty. \end{cases} \tag{7.2}$$

We have the relation:

$$\psi(\mathbf{r},\hat{\mathbf{s}})=\begin{cases} I(\mathbf{r},\hat{\mathbf{s}}), & z>0, \\ 0, & z<0. \end{cases} \tag{7.3}$$

We note that the jump condition below is derived from the above RTE.

$$\psi(\boldsymbol{\rho},0^+,\hat{\mathbf{s}})-\psi(\boldsymbol{\rho},0^-,\hat{\mathbf{s}})=I(\boldsymbol{\rho},0,\hat{\mathbf{s}}). \tag{7.4}$$

Using the fundamental solution in (4.66), (4.71), we obtain

$$\psi(\mathbf{r},\hat{\mathbf{s}})=\int_{\mathbb{S}^2}\int_{\mathbb{R}^2}G(\mathbf{r},\hat{\mathbf{s}};\boldsymbol{\rho}',0,\hat{\mathbf{s}}')\mu' I(\boldsymbol{\rho}',0,\hat{\mathbf{s}}')\,d\boldsymbol{\rho}'d\hat{\mathbf{s}}' + \int_{\mathbb{S}^2}\int_0^\infty\int_{\mathbb{R}^2}G(\mathbf{r},\hat{\mathbf{s}};\boldsymbol{\rho}',z',\hat{\mathbf{s}}')S(\boldsymbol{\rho}',z',\hat{\mathbf{s}}')d\boldsymbol{\rho}'dz'd\hat{\mathbf{s}}' \tag{7.5}$$

for $\mathbf{r}\in\mathbb{R}^3$, $\hat{\mathbf{s}}\in\mathbb{S}^2$. By letting $z\to 0^+$, we have

$$I(\boldsymbol{\rho},0,\hat{\mathbf{s}}) = \int_{\mathbb{S}^2}\int_{\mathbb{R}^2} G(\boldsymbol{\rho},0^+,\hat{\mathbf{s}};\boldsymbol{\rho}',0,\hat{\mathbf{s}}')\mu' I(\boldsymbol{\rho}',0,\hat{\mathbf{s}}')\,d\boldsymbol{\rho}'d\hat{\mathbf{s}}'$$
$$+\int_{\mathbb{S}^2}\int_0^\infty\int_{\mathbb{R}^2} G(\boldsymbol{\rho},0,\hat{\mathbf{s}};\boldsymbol{\rho}',z',\hat{\mathbf{s}}')S(\boldsymbol{\rho}',z',\hat{\mathbf{s}}')d\boldsymbol{\rho}'dz'd\hat{\mathbf{s}}' \ . \tag{7.6}$$

In the Fourier space,

$$\tilde{I}(\mathbf{q},0,\hat{\mathbf{s}}) = \int_{\mathbb{S}^2} \tilde{G}(\mathbf{q};0^+,\hat{\mathbf{s}};0,\hat{\mathbf{s}}')\mu'\tilde{I}(\mathbf{q};\ 0,\hat{\mathbf{s}}')\,d\hat{\mathbf{s}}'$$
$$+\int_{\mathbb{S}^2}\int_0^\infty \tilde{G}(\mathbf{q};0,\hat{\mathbf{s}};\ z',\hat{\mathbf{s}}')\tilde{S}(\mathbf{q};\ z',\hat{\mathbf{s}}')dz'd\hat{\mathbf{s}}' \ , \tag{7.7}$$

where $\tilde{I},\tilde{S}$ were introduced similar to $\tilde{G}$ in (4.66). We note that $\tilde{G}(\mathbf{q};0^+,\hat{\mathbf{s}};0,\hat{\mathbf{s}}')$ contains only eigenmodes whose eigenvalues are positive and $\tilde{G}(\mathbf{q};0,\hat{\mathbf{s}};\ z',\hat{\mathbf{s}}')$ contains only eigenmodes whose eigenvalues are negative. Hence, if we multiply (7.7) by $\mu\mathcal{R}_{\hat{\mathbf{k}}(-\nu,\mathbf{q})}\Phi_{-\nu}^{m*}(\hat{\mathbf{s}})$, where $\nu>0$, and integrate over $\hat{\mathbf{s}}$, we obtain

$$\int_{\mathbb{S}^2} \mu\big(\mathcal{R}_{\hat{\mathbf{k}}(-\nu,\mathbf{q})}\Phi_{-\nu}^{m*}(\hat{\mathbf{s}})\,\big)\tilde{I}(\mathbf{q},0,\hat{\mathbf{s}})\,d\hat{\mathbf{s}}$$
$$=\int_{\mathbb{S}^2}\int_0^\infty \left(\mathcal{R}_{\hat{\mathbf{k}}(-\nu,\mathbf{q})}\Phi_{-\nu}^{m*}(\hat{\mathbf{s}})\right)e^{-\mu_t\hat{k}_z(\nu q)z/\nu}\tilde{S}(\mathbf{q};\ z,\hat{\mathbf{s}})dzd\hat{\mathbf{s}} \ . \tag{7.8}$$

By slightly rewriting (7.8), we obtain

$$\int_{\mathbb{S}^2_+} \mu\big(\mathcal{R}_{\hat{\mathbf{k}}(-\nu,\mathbf{q})}\Phi_{-\nu}^{m*}(-\hat{\mathbf{s}})\,\big)\tilde{I}(\mathbf{q},0,-\hat{\mathbf{s}})\,d\hat{\mathbf{s}}$$
$$=-\int_{\mathbb{S}^2}\int_0^\infty \left(\mathcal{R}_{\hat{\mathbf{k}}(-\nu,\mathbf{q})}\Phi_{-\nu}^{m*}(\hat{\mathbf{s}})\right)e^{-\mu_t\hat{k}_z(\nu q)z/\nu}\tilde{S}(\mathbf{q};\ z,\hat{\mathbf{s}})dzd\hat{\mathbf{s}} \ . \tag{7.9}$$

Let us approximate $\tilde{I}(\mathbf{q},0,-\hat{\mathbf{s}})$ as

$$\tilde{I}(\mathbf{q},0,-\hat{\mathbf{s}}) \approx \sum_{m=-l_{\max}}^{l_{\max}}\ \sum_{\alpha=0}^{\lfloor (N-|m|)/2\rfloor} c_{|m|+2\alpha,m}(\mathbf{q})Y_{|m|+2\alpha,m}(\hat{\mathbf{s}}) \ , \tag{7.10}$$

where $N$ is the largest degree of spherical harmonics used in the expansion in (7.10). In general, $N$ is taken to be much greater than $l_{\max}$. Furthermore, $N$ can be chosen depending on $m$. See Garcia and Siewert, 1998. In this section, we set $N=l_{\max}$. In (7.10), only same-parity degrees are taken because the three-term recurrence relation of associated Legendre polynomials implies that $Y_{lm}(\hat{\mathbf{s}})$ of opposite-parity $l$ are not independent. By substituting this approximate form for $\tilde{I}(\mathbf{q},0,-\hat{\mathbf{s}})$ in (7.9), we arrive at the following key $F_N$ equation:

$$\sum_{m=-l_{\max}}^{l_{\max}}\ \sum_{l=|m|,|m|+2,\cdots} A_{lm}^{m'}(\nu,\mathbf{q})c_{lm}(\mathbf{q}) = K^{m'}(\nu,\mathbf{q}) \ , \tag{7.11}$$

where $-l_{\max}\le m'\le l_{\max}$,

$$A_{lm}^{m'}(\nu,\mathbf{q}) = \int_{\mathbb{S}^2_+} \mu Y_{lm}(\hat{\mathbf{s}})\ \left(\mathcal{R}_{\hat{\mathbf{k}}(-\nu,\mathbf{q})}\Phi_{-\nu}^{m*}(-\hat{\mathbf{s}})\right)d\hat{\mathbf{s}} \ , \tag{7.12}$$

$$K^{m\prime}(\nu,\mathbf{q}) = \int_{\mathbb{S}^2_+} \mu\tilde{f}(\mathbf{q},,\hat{\mathbf{s}})\ \big(\mathcal{R}_{\hat{\mathbf{k}}(-\nu,\mathbf{q})}\Phi_{-\nu}^{m*}(-\hat{\mathbf{s}})\,\big)\,d\hat{\mathbf{s}}$$

$$-\int_{\mathbb{S}^2}\int_0^\infty \left(\mathcal{R}_{\hat{\mathbf{k}}(-\nu,\mathbf{q})}\Phi_{-\nu}^{m*}(\hat{\mathbf{s}})\right) e^{-\mu_t\hat{k}_z(\nu q)z/\nu}\tilde{S}(\mathbf{q};\, z,\hat{\mathbf{s}})dzd\hat{\mathbf{s}}\ . \tag{7.13}$$

Thus we can compute $\tilde{I}(\mathbf{q},0,-\hat{\mathbf{s}})$ by obtaining $c_{lm}(\mathbf{q})$ in the linear equation (7.11). We can write the matrix $A_{lm}^{m'}(\nu,\mathbf{q})$ as

$$A_{lm}^{m'}(\nu,\mathbf{q}) = A_{lm}^{m'}(\nu,q)e^{im\varphi_{\mathbf{q}}}\ . \tag{7.14}$$

After some calculation, we find that $A_{lm}^{m'}(\nu,q)$ is given as follows.

$$\begin{aligned}
A_{lm}^{m'}(\nu,q) &= (-1)^m\hat{k}_z(\nu q)\sqrt{\frac{\pi}{2l+1}}d_{mm'}^{l}[i\tau(\nu q)] \\
&\times\left(\sqrt{(l+1)^2-m'^2}g_{l+1}^{m'}(\nu)+\sqrt{l^2-m'^2}g_{l-1}^{m'}(\nu)\right) \\
&-i\frac{|\nu q|}{2}\sqrt{\frac{\pi}{2l+1}}(-1)^m\sum_{m''=-l}^{l}d_{mm''}^{l}[i\tau(\nu q)] \\
&\times\Big[\delta_{m'',m'-1}\left(\sqrt{(l-m'')(l-m')}g_{l-1}^{m'}(\nu)-\sqrt{(l+m'+1)(l+m')}g_{l+1}^{m'}(\nu)\right) \\
&\quad\delta_{m'',m'+1}\left(\sqrt{(l-m'+1)(l-m')}g_{l+1}^{m'}(\nu)-\sqrt{(l+m'')(l+m')}g_{l-1}^{m'}(\nu)\right)\Big] \\
&+\frac{\varpi\nu}{2}(-1)^l\sqrt{\frac{2l+1}{4\pi}\frac{(l-m)!}{(l+m)!}}[\mathrm{sgn}(m')]^{m'}\frac{\sqrt{(2|m'|)!}}{(2|m'|-1)!!} \\
&\times\sum_{m''=-|m'|}^{|m'|}(-1)^{m''}\sqrt{\frac{(|m'|-m'')!}{(|m'|+m'')!}}d_{m'',-m'}^{|m'|}[i\tau(\nu q)] \\
&\times\int_{\mathbb{S}_+^2}\frac{g^{m'}\left(-\nu,\hat{k}_z(\nu q)\mu-i\nu q\sqrt{1-\mu^2}\cos\varphi\right)}{\nu+\hat{k}_z(\nu q)\mu-i\nu q\sqrt{1-\mu^2}\cos\varphi}\mu P_{|m'|}^{m''}(\mu)P_l^m(\mu)e^{i(m+m'')\varphi}\,d\hat{\mathbf{s}}\ .
\end{aligned} \tag{7.15}$$

The number of columns of $A_{lm}^{m'}(\nu,q)$ is $N_{\mathrm{tot}}$, where

$$N_{\mathrm{tot}} = \begin{cases}(l_{\max}+2)^2/4, & l_{\max}\ \text{even}\,,\\ (l_{\max}+1)(l_{\max}+3)/4, & l_{\max}\ \text{odd}\,.\end{cases} \tag{7.16}$$

We choose the number of rows so that $A_{lm}^{m'}(\nu,q)$ becomes square. For this purpose, different collocation schemes have been proposed (Garcia and Siewert, 1981; Garcia and Siewert, 1985; Garcia and Siewert, 1998; McCormick and Sanchez, 1981).

Let $\xi$ be one of $N_{\mathrm{tot}}$ chosen values for $\nu$. If $K^{m'}(\xi,\mathbf{q})$ is independent of $\varphi_{\mathbf{q}}$ and $K^{-m'}=K^{m'}$, then

$$c_{lm}(\mathbf{q}) = c_{lm}(q)e^{-im\varphi_{\mathbf{q}}}\,,\quad c_{l,-m}(q) = (-1)^m c_{lm}(q). \tag{7.17}$$

After calculating $I(\boldsymbol{\rho},0,\hat{\mathbf{s}})$, we can obtain the specific intensity $I(\boldsymbol{\rho},z,\hat{\mathbf{s}})$ inside the medium ($z>0$) by using $I(\boldsymbol{\rho},0,\hat{\mathbf{s}})$. In this way, $I(\mathbf{r},\hat{\mathbf{s}})$ for the half-space can be computed by the three-dimensional $F_N$ method. It is also possible to consider the slab geometry.

### 7.2. Light from a spatially modulated source

Let us consider the following problem.

$$\begin{cases} \hat{\mathbf{s}} \cdot \nabla I(\mathbf{r}, \hat{\mathbf{s}}) + \mu_t I(\mathbf{r}, \hat{\mathbf{s}}) = \mu_s \int_{\mathbb{S}^2} p(\hat{\mathbf{s}}, \hat{\mathbf{s}}') I(\mathbf{r}, \hat{\mathbf{s}}') d\hat{\mathbf{s}}' \ , \ z > 0 \, , \\ I(\mathbf{r}, \hat{\mathbf{s}}) = f(\mathbf{r}, \hat{\mathbf{s}}) \, , \ \ z = 0 \, , \ \ 0 < \mu \le 1 \, , \\ I(\mathbf{r}, \hat{\mathbf{s}}) \to 0 \, , \ \ z \to \infty \, . \end{cases} \tag{7.18}$$

Here, the incident beam is given by

$$f(\mathbf{r}, \hat{\mathbf{s}}) = e^{-i\mu_t \mathbf{q}_0 \cdot \boldsymbol{\rho}} \delta(\hat{\mathbf{s}} - \hat{\mathbf{s}}_0) \, , \tag{7.19}$$

where $\mathbf{q}_0 \in \mathbb{R}^2$ is a given two-dimensional real vector. The cosine $\mu_0$ of the polar angle of $\hat{\mathbf{s}}_0$ is assumed to be positive. Although the formulation for general $\mathbf{q}_0, \hat{\mathbf{s}}_0$ is possible, in particular we set

$$\mathbf{q}_0 = \mu_t q_0 \hat{\mathbf{x}} \ , \ \ \hat{\mathbf{s}}_0 = \hat{\mathbf{z}} \, , \tag{7.20}$$

where $\hat{\mathbf{x}} = (1,0,0)^T$. Here, $T$ means transpose.

We begin by decomposing $I(\mathbf{r}, \hat{\mathbf{s}})$ as

$$I(\mathbf{r}, \hat{\mathbf{s}}) = I_b(\mathbf{r}, \hat{\mathbf{s}}) + I_s(\mathbf{r}, \hat{\mathbf{s}}), \tag{7.21}$$

where

$$\begin{cases} \hat{\mathbf{s}} \cdot \nabla I_b(\mathbf{r}, \hat{\mathbf{s}}) + \mu_t I_b(\mathbf{r}, \hat{\mathbf{s}}) = 0 \, , \ \ z > 0 \, , \\ I_b(\mathbf{r}, \hat{\mathbf{s}}) = f(\mathbf{r}, \hat{\mathbf{s}}), \ \ z = 0 \, , \ \ 0 < \mu \le 1 \, , \end{cases} \tag{7.22}$$

$$\begin{cases} \hat{\mathbf{s}} \cdot \nabla I_s(\mathbf{r}, \hat{\mathbf{s}}) + \mu_t I_s(\mathbf{r}, \hat{\mathbf{s}}) = \mu_s \int_{\mathbb{S}^2} p(\hat{\mathbf{s}}, \hat{\mathbf{s}}') I_s(\mathbf{r}, \hat{\mathbf{s}}') d\hat{\mathbf{s}}' + S[I_b](\mathbf{r}, \hat{\mathbf{s}}), \ \ z > 0 \, , \\ I_s(\mathbf{r}, \hat{\mathbf{s}}) = 0 \, , \ \ z = 0 \, , \ \ 0 < \mu \le 1 \, . \end{cases} \tag{7.23}$$

Here, the source term for $I_s(\mathbf{r}, \hat{\mathbf{s}})$ is given by

$$S[I_b](\mathbf{r}, \hat{\mathbf{s}}) = \mu_s \int_{\mathbb{S}^2} p(\hat{\mathbf{s}}, \hat{\mathbf{s}}') I_b(\mathbf{r}, \hat{\mathbf{s}}') d\hat{\mathbf{s}}' \ . \tag{7.24}$$

Since the ballistic term $I_b(\mathbf{r}, \hat{\mathbf{s}})$ is obtained as

$$I_b(\mathbf{r}, \hat{\mathbf{s}}) = e^{-i\mu_t q_0 x} \delta(\hat{\mathbf{s}} - \hat{\mathbf{z}}), \tag{7.25}$$

we have

$$S[I_b](\mathbf{r}, \hat{\mathbf{s}}) = \frac{\mu_s}{4\pi} e^{-i\mu_t q_0 x} e^{-\mu_t z} \sum_{l=0}^{l_{\max}} \beta_l P_l(\mu) \ . \tag{7.26}$$

Let us consider how $I_s(\mathbf{r}, \hat{\mathbf{s}})$ is computed. In this case of the spatially modulated source, we obtain

$$K^{m'}(\xi, \mathbf{q}) = \breve{K}^{m'}(\xi, \mathbf{q}_0) \delta(\mathbf{q} - \mathbf{q}_0) \, , \tag{7.27}$$

where

$$\breve{K}^{m'}(\xi, \mathbf{q}_0) = \frac{-2\pi^2 \varpi \xi}{\xi + \hat{k}_z(\xi q_0)} \sum_{l=|m'|}^{l_{\max}} (-1)^l \beta_l d^l_{0m'}[i\tau(\xi q_0)] g^{m'}_l(\xi) \ . \tag{7.28}$$

Correspondingly we define $\check{c}_{lm}(\mathbf{q}_0)$ as

$$c_{lm}(\mathbf{q}) = \check{c}_{lm}(\mathbf{q}_0) \delta(\mathbf{q} - \mathbf{q}_0) \, . \tag{7.29}$$

Then the key $F_N$ equation is written as

$$\sum_{m=-l_{\max}}^{l_{\max}} \sum_{l=|m|,|m|+2,\cdots} A^{m'}_{lm}(\xi, \mathbf{q}_0) \check{c}_{lm}(\mathbf{q}_0) = \breve{K}^{m'}(\xi, \mathbf{q}_0) \ . \tag{7.30}$$

To establish the key $F_N$ equation, we need to take $\lfloor (l_{\max} - m)/2 \rfloor + 1$ values for $\nu$ for each $m$ ($N_{\text{tot}}$ values in total). In addition to eigenvalues

$$\xi_j = \nu_{j-1}^m, \quad j = 1, \cdots, M^m, \tag{7.31}$$

we choose

$$\xi_j = \cos\left(\frac{\pi}{2} \frac{j - M^m}{2\left\lfloor \frac{(l_{\max} - m)}{2} \right\rfloor + 3 - M^m}\right), \quad j = M^m + 1, \cdots, \left\lfloor \frac{l_{\max} - m}{2} \right\rfloor + 1. \tag{7.32}$$

Thus $I_s(\mathbf{r}, \hat{\mathbf{s}})$ is obtained using the formulation which was developed in Sec. 7.1.

Let us calculate the hemispheric flux:

$$J_+(\boldsymbol{\rho}; \mathbf{q}_0) = \int_0^{2\pi} \int_0^1 \mu I(\mathbf{r}, -\hat{\mathbf{s}}) d\mu d\varphi, \quad z = 0. \tag{7.33}$$

We obtain

$$J_+(\boldsymbol{\rho}; \mathbf{q}_0) = \frac{1}{4\pi^{3/2}} e^{-iq_0 x} \sum_{l=0,2,\cdots} \frac{\sqrt{2l+1}(-1)^{\frac{l}{2}+1} l!}{2^l (l-1)(l+2)\left[\left(\frac{l}{2}\right)!\right]^2} \check{c}_{l0}(q_0). \tag{7.34}$$

See Machida, 2015 for more details.

# 8. Three-dimensional analytical discrete ordinates

One way to approximate the fundamental solution in (4.66), (4.71) is to approximate the RTE from the beginning instead of finding approximate solutions to the RTE. Here, we replace the integral in the RTE by a sum. That is, we apply rotated reference frames to the discrete ordinates method.

## 8.1. Gauss-Legendre quadrature

Let $w_i$ $(i = 1, \cdots, N)$ be weights of the Gauss-Legendre quadrature. We discretize $\mu = \cos\theta$ as

$$0 < \mu_1 < \mu_2 < \cdots < \mu_N < 1, \quad -1 < \mu_{2N} < \cdots < \mu_{N+2} < \mu_{N+1} < 0, \tag{8.1}$$

where $\mu_{N+i} = -\mu_i$ $(i = 1, \cdots, N)$. The pairs $(w_i, \mu_i)$ can be calculated, for example, using the Golub-Welsch algorithm (Golub and Welsch, 1969). Instead of $\hat{\mathbf{s}}$, we introduce $\hat{\mathbf{s}}_i$ as

$$\hat{\mathbf{s}}_i = \begin{pmatrix} \boldsymbol{\omega}_i \\ \mu_i \end{pmatrix}, \quad \boldsymbol{\omega}_i = \begin{pmatrix} \sqrt{1-\mu_i^2}\cos\varphi \\ \sqrt{1-\mu_i^2}\sin\varphi \end{pmatrix}, \quad 0 \le \varphi < 2\pi, \quad i = 1, \cdots, 2N. \tag{8.2}$$

Then we have

$$\boldsymbol{\omega}_i \cdot \mathbf{q} = q\sqrt{1-\mu_i^2}\cos(\varphi - \varphi_{\mathbf{q}}). \tag{8.3}$$

We define

$$\omega_i^0(q, \varphi) = q\sqrt{1-\mu_i^2}\cos\varphi. \tag{8.4}$$

## 8.2. Discrete ordinates

Let us consider the numerical evaluation of the fundamental solution (4.66), (4.71). We begin by performing the ballistic subtraction for the RTE (4.62):

$$I_s(\mathbf{r},\hat{\mathbf{s}}) = G(\mathbf{r},\hat{\mathbf{s}};\mathbf{r}_0,\hat{\mathbf{s}}_0) - I_b(\mathbf{r},\hat{\mathbf{s}}). \tag{8.5}$$

We have

$$\hat{\mathbf{s}}\cdot\nabla I_b(\mathbf{r},\hat{\mathbf{s}}) + \mu_t I_b(\mathbf{r},\hat{\mathbf{s}}) = \delta(\mathbf{r}-\mathbf{r}_0)\delta(\hat{\mathbf{s}}-\hat{\mathbf{s}}_0)\ , \tag{8.6}$$

and

$$\hat{\mathbf{s}}\cdot\nabla I_s(\mathbf{r},\hat{\mathbf{s}}) + \mu_t I_s(\mathbf{r},\hat{\mathbf{s}}) = \mu_s\int_{\mathbb{S}^2} p(\hat{\mathbf{s}},\hat{\mathbf{s}}')I_s(\mathbf{r},\hat{\mathbf{s}}')d\hat{\mathbf{s}}' + S[I_b](\mathbf{r},\hat{\mathbf{s}})\,. \tag{8.7}$$

Here, $S[I_b](\mathbf{r},\hat{\mathbf{s}})$ was given in (7.24) but the ballistic term is obtained as

$$I_b(\mathbf{r},\hat{\mathbf{s}}) = \Theta\big(\mu_0(z-z_0)\big)\frac{\mu_t}{|\mu_0|}\delta\left(\boldsymbol{\rho}-\boldsymbol{\rho}_0-\frac{z-z_0}{\mu_0}\boldsymbol{\omega}_0\right)e^{-\mu_t(z-z_0)/\mu_0}\delta(\hat{\mathbf{s}}-\hat{\mathbf{s}}_0)\,, \tag{8.8}$$

where $\boldsymbol{\omega}_0,\mu_0$ are components of $\hat{\mathbf{s}}_0$. Alternatively, $I_b(\mathbf{r},\hat{\mathbf{s}})$ can be expressed as

$$I_b(\mathbf{r},\hat{\mathbf{s}}) = \frac{1}{|\mathbf{r}-\mathbf{r}_0|^2}e^{-\mu_t|\mathbf{r}-\mathbf{r}_0|}\delta\left(\frac{\mathbf{r}-\mathbf{r}_0}{|\mathbf{r}-\mathbf{r}_0|}-\hat{\mathbf{s}}_0\right)\delta(\hat{\mathbf{s}}-\hat{\mathbf{s}}_0)\,. \tag{8.9}$$

Let us set

$$\mathbf{r}_0 = \mathbf{0}. \tag{8.10}$$

We obtain

$$S[I_b](\mathbf{r},\hat{\mathbf{s}}) = \Theta(\mu_0 z)\frac{\mu_s}{|\mu_0|}p(\hat{\mathbf{s}},\hat{\mathbf{s}}_0)\delta\left(\boldsymbol{\rho}-\frac{z}{\mu_0}\boldsymbol{\omega}_0\right)e^{-\mu_t z/\mu_0}\,. \tag{8.11}$$

Let us consider the Fourier transform of the scattering part:

$$\tilde{I}_s(\mathbf{q},z,\hat{\mathbf{s}}) = \int_{\mathbb{R}^2} e^{-i\mu_t\mathbf{q}\cdot\boldsymbol{\rho}}I_s(\mathbf{r},\hat{\mathbf{s}})\,d\boldsymbol{\rho}\ . \tag{8.12}$$

This $\tilde{I}_s(\mathbf{q},z,\hat{\mathbf{s}})$ satisfies

$$\left(\mu\frac{\partial}{\partial z}+(1+i\boldsymbol{\omega}\cdot\boldsymbol{q})\mu_t\right)\tilde{I}_s(\mathbf{q},z,\hat{\mathbf{s}}) = \mu_s\int_{\mathbb{S}^2} p(\hat{\mathbf{s}},\hat{\mathbf{s}}')\tilde{I}_s(\mathbf{q},z,\hat{\mathbf{s}}')d\hat{\mathbf{s}}' + \tilde{S}[I_b](\mathbf{q},z,\hat{\mathbf{s}})\,, \tag{8.13}$$

where $\tilde{S}[I_b](\mathbf{q},z,\hat{\mathbf{s}})$ is the Fourier transform of $S[I_b](\mathbf{r},\hat{\mathbf{s}})$ which is defined similar to (8.12).

Now we replace the integral in (8.13) with a sum:

$$\left(\mu_i\frac{\partial}{\partial z}+(1+i\boldsymbol{\omega}_i\cdot\boldsymbol{q})\mu_t\right)\check{I}_s(\mathbf{q},z,\hat{\mathbf{s}}_i) = \mu_s\sum_{i'=1}^{2N} w_{i'}\int_0^{2\pi} p(\hat{\mathbf{s}}_i,\hat{\mathbf{s}}_{i'})\check{I}_s(\mathbf{q},z,\hat{\mathbf{s}}_{i'})d\varphi' + \tilde{S}(\mathbf{q},z,\hat{\mathbf{s}}_i)\,. \tag{8.14}$$

The following equality numerically holds for sufficiently large $N$.

$$\tilde{I}_s(\mathbf{q},z,\hat{\mathbf{s}}) = \check{I}_s(\mathbf{q},z,\hat{\mathbf{s}}_i)\,. \tag{8.15}$$

### 8.3. Eigenmodes

Let us consider the homogeneous equation:

$$\left(\mu_i\frac{\partial}{\partial z}+(1+i\boldsymbol{\omega}_i\cdot\boldsymbol{q})\mu_t\right)\check{I}(\mathbf{q},z,\hat{\mathbf{s}}_i) = \mu_s\sum_{i'=1}^{2N} w_{i'}\int_0^{2\pi} p(\hat{\mathbf{s}}_i,\hat{\mathbf{s}}_{i'})\check{I}(\mathbf{q},z,\hat{\mathbf{s}}_{i'})d\varphi'\ . \tag{8.16}$$

As seen below, $\check{I}(\mathbf{q},z,\hat{\mathbf{s}}_i)$ does not contain generalized functions. Let us assume the following separated solution.

$$\check{I}(\mathbf{q},z,\hat{\mathbf{s}}_i) = e^{-\mu_t\hat{k}_z(\xi q)z/\xi}\mathcal{R}_{\hat{\mathbf{k}}(\xi,\mathbf{q})}\check{\Phi}_\xi^m(\hat{\mathbf{s}}_i), \tag{8.17}$$

where

$$\check{\Phi}_\xi^m(\hat{\mathbf{s}}_i) = \check{\phi}^m(\xi,\mu_i)(1-\mu_i{}^2)^{|m|/2}e^{im\varphi}\,. \tag{8.18}$$

The separation constant $\xi$ is not necessarily the same as $\xi$ used in Sec. 7. We will determine $\xi$ below. We impose the normalization condition:

$$\sum_{i=1}^{2N} w_i \check{\phi}^m(\xi,\mu_i)(1-\mu_i{}^2)^{|m|/2} = 1\,. \tag{8.19}$$

A remark is necessary. The normalization condition implies that the following equality numerically holds for a unit vector $\hat{\mathbf{k}} \in \mathbb{C}$.

$$\frac{1}{2\pi}\sum_{i=1}^{2N} w_i \int_0^{2\pi} \check{\phi}^m(\xi,\hat{\mathbf{s}}_i\cdot\hat{\mathbf{k}}) \left[1-\left(\hat{\mathbf{s}}_i\cdot\hat{\mathbf{k}}\right)^2\right]^{|m|} d\varphi = 1\,. \tag{8.20}$$

Similar to the calculation in Sec. 4, we arrive at the following equation within the approximation of discrete ordinates.

$$\left(1-\frac{\mu_i}{\xi}\right)\check{\Phi}_\xi^m(\hat{\mathbf{s}}_i) = \varpi \sum_{l'=0}^{l_{\max}} \sum_{m'=-l'}^{l'} \frac{\beta_{l'}}{2l'+1} Y_{l'm'}(\hat{\mathbf{s}}_i) \sum_{i'=1}^{2N} w_{i'} \int_0^{2\pi} Y_{l'm'}^*(\hat{\mathbf{s}}_{i'})\check{\Phi}_\xi^m(\hat{\mathbf{s}}_{i'})\,d\varphi' \tag{8.21}$$

for $i = 1,\cdots,2N$. By a straightforward calculation, we can show

$$\sum_{i'=1}^{2N} w_{i'} \int_0^{2\pi} Y_{l'm'}^*(\hat{\mathbf{s}}_{i'})\check{\Phi}_\xi^m(\hat{\mathbf{s}}_{i'})d\varphi' = \delta_{mm'}(-1)^m\sqrt{(2l'+1)\pi}\,g_{l'}^m(\xi)\,, \tag{8.22}$$

where we used (see (4.21) and (4.27))

$$g_{l'}^m(\xi) = \sum_{i'=1}^{2N} w_{i'} p_{l'}^m(\mu_i)\check{\phi}^m(\xi,\mu_i)(1-\mu_i^2)^{|m|}\,. \tag{8.23}$$

We have

$$(\xi-\mu_i)\check{\phi}(\xi,\mu_i) = \frac{\varpi\xi}{2}\sum_{l'=|m|}^{l_{\max}} \beta_{l'} p_{l'}^m(\mu_i) g_{l'}^m(\xi)\,. \tag{8.24}$$

We obtain (Machida and Das, 2022)

$$\check{\Phi}_\xi^m(\hat{\mathbf{s}}_i) = \frac{(-1)^m\varpi\xi}{\xi-\mu_i}\sum_{l=|m|}^{l_{\max}} \sqrt{\frac{\pi}{2l+1}}\beta_l g_l^m(\xi) Y_{lm}(\hat{\mathbf{s}}_i)\,, \tag{8.25}$$

and thus

$$\check{I}(\mathbf{q},z,\hat{\mathbf{s}}_i) = \frac{(-1)^m\varpi\xi}{\xi-\hat{\mathbf{s}}_i\cdot\hat{\mathbf{k}}(\xi,\mathbf{q})} e^{-\mu_t\hat{k}_z(\xi q)z/\xi} \sum_{l=|m|}^{l_{\max}} \sqrt{\frac{\pi}{2l+1}}\beta_l g_l^m(\xi)\mathcal{R}_{\hat{\mathbf{k}}(\xi,\mathbf{q})} Y_{lm}(\hat{\mathbf{s}}_i)\,. \tag{8.26}$$

The separation constant $\xi$ is shown to be an eigenvalue of a matrix. We can calculate $\xi$ as follows. For this purpose we rewrite (8.24) as the following two equations:

$$\left(1-\frac{\mu_i}{\xi}\right)\check{\phi}^m(\xi,\mu_i) = \frac{\varpi}{2}\sum_{l'=|m|}^{l_{\max}} \beta_{l'} p_{l'}^m(\mu_i) \sum_{i'=1}^{N} w_{i'} p_{l'}^m(\mu_{i'})(1-\mu_{i'}^2)^{|m|} \times\left[\check{\phi}^m(\xi,\mu_{i'}) + (-1)^{l'+m}\check{\phi}^m(\xi,-\mu_{i'})\right], \tag{8.27}$$

$$\left(1+\frac{\mu_i}{\xi}\right)\check{\phi}^m(\xi,-\mu_i) = \frac{\varpi}{2}\sum_{l'=|m|}^{l_{\max}} \beta_{l'} p_{l'}^m(\mu_i) \sum_{i'=1}^{N} w_{i'} p_{l'}^m(\mu_{i'})(1-\mu_{i'}^2)^{|m|} \times\left[(-1)^{l'+m}\check{\phi}^m(\xi,\mu_{i'}) + \check{\phi}^m(\xi,-\mu_{i'})\right] \tag{8.28}$$

for $i = 1, \cdots, N$. Thus we arrive at matrix-vector equations:

$$\left(\mathbb{I}_N - \frac{1}{\nu}\Xi\right)\boldsymbol{\Phi}_+^m(\xi) = \frac{\varpi}{2}\left[W_+^m\boldsymbol{\Phi}_+^m(\xi) + W_-^m\boldsymbol{\Phi}_-^m(\xi)\right],$$
$$\left(\mathbb{I}_N + \frac{1}{\nu}\Xi\right)\boldsymbol{\Phi}_-^m(\xi) = \frac{\varpi}{2}\left[W_-^m\boldsymbol{\Phi}_+^m(\xi) + W_+^m\boldsymbol{\Phi}_-^m(\xi)\right], \tag{8.29}$$

where $\mathbb{I}_N$ is the $N$-dimensional identity matrix. Matrix $\Xi$ and vectors $\boldsymbol{\Phi}_\pm^m(\xi)$ are defined as

$$\Xi = \begin{pmatrix} \mu_1 & & \\ & \ddots & \\ & & \mu_N \end{pmatrix}, \quad \boldsymbol{\Phi}_\pm^m(\xi) = \begin{pmatrix} \breve{\phi}^m(\xi, \pm\mu_1) \\ \vdots \\ \breve{\phi}^m(\xi, \pm\mu_N) \end{pmatrix}. \tag{8.30}$$

Elements of the matrices $W_\pm^m$ are given by

$$\{W_\pm^m\}_{ij} = w_j \sum_{l=|m|}^{l_{\max}} \beta_l p_m^l(\pm\mu_i) p_m^l(\mu_j)\left(1 - \mu_j^2\right)^{|m|}. \tag{8.31}$$

The fact that $\boldsymbol{\Phi}_\pm^m(-\xi) = \boldsymbol{\Phi}_\mp^m(\xi)$ implies that $-\xi$ is an eigenvalue if $\xi$ is an eigenvalue. According to Barichello, 2011; Siewert, 2000, we introduce

$$\mathbf{U}^m(\xi) = \boldsymbol{\Phi}_+^m(\xi) + \boldsymbol{\Phi}_-^m(\xi), \quad \mathbf{V}^m(\xi) = \boldsymbol{\Phi}_+^m(\xi) - \boldsymbol{\Phi}_-^m(\xi). \tag{8.32}$$

By adding and subtracting two equations, we obtain

$$\mathbf{U}^m(\xi) - \frac{1}{\nu}\Xi\mathbf{V}^m(\xi) = \frac{\varpi}{2}(W_+^m + W_-^m)\mathbf{U}^m(\xi), \tag{8.33}$$

$$\mathbf{V}^m(\xi) - \frac{1}{\nu}\Xi\mathbf{U}^m(\xi) = \frac{\varpi}{2}(W_+^m - W_-^m)\mathbf{V}^m(\xi). \tag{8.34}$$

Hence we arrive at (Barichello and Siewert, 2000; Siewert and Wright, 1999)

$$E_-^m E_+^m \Xi \mathbf{U}^m = \frac{1}{\xi^2}\Xi\mathbf{U}^m, \tag{8.35}$$

where

$$E_\pm^m = \left[\mathbb{I}_N - \frac{\varpi}{2}(W_+^m \pm W_-^m)\right]\Xi^{-1}. \tag{8.36}$$

The eigenvalues $\xi$ can be calculated using (8.35):

$$\xi = \xi_n^m > 0, \quad n = 1, \cdots, N, \quad m = 0, \cdots, l_{\max}. \tag{8.37}$$

We have $\xi_n^{-m} = \xi_n^m$. In the case of isotropic scattering, it is proved that $\xi_n^0$ never coincide with $\mu_i$ (Siewert and Wright, 1999).

We note that the orthogonality relation (4.35) holds within the approximation of discrete ordinates. This implies that the fundamental solution (4.71) is obtained in the same manner:

$$\tilde{G}(\mathbf{q}; z, \hat{\mathbf{s}}; z_0, \hat{\mathbf{s}}_{i_0}) = \sum_{m=-l_{\max}}^{l_{\max}} \sum_{n=1}^{N} \frac{w_{i_0}}{2\pi\hat{k}_z(\xi_n^m q)\mathcal{N}(\xi_n^m, q)} e^{-\hat{k}_z(\xi_n^m q)|z-z_0|/\xi_n^m}$$
$$\times \mathcal{R}_{\hat{\mathbf{k}}(\pm\xi_n^m,\mathbf{q})} \breve{\Phi}_{\pm\xi_n^m}^m(\hat{\mathbf{s}}_i) \breve{\Phi}_{\pm\xi_n^m}^{m*}(\hat{\mathbf{s}}_{i_0}). \tag{8.38}$$

Within the approximation of discrete ordinates, the normalization factor is given by

$$\mathcal{N}(\xi, q) = 2\pi\hat{k}_z(\xi q)\sum_{i=1}^{2N} w_i\mu_i\left|\breve{\Phi}_\xi^m(\hat{\mathbf{s}}_i)\right|^2. \tag{8.39}$$

Although $\hat{\mathbf{s}}_{i_0}$ is arbitrary, let us assume

$$\hat{\mathbf{s}}_{i_0} = \hat{\mathbf{z}}. \tag{8.40}$$

Using (8.11), (8.38), we obtain

$$\check{I}_s(\mathbf{q},z,\hat{\mathbf{s}}_i) = \sum_{i'=1}^{2N} \int_0^{2\pi} \int_{-\infty}^{\infty} \tilde{G}(\mathbf{q};z,\hat{\mathbf{s}};z',\hat{\mathbf{s}}_{i'}) \tilde{S}(\mathbf{q},z',\hat{\mathbf{s}}_{i'})\, dz' d\varphi'$$

$$= \frac{\mu_s {\mu_t}^2}{2} \sum_{m=-l_{\max}}^{l_{\max}} \sum_{n=1}^{N} \frac{(-1)^m \xi_n^m}{\hat{k}_z(\xi_n^m q)\mathcal{N}(\xi_n^m, q)} \left[\Theta(z)(\psi_1+\psi_2) + \big(1-\Theta(z)\big)\psi_3\right] , \tag{8.41}$$

where

$$\psi_1 = C\left(z; 1, \frac{\xi_n^m}{\hat{k}_z(\xi_n^m q)}\right)\left(\mathcal{R}_{\hat{\mathbf{k}}(\xi_n^m,\mathbf{q})} \breve{\Phi}_{\xi_n^m}^m(\hat{\mathbf{s}}_i)\right) \sum_{l=|m|}^{l_{\max}} \beta_l g_l^m(\xi_n^m) d_{0m}^l[i\tau(\xi_n^m q)] , \tag{8.42}$$

$$\psi_2 = \frac{e^{-\mu_t z}\hat{k}_z(\xi_n^m q)}{\xi_n^m + \hat{k}_z(\xi_n^m q)} \left(\mathcal{R}_{\hat{\mathbf{k}}(-\xi_n^m,\mathbf{q})} \breve{\Phi}_{-\xi_n^m}^m(\hat{\mathbf{s}}_i)\right) \sum_{l=|m|}^{l_{\max}} \beta_l g_l^m(-\xi_n^m) d_{0m}^l[i\tau(\xi_n^m q)] , \tag{8.43}$$

$$\psi_3 = \frac{e^{\hat{k}_z(\xi_n^m q)\mu_t z/\xi_n^m}\hat{k}_z(\xi_n^m q)}{\xi_n^m + \hat{k}_z(\xi_n^m q)} \left(\mathcal{R}_{\hat{\mathbf{k}}(-\xi_n^m,\mathbf{q})} \breve{\Phi}_{-\xi_n^m}^m(\hat{\mathbf{s}}_i)\right) \sum_{l=|m|}^{l_{\max}} \beta_l g_l^m(-\xi_n^m) d_{0m}^l[i\tau(\xi_n^m q)] . \tag{8.44}$$

Here, the function $C(\tau;\varsigma,\eta)$ is defined as

$$C(\tau;\varsigma,\eta) = \frac{e^{-\tau/\varsigma} - e^{-\tau/\eta}}{\varsigma - \eta} . \tag{8.45}$$

## 8.4. Numerical examples

Let us compute the energy density (4.73). Since the energy density is symmetric about the $x$- $y$ plane due to the assumption (8.40), we can write $U(\mathbf{r}) = U(\rho, z)$ and obtain

$$U(\rho,z) = \frac{\mu_t}{2\pi c} \int_0^{\infty} q J_0(q\rho) \sum_{i=1}^{2N} w_i \int_0^{2\pi} \check{I}_s(\mathbf{q},\mathrm{z},\hat{\mathbf{s}}_i)\, d\varphi\, dq . \tag{8.46}$$

Since the following equality numerically holds for sufficiently large $N$:

$$\sum_{i=1}^{2N} w_i \int_0^{2\pi} \mathcal{R}_{\hat{\mathbf{k}}(\xi_n^m,\mathbf{q})} \breve{\Phi}_{\xi_n^m}^m(\hat{\mathbf{s}}_i)\, d\varphi = 2\pi\delta_{m0} , \tag{8.47}$$

we obtain

$$U(\rho,z) = \frac{\mu_s\mu_t}{2c} \int_0^{\infty} q J_0(q\rho) F(q,z)\, dq , \tag{8.48}$$

where

$$F(q,z) = \sum_{n=1}^{N} \frac{\xi_n^0}{\mathcal{N}(\xi_n^0, q)} \Bigg[\Theta(z) \frac{e^{-\mu_t z} - e^{-\hat{k}_z(\xi_n^0 q)\mu_t z/\xi_n^0}}{\hat{k}_z(\xi_n^0 q) - \xi_n^0} \sum_{l=0}^{l_{\max}} \beta_l g_l^0(\xi_n^0) d_{00}^l[i\tau(\xi_n^0 q)]$$
$$+\Theta(z) \frac{e^{-\mu_t z}}{\hat{k}_z(\xi_n^0 q) + \xi_n^0} \sum_{l=0}^{l_{\max}} (-1)^l \beta_l g_l^0(\xi_n^0) d_{00}^l[i\tau(\xi_n^0 q)]$$
$$+\big(1-\Theta(z)\big) \frac{e^{\hat{k}_z(\xi_n^0 q)\mu_t z/\xi_n^0}}{\hat{k}_z(\xi_n^0 q) + \xi_n^0} \sum_{l=0}^{l_{\max}} (-1)^l \beta_l g_l^0(\xi_n^0) d_{00}^l[i\tau(\xi_n^0 q)]\Bigg] . \tag{8.49}$$

To numerically evaluate $U(\rho, z)$, we can rewrite the integral in (8.48) using the asymptotic expression of the Bessel function $J_0$ (Machida and Das). The integrals of sinusoidally oscillating integrands can be evaluated with the double-exponential formula (Ooura and Mori, 1991; Ooura and Mori, 1999; Ogata, 2005).

Let us set

$$\mu_a = 0.01\ \text{mm}^{-1}, \quad \mu_s = 10\ \text{mm}^{-1}. \tag{8.50}$$

Moreover, in the scattering phase function we set

$$\beta_l = (2l + 1)\text{g}^l \tag{8.51}$$

with

$$\text{g} = 0.0 \text{ or } 0.9\,. \tag{8.52}$$

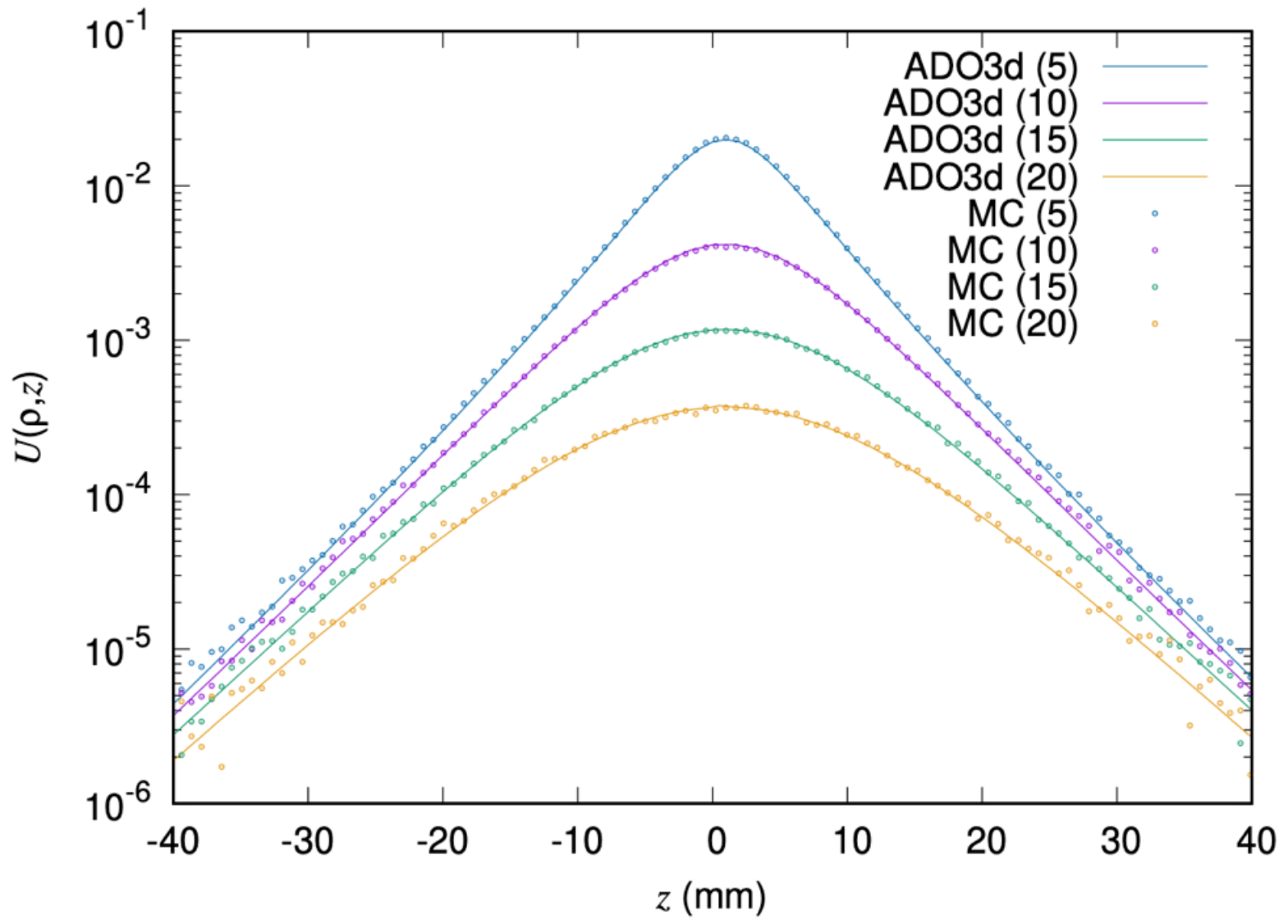

Figure 8.1. (Fig. 3 in Machida and Das) The energy density $U(\rho, z)$ in (8.48) is plotted (the factor $c$ is multiplied). Results are compared to Monte Carlo simulation. From the top, $\rho = 5$, 10, 15, and 20 mm. The parameters are $\mu_a = 0.01\ \text{mm}^{-1}$, $\mu_s = 10\ \text{mm}^{-1}$, $\text{g} = 0.9$.

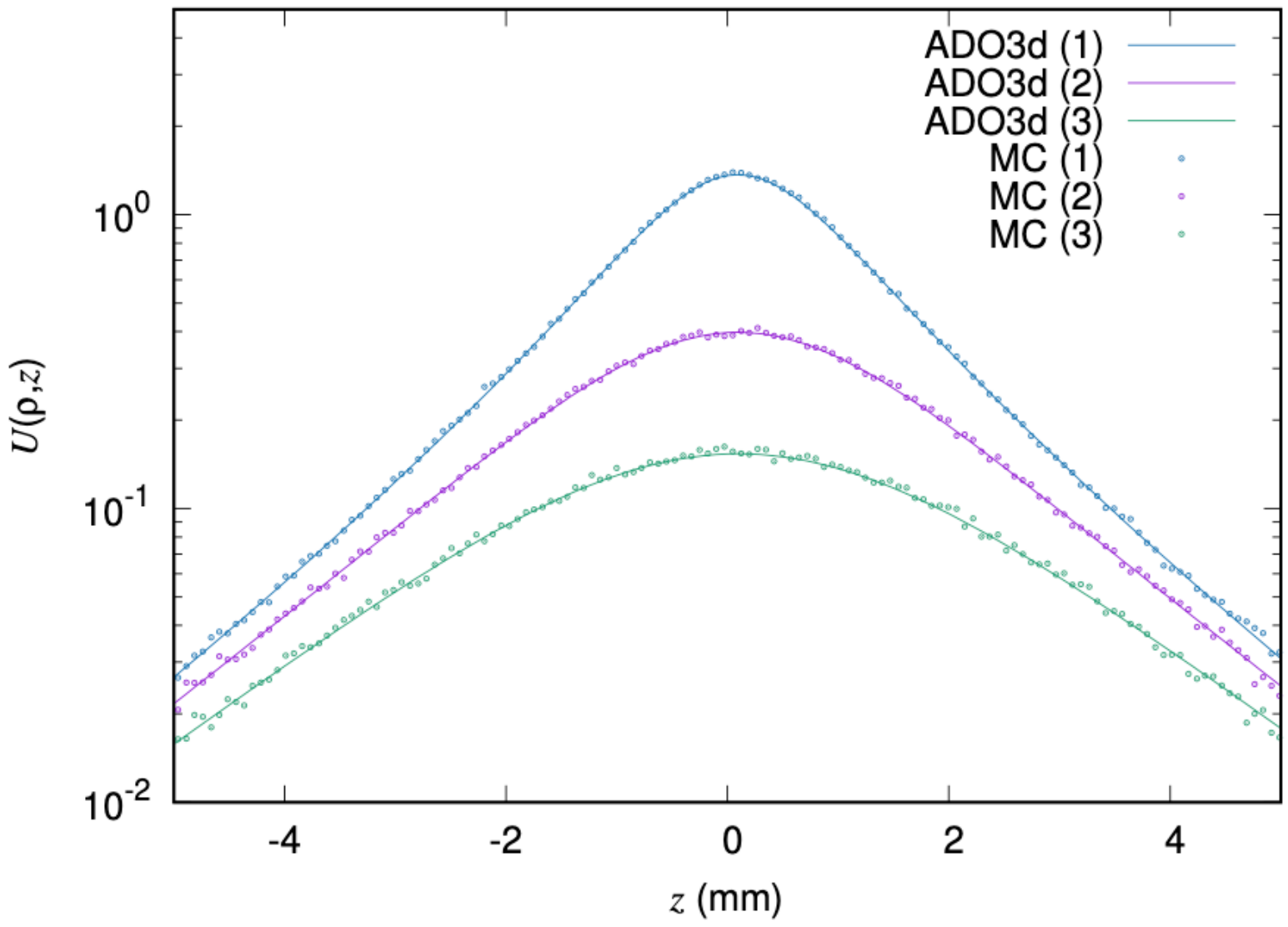

Figure 8.2. (Fig. 4 in Machida and Das) The energy density $U(\rho, z)$ in (8.48) is plotted (the factor $c$ is multiplied). Results are compared to Monte Carlo simulation. From the top, $\rho = 1$, 2, and 3 mm. The parameters are $\mu_a = 0.01\ \mathrm{mm}^{-1}$, $\mu_s = 10\ \mathrm{mm}^{-1}$, $\mathrm{g} = 0.0$.

### 8.5. Remarks

It is possible to develop the three-dimensional ADO in the presence of boundaries. In particular, by the three-dimensional ADO we can solve the RTE for the half space of $\mathbb{R}^3$, that is, for the domain $\{\mathbf{r} \in \mathbb{R}^3\,;\ \ -\infty < x < \infty\,,\ \ -\infty < y < \infty\,,\ \ 0 < z < \infty\}$. In this case, the specific intensity$\check{I}(\mathbf{r}, \hat{\mathbf{s}}_i)$ satisfies

$$\hat{\mathbf{s}} \cdot \nabla \check{I}(\mathbf{r}, \hat{\mathbf{s}}_i) + \mu_t \check{I}(\mathbf{r}, \hat{\mathbf{s}}_i) = \mu_s \sum_{i'=1}^{2N} w_{i'} \int_0^{2\pi} p(\hat{\mathbf{s}}_i, \hat{\mathbf{s}}_{i'}) \check{I}(\mathbf{r}, \hat{\mathbf{s}}_{i'}) d\varphi' \tag{8.53}$$

for $i = 1, \cdots, 2N$, $0 \le \varphi < 2\pi$ with the boundary at $z = 0$,

$$\check{I}(\mathbf{r}, \hat{\mathbf{s}}_i) = g(\mathbf{r}, \hat{\mathbf{s}}_i) \tag{8.54}$$

for $i = 1, \cdots, N$, $0 \le \varphi < 2\pi$. Here, $g(\mathbf{r}, \hat{\mathbf{s}}_i)$ is the incident beam. Note that $\check{I} \to 0$ as $|\mathbf{r}| \to \infty$.

The solutions with the three-dimensional ADO have the form of plane-wave decomposition. That is, the solutions depend on $x, y$ through the factor $e^{i\mu_t \mathbf{q}\cdot\boldsymbol{\rho}}$. Kim, 2004, Kim and Keller, 2003, and Kim and Moscoso, 2004 established a numerical method of discrete-ordinates with plane-wave decomposition for the three-dimensional RTE. Their method can be viewed as separation of variable assuming the form

$$I_{\mathbf{q},\lambda}(\mathbf{r}, \hat{\mathbf{s}}_i) = e^{i\mu_t \mathbf{q}\cdot\boldsymbol{\rho}} e^{-\lambda(\mathbf{q})z} V_{\lambda(\mathbf{q})}(\hat{\mathbf{s}}_i)\,, \tag{8.55}$$

where $\lambda(\mathbf{q})$ is a separation constant for a given $\mathbf{q} \in \mathbb{R}^2$ and $V_{\lambda(\mathbf{q})}$ is a function labeled by $\lambda(\mathbf{q})$. These $\lambda, V_\lambda$ are obtained as eigenvalues and eigenvectors of a matrix which depends on $\mathbf{q}$. By the comparison to the three-dimensional ADO, we notice that $\lambda$ is not necessarily calculated for each **q** because the **q**-dependence of the matrix can be eliminated as the matrices $E_{\pm}^{m}$ in (8.36) are independent of **q**. Indeed, we have the relation:

$$\lambda(\mathbf{q}) = \frac{\mu_t \hat{k}_z(\xi q)}{\xi} \ , \tag{8.56}$$

where $\xi$ is obtained from just one diagonalization.

Finally, we note that the use of rotated reference frames is not the only way to extend ADO in the one-dimensional transport theory to multidimensions with anisotropic scattering. The ADO-Nodal formulation makes it possible to solve the two-dimensional RTE for anisotropic scattering. See Barichello, Pazinatto, and Rui, 2020, Barichello, Rui, da Cunha, 2022, and references therein. An inverse source problem for the two-dimensional RTE was solved with the ADO-Nodal method (Pazinatto and Barichello, 2021).

# 9. Rotated reference frames with the spherical-harmonic expansion

## 9.1. Rotated spherical harmonics

The $P_N$ method is a well-known numerical scheme for the RTE. For the one-dimensional RTE, the numerical solution is given in terms of Legendre polynomials:

$$I(z,\mu) = \sum_{l=0}^{l_{\max}} \frac{2l+1}{4\pi} \psi_l(z) P_l(\mu) \ , \tag{9.1}$$

where $\psi_l$ is a function of $z$. More generally, the specific intensity can be expanded by spherical harmonics $Y_{lm}(\hat{\mathbf{s}})$ as

$$I(z,\hat{\mathbf{s}}) = \sum_{l=0}^{l_{\max}} \sum_{m=-l}^{l} \sqrt{\frac{2l+1}{4\pi}} \psi_{lm}(z) Y_{lm}(\hat{\mathbf{s}}) \ , \tag{9.2}$$

where $\psi_{lm}$ is a function of $z$. A naive way of extending the method to three dimensions is to write

$$I(\mathbf{r},\hat{\mathbf{s}}) = \sum_{l=0}^{l_{\max}} \sum_{m=-l}^{l} \sqrt{\frac{2l+1}{4\pi}} \psi_{lm}(\mathbf{r}) Y_{lm}(\hat{\mathbf{s}}) \ . \tag{9.3}$$

However, the above expansion does not work well. It is written "This rather awe-inspiring set of equations ... has perhaps only academic interest" (p. 219 in Case and Zweifel, 1967). The trick to use spherical harmonics for three dimensions is to use them in rotated reference frames. Let us write

$$I(\mathbf{r},\hat{\mathbf{s}}) = \frac{{\mu_t}^2}{(2\pi)^2} \sum_{m=-l_{\max}}^{l_{\max}} \int_{\mathbb{R}^2} e^{i\mu_t \mathbf{q}\cdot\boldsymbol{\rho}} \tilde{I}_m(\mathbf{q},z,\hat{\mathbf{s}}) \, d\mathbf{q} \ , \tag{9.4}$$

where

$$\tilde{I}_m(\mathbf{q},z,\hat{\mathbf{s}}) = \sum_{l'=0}^{l_{\max}} \sum_{m'=-l'}^{l'} \sum_{\lambda} \psi_{l'm'}^{(\lambda,m)}(\mathbf{q},z) \mathcal{R}_{\hat{\mathbf{k}}(\lambda,\mathbf{q})} Y_{l'm'}(\hat{\mathbf{s}}) \ . \tag{9.5}$$

The parameter $\lambda$ and the sum $\Sigma_\lambda$ will become clear below.

The expansion by $\mathcal{R}_{\hat{\mathbf{k}}}Y_{lm}(\hat{\mathbf{s}})$ was first done by Dede in 1964 (Dede, 1964). Then Kobayashi extended Dede's calculation (Kobayashi, 1977). However, until 2004, the expansion by $\mathcal{R}_{\hat{\mathbf{k}}}Y_{lm}(\hat{\mathbf{s}})$ did not bring efficient numerical methods even though the expansion was interesting. In 2004, Markel independently arrived at $\mathcal{R}_{\hat{\mathbf{k}}}Y_{lm}(\hat{\mathbf{s}})$ (Markel, 2004). He derived an eigenvalue problem and devised an efficient algorithm for the three-dimensional RTE.

Since then, as a numerical method to solve the three-dimensional RTE, the usefulness of the method of rotated reference frames has been numerically justified for an infinite medium (Markel, 2004; Panasyuk, Schotland, and Markel, 2006) , the slab geometry (Machida, Panasyuk, Schotland, and Markel, 2010), the half-space geometry (Liemert and Kienle, 2012b; Liemert and Kienle, 2012c; Liemert and Kienle, 2013b; Liemert and Kienle, 2013c). The method was developed for the three-dimensional time-dependent RTE for an infinite medium (Liemert and Kienle, 2012d; Liemert and Kienle, 2012e). The spherical-harmonic expansion was applied to the RTE in two dimensions (Liemert and Kienle, 2011; Liemert and Kienle, 2012a; Liemert and Kienle, 2012f; Liemert and Kienle, 2013a). The method was also used to experimentally determine optical properties of turbid media (Xu and Patterson, 2006a; Xu and Patterson, 2006b). The usefulness of the method of rotated reference frames has also been numerically justified in flatland (Liemert and Kienle, 2011; Liemert and Kienle, 2012a; Liemert and Kienle, 2012f; Liemert and Kienle, 2013a).

## 9.2. Separated solutions

Let us give $\psi_{l'm'}^{(\lambda,m)}(\mathbf{q},z)$ in (9.5) as

$$\psi_{l'm'}^{(\lambda,m)}(\mathbf{q},z) = c_{l'm'}^{m}(\lambda)e^{-\mu_t\hat{k}_z(\lambda q)z/\lambda}\,. \tag{9.6}$$

We take the Fourier transform for the homogeneous RTE (4.12) and substitute $\tilde{I}(\mathbf{q},z,\hat{\mathbf{s}})$ in (9.5). By this, we arrive at an eigenvalue problem. We obtain $\lambda$ after diagonalizing the resulting matrix.

Since we already know the analytical solution from Sec. 4, we see the relation,

$$\sum_{l'=0}^{l_{\max}}\sum_{m'=-l'}^{l'} c_{l'm'}^{m}\mathcal{R}_{\hat{\mathbf{k}}(\lambda,\mathbf{q})}Y_{l'm'}(\hat{\mathbf{s}}) = \mathcal{R}_{\hat{\mathbf{k}}(\lambda,\mathbf{q})}\sum_{l'=0}^{l_{\max}}\sum_{m'=-l'}^{l'} c_{l'm'}^{m}Y_{l'm'}(\hat{\mathbf{s}}) \approx \mathcal{R}_{\hat{\mathbf{k}}(\lambda,\mathbf{q})}\Phi_{\lambda}^{m}(\hat{\mathbf{s}}). \tag{9.7}$$

From (9.7), we see that (9.5) is obtained if the singular eigenfunction is expanded by spherical harmonics. The coefficients $c_{l'm'}^{m}$ can be calculated as follows (Machida, 2015). We multiply $\mathcal{R}_{\hat{\mathbf{k}}(\lambda,\mathbf{q})}Y_{l''m''}^{*}(\hat{\mathbf{s}})$ $(0 \le l'' \le l_{\max}, -l'' \le m'' \le l'')$ on both sides of (9.7) and integrate over $\hat{\mathbf{s}}$. We have

$$\int_{\mathbb{S}^2} \mathcal{R}_{\hat{\mathbf{k}}(\lambda,\mathbf{q})}Y_{l''m''}^{*}(\hat{\mathbf{s}})\sum_{l'=0}^{l_{\max}}\sum_{m'=-l'}^{l'} c_{l'm'}^{m}\mathcal{R}_{\hat{\mathbf{k}}(\lambda,\mathbf{q})}Y_{l'm'}(\hat{\mathbf{s}})\;d\hat{\mathbf{s}}$$

$$= \int_{\mathbb{S}^2} \mathcal{R}_{\hat{\mathbf{k}}(\lambda,\mathbf{q})}Y_{l''m''}^{*}(\hat{\mathbf{s}})\mathcal{R}_{\hat{\mathbf{k}}(\lambda,\mathbf{q})}\Phi_{\lambda}^{m}(\hat{\mathbf{s}})d\hat{\mathbf{s}}\,. \tag{9.8}$$

The left-hand side is computed as

$$\text{LHS of (9.8)} = \int_{\mathbb{S}^2} Y_{l''m''}^{*}(\hat{\mathbf{s}})\sum_{l'=0}^{l_{\max}}\sum_{m'=-l'}^{l'} c_{l'm'}^{m}(\lambda)Y_{l'm'}(\hat{\mathbf{s}})\;d\hat{\mathbf{s}} = c_{l''m''}^{m}(\lambda)\,. \tag{9.9}$$

The right-hand side is obtained as

$$\text{RHS of (9.8)} = \int_{\mathbb{S}^2} Y_{l''m''}^{*}(\hat{\mathbf{s}})\Phi_\lambda^m(\hat{\mathbf{s}})d\hat{\mathbf{s}} = (-1)^m\delta_{mm''}\sqrt{(2l''+1)\pi}g_{l''}^m(\lambda)\,. \tag{9.10}$$

Thus we obtain

$$c_{l'm'}^m(\lambda) = (-1)^m\delta_{mm'}\sqrt{(2l+1)\pi}g_l^m(\lambda)\,. \tag{9.11}$$

One remark is necessary for the case that $l_{\max} \to \infty$ in (9.7). Unlike the one-dimensional transport theory, which corresponds to the case of $\hat{\mathbf{k}} = \hat{\mathbf{z}}$, the first equality in (9.7) does not hold if $l_{\max} = \infty$ because the leftmost side diverges as $l_{\max} \to \infty$. To see this, we begin with the following Legendre expansion:

$$\frac{1}{\nu - \mu} = \sum_{l=0}^{\infty} c_l(\nu)P_l(\mu)\ , \qquad \nu > 1\,,\ \ -1 \le \mu \le 1\,, \tag{9.12}$$

where the first a few $c_l(\nu)$ are obtained as

$$c_0 = \frac{1}{2}\log\frac{\nu+1}{\nu-1}\,,\ \ c_1 = 3\nu\coth^{-1}(\nu) - 3\,,\ \ c_2 = \frac{5}{2}((3\nu^2-1)\coth^{-1}(\nu) - 3\nu)\,. \tag{9.13}$$

Now, suppose $|\mu| > 1$. Since Legendre polynomials are given by

$$P_0(\mu) = 1\,,\ \ P_1(\mu) = \mu\,,\ \ P_2(\mu) = \frac{1}{2}(3\mu^2 - 1)\,,\ \ \cdots\,, \tag{9.14}$$

the right-hand side of (9.12) diverges if $|\mu| > 1$ although the left-hand side remains finite. For simplicity, let us consider the case of isotropic scattering. For $\nu \notin [-1,1]$, we have

$$\mathcal{R}_{\hat{\mathbf{k}}(\nu,\mathbf{q})}\Phi_\nu^0(\hat{\mathbf{s}}) = \frac{\varpi\nu}{2}\frac{1}{\nu - \hat{\mathbf{s}}\cdot\hat{\mathbf{k}}(\nu,\mathbf{q})}\,. \tag{9.15}$$

Recall

$$\hat{\mathbf{s}}\cdot\hat{\mathbf{k}}(\nu,\mathbf{q}) = -i\nu q\sqrt{1-\mu^2}\cos(\varphi - \varphi_{\mathbf{q}}) + \mu\hat{k}_z(\nu q), \tag{9.16}$$

where $\hat{k}_z(\nu q)$ is given in (4.34). Since $\hat{k}_z(\nu q)$ is greater than 1, the Legendre expansion of $\mathcal{R}_{\hat{\mathbf{k}}(\nu,\mathbf{q})}\Phi_\nu^0(\hat{\mathbf{s}})$ does not converge. In the spherical-harmonic expansion for the three-dimensional RTE, usually $l_{\max}$ around 10 gives good numerical solutions. In Fig. 9.1, the hemispheric flux (7.33) is computed with the $F_N$ method (see (7.34)) and the spherical-harmonic expansion.

In one dimension, it is known that the discrete ordinates and spherical harmonics methods are equivalent (Barichello and Siewert, 1998).

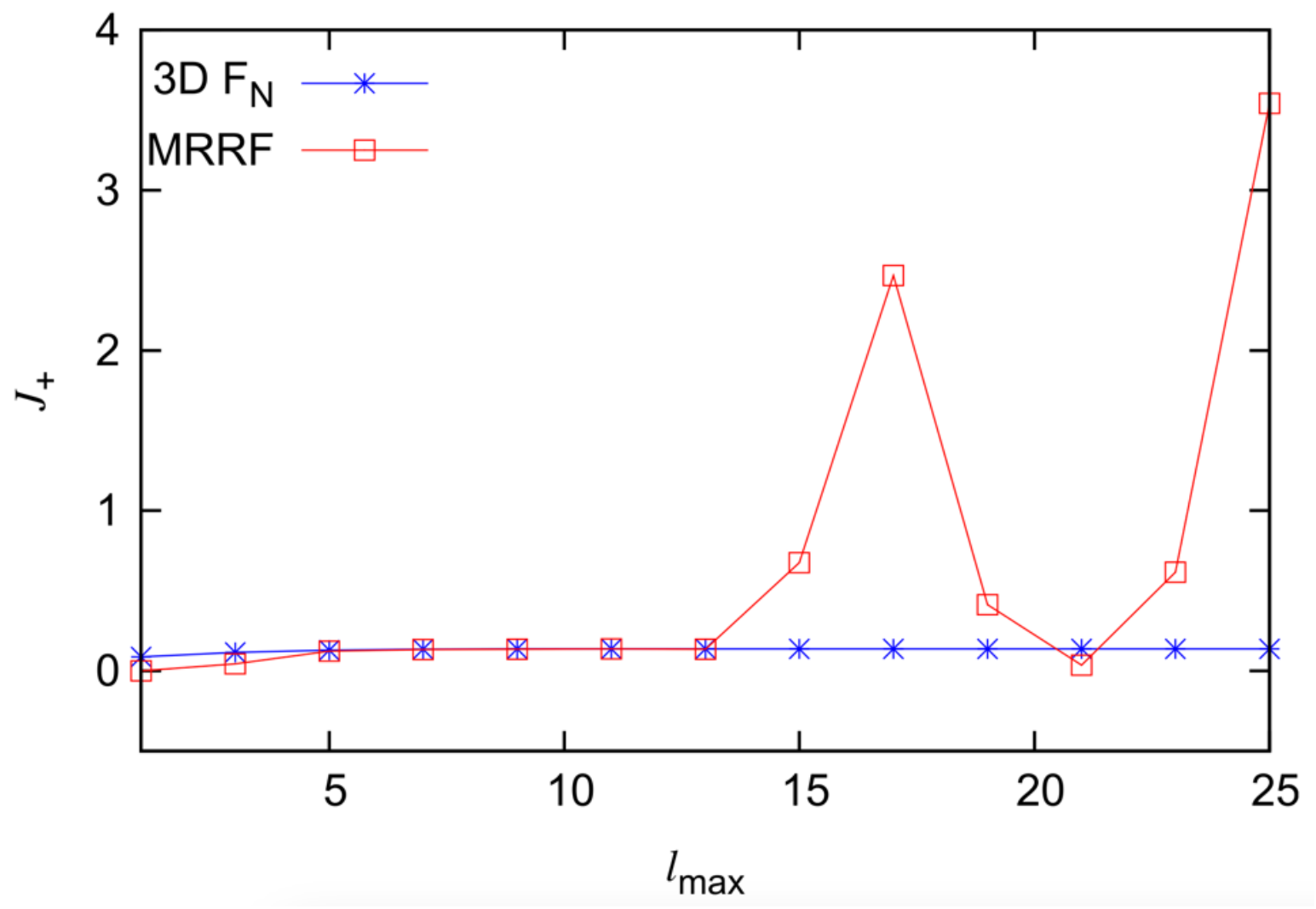

Figure 9.1. (Fig.1 in Machida, 2015) The exitance $|J_+(\boldsymbol{\rho};\mathbf{q}_0)|$, where $J_+(\boldsymbol{\rho};\mathbf{q}_0)$ given in (7.33), is plotted as a function of $l_{\max}$. The parameters are set to $\mu_a = 0.05$, $\mu_s = 100$, $\mathrm{g} = 0.001$. The unit of length is normalized using $\mu_t$. MRRF means the spherical-harmonic expansion in (9.4), (9.5).

### 9.3. Green's function for an infinite medium

Let us compute the Green's function for an infinite medium or the fundamental solution for the RTE (4.62). It is convenient to introduce

$$\sigma_l = \frac{\mu_t h_l}{2l+1}, \tag{9.17}$$

where $h_l$ was defined in (4.20). We also introduce

$$\psi_l^m(\lambda) = \sqrt{\frac{(2l+1)\sigma_l\lambda}{2\mathcal{N}^m(\nu)\mu_t}}\, g_l^m(\lambda). \tag{9.18}$$

Recall that (4.67) is the RTE in the Fourier space. Let us substitute the expression on the leftmost side of (9.7) for $\mathcal{R}_{\hat{\mathbf{k}}(\nu,\mathbf{q})}\Phi_\nu^m(\hat{\mathbf{s}})$ in the fundamental solution in (4.66), (4.71). Within the approximation of the spherical-harmonic expansion, we obtain

$$\begin{aligned}
&\tilde{G}(\mathbf{q};\ z,\hat{\mathbf{s}};\mathbf{r}_0,\hat{\mathbf{s}}_0) \\
&= \frac{1}{2\pi}e^{-i\mu_t\mathbf{q}\cdot\boldsymbol{\rho}_0}\sum_{m=-l_{\max}}^{l_{\max}}\sum_{\lambda>0}\frac{e^{-\mu_t\hat{k}_z(\lambda q)|z-z_0|/\lambda}}{\hat{k}_z(\lambda q)\mathcal{N}^m(\lambda)}\left[\sum_{l=|m|}^{l_{\max}}\sqrt{(2l+1)\pi}\,g_l^m(\lambda)\mathcal{R}_{\hat{\mathbf{k}}(\pm\lambda,\mathbf{q})}Y_{lm}(\hat{\mathbf{s}})\right] \\
&\qquad\times\left[\sum_{l'=|m|}^{l_{\max}}\sqrt{(2l'+1)\pi}\,g_{l'}^m(\lambda)\mathcal{R}_{\hat{\mathbf{k}}(\pm\lambda,\mathbf{q})}Y_{l'm}^*(\hat{\mathbf{s}}_0)\right].
\end{aligned} \tag{9.20}$$

Hence the Green's function is written as

$$G(\mathbf{r},\hat{\mathbf{s}};\mathbf{r}_0,\hat{\mathbf{s}}_0)=\frac{{\mu_t}^3}{(2\pi)^2}\int_{\mathbb{R}^2}\sum_{m=-l_{\max}}^{l_{\max}}\sum_{\lambda>0}\frac{1}{\lambda\hat{k}_z(\lambda q)}I_\lambda^{(\pm)m}(\mathbf{r},\hat{\mathbf{s}};\mathbf{q})I_\lambda^{(\mp)m}(\mathbf{r}_0,-\hat{\mathbf{s}}_0;-\mathbf{q})\,d\mathbf{q}\;. \quad (9.21)$$

where

$$I_\lambda^{(+)m}(\mathbf{r},\hat{\mathbf{s}};\mathbf{q})=e^{i\mu_t\mathbf{q}\cdot\boldsymbol{\rho}}e^{-\mu_t\hat{k}_z(\lambda q)z/\lambda}\sum_{l=|m|}^{l_{\max}}\frac{\psi_l^m(\lambda)}{\sqrt{\sigma_l}}\mathcal{R}_{\hat{\mathbf{k}}(\lambda,\mathbf{q})}Y_{lm}(\hat{\mathbf{s}})\;, \quad (9.22)$$

and

$$I_\lambda^{(-)m}(\mathbf{r},\hat{\mathbf{s}};\mathbf{q})=e^{i\mu_t\mathbf{q}\cdot\boldsymbol{\rho}}e^{\mu_t\hat{k}_z(\lambda q)z/\lambda}\sum_{l=|m|}^{l_{\max}}\frac{\psi_l^m(\lambda)}{\sqrt{\sigma_l}}\mathcal{R}_{\hat{\mathbf{k}}(\lambda,-\mathbf{q})}Y_{lm}^*(-\hat{\mathbf{s}})$$

$$=e^{i\mu_t\mathbf{q}\cdot\boldsymbol{\rho}}e^{\mu_t\hat{k}_z(\lambda q)z/\lambda}\sum_{l=|m|}^{l_{\max}}\frac{\psi_l^m(\lambda)}{\sqrt{\sigma_l}}(-1)^{l+m}\mathcal{R}_{\hat{\mathbf{k}}(\lambda,-\mathbf{q})}Y_{l,-m}(\hat{\mathbf{s}})$$

$$=(-1)^m e^{i\mu_t\mathbf{q}\cdot\boldsymbol{\rho}}e^{\mu_t\hat{k}_z(\lambda q)z/\lambda}\sum_{l=|m|}^{l_{\max}}\frac{\psi_l^m(\lambda)}{\sqrt{\sigma_l}}\sum_{m'=-l}^{l}(-1)^{l+m'}e^{-im'(\varphi_{\mathbf{q}}+\pi)}d_{m',-m}^{l}[i\tau(\lambda q)]Y_{lm'}(\hat{\mathbf{s}})\;. \quad (9.23)$$

Here, we used $Y_{lm}^*(\hat{\mathbf{s}})=(-1)^mY_{l,-m}(\hat{\mathbf{s}})=(-1)^{l+m}Y_{l,-m}(-\hat{\mathbf{s}})$.

As follows, we can see that $\psi_l^m$ are normalized as

$$\sum_{l=|m|}^{l_{\max}}|\psi_l^m(\lambda)|^2=1\,. \quad (9.24)$$

To see this, first we note that the expressions of the fundamental solution imply

$$\frac{1}{2\pi\mathcal{N}^m(\lambda)}\mathcal{R}_{\hat{\mathbf{k}}(\lambda,\mathbf{q})}\Phi_\lambda^m(\hat{\mathbf{s}})\Phi_\lambda^{m*}(\hat{\mathbf{s}}_0)$$

$$=\frac{\mu_t}{\lambda}\left(\sum_{l=|m|}^{l_{\max}}\frac{\psi_l^m(\lambda)}{\sqrt{\sigma_l}}\mathcal{R}_{\hat{\mathbf{k}}(\lambda,\mathbf{q})}Y_{lm}(\hat{\mathbf{s}})\right)\left(\sum_{l=|m|}^{l_{\max}}\frac{\psi_l^m(\lambda)}{\sqrt{\sigma_l}}\mathcal{R}_{\hat{\mathbf{k}}(\lambda,\mathbf{q})}Y_{lm}^*(\hat{\mathbf{s}}_0)\right). \quad (9.25)$$

Hence, by setting $\hat{\mathbf{s}}_0=\hat{\mathbf{s}}$,

$$|\Phi_\lambda^m(\hat{\mathbf{s}})|^2=\frac{2\pi\mathcal{N}^m(\lambda)\mu_t}{\lambda}\sum_{l=|m|}^{l_{\max}}\sum_{l'=|m|}^{l_{\max}}\frac{\psi_l^m(\lambda)\psi_{l'}^m(\lambda)}{\sqrt{\sigma_l\sigma_{l'}}}Y_{lm}(\hat{\mathbf{s}})Y_{l'm}^*(\hat{\mathbf{s}})\;. \quad (9.26)$$

We have

$$2\pi\mathcal{N}^m(\lambda)=\int_{\mathbb{S}^2}\mu|\Phi_\lambda^m(\hat{\mathbf{s}})|^2\,d\hat{\mathbf{s}}$$

$$=\frac{2\pi\mathcal{N}^m(\lambda)\mu_t}{\lambda}\sum_{l=|m|}^{l_{\max}}\sum_{l'=|m|}^{l_{\max}}\frac{\psi_l^m(\lambda)\psi_{l'}^m(\lambda)}{\sqrt{\sigma_l\sigma_{l'}}}\int_{\mathbb{S}^2}\mu Y_{lm}(\hat{\mathbf{s}})Y_{l'm}^*(\hat{\mathbf{s}})\,d\hat{\mathbf{s}}\;. \quad (9.27)$$

We note that

$$\int_{\mathbb{S}^2}\mu Y_{lm}^*(\hat{\mathbf{s}})Y_{l'm'}(\hat{\mathbf{s}})d\hat{\mathbf{s}}=\frac{1}{2}\delta_{mm'}\sqrt{(2l+1)(2l'+1)\frac{(l-m)!}{(l+m)!}\frac{(l'-m)!}{(l'+m)!}}\int_{-1}^{1}\mu P_l^m(\mu)P_{l'}^m(\mu)\,d\mu$$

$$= \delta_{mm'}\left(\sqrt{\frac{(l+1)^2-m^2}{4(l+1)^2-1}}\delta_{l+1,l'} + \sqrt{\frac{l^2-m^2}{4l^2-1}}\delta_{l-1,l'}\right). \tag{9.28}$$

We obtain

$$1 = \frac{\mu_t}{\lambda}\sum_{l=|m|}^{l_{\max}} \frac{\psi_l^m(\nu)}{\sqrt{\sigma_l}}\left(\frac{\psi_{l+1}^m(\lambda)}{\sqrt{\sigma_{l+1}}}\sqrt{\frac{(l+1)^2-m^2}{4(l+1)^2-1}} + \frac{\psi_{l-1}^m(\lambda)}{\sqrt{\sigma_{l-1}}}\sqrt{\frac{l^2-m^2}{4l^2-1}}\right)$$

$$= \frac{1}{2\mathcal{N}^m(\lambda)}\sum_{l=0}^{l_{\max}} \sqrt{2l+1}g_l^m(\lambda)$$

$$\times\left(\sqrt{2l+3}g_{l+1}^m(\lambda)\sqrt{\frac{(l+1)^2-m^2}{4(l+1)^2-1}} + \sqrt{2l-1}g_{l-1}^m(\lambda)\sqrt{\frac{l^2-m^2}{4l^2-1}}\right). \tag{9.29}$$

We have (see (4.19) and (4.61))

$$b_{l+1}(m)\sqrt{2(l+1)+1}g_{l+1}^m(\lambda) + b_l(m)\sqrt{2(l-1)+1}g_{l-1}^m(\lambda)$$

$$= \frac{\lambda h_l}{2l+1}\sqrt{2l+1}g_l^m(\lambda), \tag{9.30}$$

where

$$b_l(m) = \sqrt{\frac{l^2-m^2}{4l^2-1}}. \tag{9.31}$$

Hence we obtain

$$1 = \frac{\lambda}{2\mathcal{N}^m(\lambda)}\sum_{l=|m|}^{l_{\max}} h_l g_l^m(\lambda)^2 . \tag{9.32}$$

Therefore,

$$\sum_{l=|m|}^{l_{\max}} |\psi_l^m(\lambda)|^2 = 1. \tag{9.33}$$

The relation (9.7) implies

$$\sum_{l'=0}^{\infty}\sum_{m'=-l'}^{l'} c_{l'm'}^m(\lambda)Y_{l'm'}(\hat{\mathbf{s}}) = \Phi_\lambda^m(\hat{\mathbf{s}}). \tag{9.34}$$

Then using $\lambda$ instead of $\nu$, (4.17) can be rewritten as

$$\left(1-\frac{\mu}{\lambda}\right)\sum_{l'=0}^{\infty}\sum_{m'=-l'}^{l'} c_{l'm'}^m(\lambda)Y_{l'm'}(\hat{\mathbf{s}}) = \varpi\sum_{l'=0}^{\infty}\sum_{m'=-l'}^{l'}\frac{\beta_{l'}}{2l'+1}c_{l'm'}^m(\lambda)Y_{l'm'}(\hat{\mathbf{s}}) . \tag{9.35}$$

Thus for $|M| \le l_{\max}$, we have

$$c_{lm}^M(\lambda) - \frac{1}{\lambda}\sum_{l'=0}^{\infty}\sum_{m'=-l'}^{l'} c_{l'm'}^M(\lambda)\int_{\mathbb{S}^2}\mu Y_{l'm'}(\hat{\mathbf{s}})Y_{lm}^*(\hat{\mathbf{s}})\,d\hat{\mathbf{s}} = \frac{\varpi\beta_l}{2l+1}c_{lm}^M(\lambda). \tag{9.36}$$

The right-hand side of (9.36) vanishes if $l > l_{\max}$ (Recall that $l_{\max}$ was introduced in (4.1)). Using the orthogonality relation for associated Legendre polynomials (4.8), we obtain (Panasyuk, Schotland, and Markel, 2006)

$$\sum_{l'=|M|}^{l_{\max}} B_{ll'}(M)\sqrt{\sigma_{l'}}c_{l'M}^{M}(\lambda) = \frac{\lambda}{\mu_t}\sqrt{\sigma_l}c_{lM}^{M}(\lambda)\,, \tag{9.37}$$

where $B_{ll'}(M)$ are elements of matrix $B(M)$ which are given by

$$B_{ll'}(M) = \beta_l(M)\delta_{l',l-1} + \beta_{l+1}(M)\delta_{l',l+1}\,, \tag{9.38}$$

where

$$\beta_l(M) = \sqrt{\frac{l^2 - M^2}{(4l^2-1)\sigma_{l-1}\sigma_l}}\,. \tag{9.39}$$

We can readily see that (9.36) is equivalent to (4.61) and we have the relation

$$B(M) = \frac{1}{\mu_t}A(M)\,, \tag{9.40}$$

where $A(M)$ was introduced in (4.50). Hence, $\lambda$ is obtained using the tridiagonal matrices.

Let us look at the fundamental solution more carefully. The expression (9.21) can be rewritten as

$$G(\mathbf{r},\hat{\mathbf{s}};\mathbf{r}_0,\hat{\mathbf{s}}_0) = \frac{1}{(2\pi)^2}\sum_{l=0}^{l_{\max}}\sum_{m=-l}^{l}\sum_{l'=0}^{l_{\max}}\sum_{m'=-l'}^{l'}$$

$$\int_{\mathbb{R}^2} e^{i\mu_t\mathbf{q}\cdot(\boldsymbol{\rho}-\boldsymbol{\rho}_0)}Y_{lm}(\hat{\mathbf{s}})\kappa_{l'm'}^{lm}(\mathbf{q},z-z_0)Y_{l'm'}^{*}(\hat{\mathbf{s}}_0)\,d\mathbf{q}\ . \tag{9.41}$$

where

$$\kappa_{l'm'}^{lm}(\mathbf{q},z-z_0) = \frac{e^{-i(m-m')\varphi_{\mathbf{q}}}}{\sqrt{\sigma_l\sigma_{l'}}}[\mathrm{sgn}(z-z_0)]^{l+l'+m+m'}\sum_{m_1=-l}^{l}\sum_{m_2=-l'}^{l'}\sum_{\lambda>0} d_{mm_1}^{l}[i\tau(\lambda q)]$$

$$\times\,\psi_l^{m_1}(\lambda)\frac{e^{-\mu_t\hat{k}_z(\nu q)|z-z_0|/\lambda}}{\lambda\hat{k}_z(\lambda q)}\psi_{l'}^{m_2}(\lambda)^{*}d_{m'm_2}^{l'}[i\tau(\lambda q)] \tag{9.42}$$

An alternative way of computing the fundamental solution is to express it as

$$G(\mathbf{r},\hat{\mathbf{s}};\mathbf{r}_0,\hat{\mathbf{s}}_0) = \frac{1}{(2\pi)^3}\sum_{l=0}^{l_{\max}}\sum_{m=-l}^{l}\sum_{l'=0}^{l_{\max}}\sum_{m'=-l'}^{l'}\int_0^{\infty}k^2$$

$$\times\int_{\mathbb{R}^2} e^{ik\hat{\mathbf{k}}\cdot(\mathbf{r}-\mathbf{r}_0)}\left(\mathcal{R}_{\hat{\mathbf{k}}}Y_{lm}(\hat{\mathbf{s}})\right)K_{l'm'}^{lm}(k)\left(\mathcal{R}_{\hat{\mathbf{k}}}Y_{l'm'}^{*}(\hat{\mathbf{s}}_0)\right)\,d\hat{\mathbf{k}}\,dk\,, \tag{9.43}$$

where $K_{l'm'}^{lm}(k)$ will be derived below. The reciprocal property $G(\mathbf{r},\hat{\mathbf{s}};\mathbf{r}_0,\hat{\mathbf{s}}_0) = G(\mathbf{r}_0,-\hat{\mathbf{s}}_0;\mathbf{r},-\hat{\mathbf{s}})$ implies

$$K_{lm}^{l'm'}(k) = (-1)^{l+l'}K_{l'm'}^{lm}(k)^{*}\,. \tag{9.44}$$

By substituting the above expression into the RTE (4.62), we obtain

$$\sum_{l''=0}^{l_{\max}}\sum_{m''=-l''}^{l''}\left(ikR_{l''m''}^{lm} + \Sigma_{l''m''}^{lm}\right)K_{l'm'}^{l''m''}(k) = \delta_{ll'}\delta_{mm'}\,, \tag{9.45}$$

where

$$R_{l''m''}^{lm} = \int_{\mathbb{S}^2}\hat{\mathbf{s}}\cdot\hat{\mathbf{k}}\left(\mathcal{R}_{\hat{\mathbf{k}}}Y_{lm}^{*}(\hat{\mathbf{s}})\right)\left(\mathcal{R}_{\hat{\mathbf{k}}}Y_{l'm'}(\hat{\mathbf{s}})\right)\,d\hat{\mathbf{s}} = \int_{\mathbb{S}^2}\mu Y_{lm}^{*}(\hat{\mathbf{s}})Y_{l'm'}(\hat{\mathbf{s}})\,d\hat{\mathbf{s}}\ , \tag{9.46}$$

and

$$\Sigma_{l''m''}^{lm} = \sigma_l \delta_{ll''} \delta_{mm''} \, . \tag{9.47}$$

The relation can be formally expressed as

$$K = \frac{1}{\sqrt{\Sigma}} (1 + ikW)^{-1} \frac{1}{\sqrt{\Sigma}} \, , \tag{9.48}$$

where

$$W = \frac{1}{\sqrt{\Sigma}} R \frac{1}{\sqrt{\Sigma}} \, . \tag{9.49}$$

We note that the matrix $W$ has a block-diagonal structure:

$$W_{l'm'}^{lm} = \frac{1}{\sqrt{\sigma_l}} R_{l'm'}^{lm} \frac{1}{\sqrt{\sigma_{l'}}} = B_{ll'}(m) \delta_{mm'} \, . \tag{9.50}$$

Recall that $\lambda/\mu_t$ and $\psi_l^m(\lambda)$ are eigenvalues and eigenvectors of $B$. We obtain

$$K_{l'm'}^{lm}(k) = \sum_{m=-l_{\max}}^{l_{\max}} \sum_{\lambda} \sum_{l=|m|}^{l_{\max}} \sum_{l'=|m|}^{l_{\max}} \frac{\psi_l^m(\lambda) \psi_{l'}^m(\lambda)^*}{\sqrt{\sigma_l} \left(1 + \frac{ik\lambda}{\mu_t}\right) \sqrt{\sigma_{l'}}} \, . \tag{9.51}$$

For the relation between $\kappa_{l'm'}^{lm}$ and $K_{l'm'}^{lm}$, we refer the reader to Panasyuk, Schotland, and Markel, 2006.

We will use the notation

$$\mathbf{R} = \mathbf{r} - \mathbf{r}_0 \, , \;\; R = |\mathbf{r} - \mathbf{r}_0| \, , \;\; \hat{\mathbf{R}} = \frac{\mathbf{r} - \mathbf{r}_0}{|\mathbf{r} - \mathbf{r}_0|} \, . \tag{9.52}$$

Let us write (9.43) as

$$G(\mathbf{r}, \hat{\mathbf{s}}; \mathbf{r}_0, \hat{\mathbf{s}}_0) = \sum_{l=0}^{l_{\max}} \sum_{m=-l}^{l} \sum_{l'=0}^{l_{\max}} \sum_{m'=-l'}^{l'} Y_{lm}(\hat{\mathbf{s}}) \chi_{l'm'}^{lm}(\mathbf{R}) Y_{l'm'}^*(\hat{\mathbf{s}}_0) \tag{9.53}$$

with function $\chi_{l'm'}^{lm}(\mathbf{R})$. In the case of $\hat{\mathbf{R}} = \hat{\mathbf{z}}$, we obtain

$$\chi_{l'm'}^{lm}(\mathbf{R}) = \frac{{\mu_t}^3}{2\pi\sqrt{\sigma_l \sigma_{l'}}} \delta_{mm'} \sum_{M=-\min(l,l')}^{\min(l,l')} (-1)^M$$

$$\times \sum_{j=0}^{\min(l,l')} C_{l,M,l',-M}^{|l-l'|-2j,0} C_{l,m,l',-m}^{|l-l'|+2j,0} \sum_{\lambda>0} \frac{\psi_l^M(\lambda) \psi_{l'}^M(\lambda)^*}{\lambda^3} k_{|l-l'|+2j} \left(\frac{\mu_t R}{\lambda}\right) , \tag{9.54}$$

where $k_n$ is the modified spherical Bessel function of the second kind (defined without the factor of $\pi/2$), $C_{j_1 m_1 j_2 m_2}^{j_3 m_3}$ are the Clebsch-Gordan coefficients. To show this, we set $\mathbf{r}_0 = \mathbf{0}$ and introduce the notation

$$\hat{\mathbf{r}} = \frac{\mathbf{r}}{|\mathbf{r}|} \, . \tag{9.55}$$

In general, we have

$$\chi_{l'm'}^{lm}(\mathbf{r}) = \frac{1}{(2\pi)^3 \sqrt{\sigma_l \sigma_{l'}}} \sum_{m_1=-l}^{l} \sum_{m_2=-l'}^{l'} \int_0^\infty k^2 \int_{\mathbb{R}^2} e^{ik\hat{\mathbf{k}}\cdot\mathbf{r}} e^{-i(m-m')\varphi_{\hat{\mathbf{k}}}} d_{mm_1}^l(\theta_{\hat{\mathbf{k}}}) d_{m'm_2}^{l'}(\theta_{\hat{\mathbf{k}}}) \, d\hat{\mathbf{k}}$$

$$\times \sum_{\lambda} \frac{\psi_l^{m_1}(\lambda) \psi_{l'}^{m_2}(\lambda)^*}{1 + ik\lambda/\mu_t} \, dk \, . \tag{9.56}$$

Since $\hat{\mathbf{r}} = \hat{\mathbf{z}}$, we have $\hat{\mathbf{k}} \cdot \hat{\mathbf{r}} = \cos\theta_{\hat{\mathbf{k}}}$. In this case,

$$\chi_{l'm'}^{lm}(\mathbf{r}) = \frac{1}{(2\pi)^2\sqrt{\sigma_l\sigma_{l'}}} \sum_{M=-\min(l,l')}^{\min(l,l')} \int_0^\infty k^2 \sum_\lambda \frac{\psi_l^M(\lambda)\psi_{l'}^M(\lambda)^*}{1+ik\lambda/\mu_t} I\, dk , \tag{9.57}$$

where

$$I = \frac{1}{2\pi}\int_0^{2\pi}\int_0^\pi e^{-i(m-m')\varphi} e^{ikr\cos\theta} d_{mM}^l(\theta) d_{m'M}^{l'}(\theta) \sin\theta\, d\theta d\varphi \ . \tag{9.58}$$

Using formulae,

$$e^{ikr\cos\theta} = \sum_{L=0}^\infty i^L(2L+1) j_L(kr) d_{00}^L(\theta) \ , \tag{9.59}$$

where $j_L$ is the spherical Bessel function of the first kind, and

$$\int_0^\pi d_{mM}^l(\theta) d_{mM}^{l'}(\theta) d_{00}^L(\theta) \sin\theta\, d\theta = \frac{2(-1)^{m-M}}{2L+1} C_{l,m,l',-m}^{L,0} C_{l,M,l',-M}^{L,0} , \tag{9.60}$$

we obtain

$$I = 2\delta_{mm'}(-1)^{m-M} \sum_{L=|l-l'|}^{l+l'} i^L j_L(kr) C_{l,m,l',-m}^{L,0} C_{l,M,l',-M}^{L,0} \ . \tag{9.61}$$

In this way, we arrive at

$$\chi_{l'm'}^{lm}(\mathbf{r}) = \frac{2(-1)^m}{(2\pi)^2\sqrt{\sigma_l\sigma_{l'}}}\delta_{mm'} \sum_{M=-\min(l,l')}^{\min(l,l')} (-1)^M \sum_{L=|l-l'|}^{l+l'} i^L C_{l,m,l',-m}^{L,0} C_{l,M,l',-M}^{L,0}$$

$$\times \sum_\lambda \psi_l^M(\lambda)\psi_{l'}^M(\lambda)^* \int_0^\infty k^2 \frac{j_L(kr)}{1+ik\lambda/\mu_t}\, dk$$

$$= \frac{2(-1)^m}{\sqrt{\sigma_l\sigma_{l'}}}\delta_{mm'} \sum_{M=-\min(l,l')}^{\min(l,l')} (-1)^M \sum_{j=0}^{\min(l,l')} i^{|l-l'|+2j} C_{l,m,l',-m}^{|l-l'|+2j,0} C_{l,M,l',-M}^{|l-l'|+2j,0}$$

$$\times \sum_{\lambda>0} \psi_l^M(\lambda)\psi_{l'}^M(\lambda)^* J \ , \tag{9.62}$$

where

$$J = \frac{1}{(2\pi)^2}\int_0^\infty k^2 j_{|l-l'|+2j}(kr) \frac{1+(-1)^{l+l'} - ik(\lambda/\mu_t)\left[1-(-1)^{l+l'}\right]}{1+(k\lambda/\mu_t)^2}\, dk \ . \tag{9.63}$$

We can perform the integral as

$$J = \pi i^{-(|l-l'|+2j)} \left(\frac{\mu_t}{\lambda}\right)^3 k_{|l-l'|+2j}\left(\frac{\mu_t r}{\lambda}\right) . \tag{9.64}$$

Thus, $\chi_{l'm'}^{lm}(\mathbf{r})$ is obtained. In the case of $\hat{\mathbf{s}}_0 = \hat{\mathbf{z}}$, we obtain

$$\chi_{l'm'}^{lm}(\mathbf{r}-\mathbf{r}_0)$$

$$= \frac{{\mu_t}^3(-1)^l}{\sqrt{\pi(2l'+1)\sigma_l\sigma_{l'}}}\delta_{m'0} \sum_{M=-\min(l,l')}^{\min(l,l')} \sum_{j=0}^{\min(l,l')} \left(\mathcal{R}_{\hat{\mathbf{s}}_0} Y_{|l-l'|+2j,m}^*\left(\frac{\mathbf{r}-\mathbf{r}_0}{|\mathbf{r}-\mathbf{r}_0|}\right)\right) C_{l,M,|l-l'|+2j,0}^{l',M} C_{l,m,l',0}^{|l-l'|+2j,m}$$

$$\times \sum_{\lambda>0} \frac{\psi_l^M(\lambda)\psi_{l'}^M(\lambda)^*}{\lambda^3} k_{|l-l'|+2j}\left(\frac{\mu_t|\mathbf{r}-\mathbf{r}_0|}{\lambda}\right). \tag{9.65}$$

See Panasyuk, Schotland, and Markel, 2006 for detail calculations.

### 9.4. Ballistic subtraction

As seen in Secs. 7 and 8, numerical instability can be reduced by the ballistic subtraction. We will explore an alternative way of the ballistic subtraction (Panasyuk, Schotland, and Markel, 2006).

One way to subtract the ballistic term is to decompose the specific intensity into two terms (see Sec. 8.2 and also Liemert and Kienle, 2015). The other way is to subtract the ballistic component term by term.

Let us write

$$G_b(\mathbf{r},\hat{\mathbf{s}};\mathbf{r}_0,\hat{\mathbf{s}}_0) = \sum_{l=0}^{l_{\max}} \sum_{m=-l}^{l} \sum_{l'=0}^{l_{\max}} \sum_{m'=-l'}^{l'} Y_{lm}(\hat{\mathbf{s}})\chi_{b,l'm'}^{lm}(\mathbf{R})Y_{l'm'}^*(\hat{\mathbf{s}}_0). \tag{9.66}$$

In the case of $\hat{\mathbf{s}}_0 = \hat{\mathbf{z}}$, we obtain

$$\chi_{b,l'm'}^{lm}(\mathbf{r}-\mathbf{r}_0) = \frac{\sqrt{(2l+1)(2l'+1)}}{4\pi|\mathbf{r}-\mathbf{r}_0|^2}\delta_{m0}e^{-\mu_t|\mathbf{r}-\mathbf{r}_0|}. \tag{9.67}$$

Therefore the ballistic term can be subtracted if $\chi_{l'm'}^{lm}$ are replaced by $\chi_{l'm'}^{lm} - \chi_{b,l'm'}^{lm}$ in the expression of $G(\mathbf{r},\hat{\mathbf{s}};\mathbf{r}_0,\hat{\mathbf{s}}_0)$ in (9.53). Since the expansion by spherical harmonics diverges in rotated reference frames as described in Sec. 9.2, the two ballistic subtractions are not equivalent.

We note that the expression of $G(\mathbf{r},\hat{\mathbf{s}};\mathbf{r}_0,\hat{\mathbf{s}}_0)$ in (9.53) with $\chi_{l'm'}^{lm}(\mathbf{r}-\mathbf{r}_0)$ in (9.65) is numerically more stable than the naive superposition of eigenmodes given in (9.21). The reason is that in (9.65) reference frames are rotated about the unit vector $\hat{\mathbf{s}}_0$ instead of $\hat{\mathbf{k}}(\lambda,\mathbf{q})$, for which analytically-continued Wigner's $d$-matrices appear when reference frames are rotated. The use of the function $\chi_{l'm'}^{lm}(\mathbf{r}-\mathbf{r}_0)$ is possible for the whole space. When the solution of the RTE is expressed as a superposition of eigenmodes in the presence of boundaries, the expansion by $\mathcal{R}_{\hat{\mathbf{k}}}Y_{lm}(\hat{\mathbf{s}})$ gets unstable when $l_{\max}$ becomes large.

We close this subsection by numerical examples. In Fig. 9.2, the computed specific intensity in (9.53) is shown for different $l_{\max}$. In addition, the effect of the ballistic subtraction is illustrated.

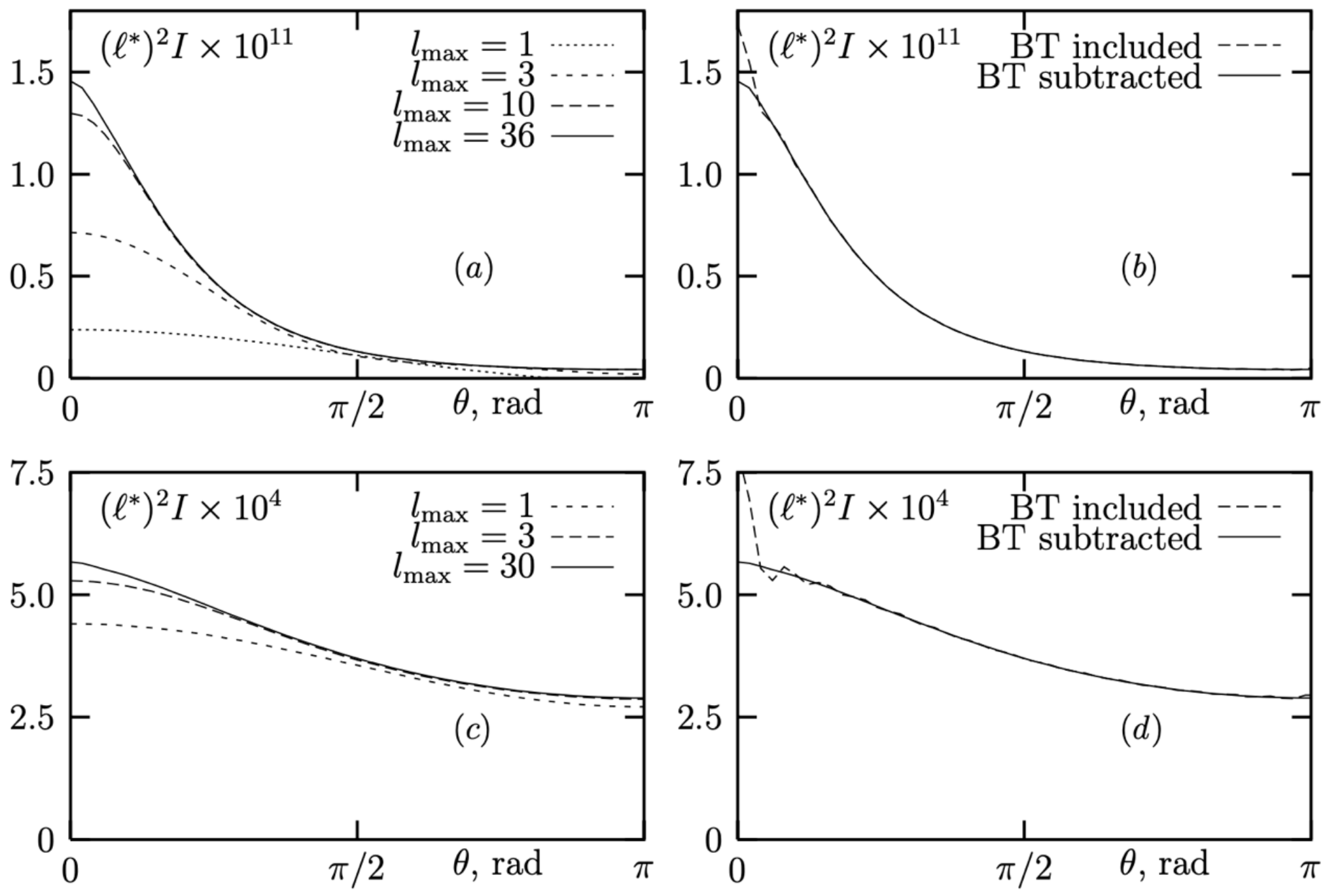

Figure 9.2. (Fig. 2 in Panasyuk, Schotland, and Markel, 2006). Angular dependence of the specific intensity (9.53) for forward propagation at the distance $z$ from the source. Left column (a), (c): convergence with parameter $l_{\max}$. Right column (b), (d): solid line shows the converged result obtained with subtraction of the ballistic term; dashed line: result obtained without subtraction of the ballistic term for the same $l_{\max}$. Top row (a), (b): $\mathrm{g} = 0.5$, $\mu_a/\mu_s = 0.5$, $z = 20\ell^*$ ($\ell^* = \mu_a + (1-\mathrm{g})\mu_s$). Bottom row (c), (d): $\mathrm{g} = 0.2$, $\mu_a/\mu_s = 0.01$, $z = 10\ell^*$.

### 9.5. Green's function for the slab geometry

Let us consider the RTE in the presence of boundaries. Here we compute the Green's function $I(\mathbf{r}, \hat{\mathbf{s}}; \mathbf{r}_0, \hat{\mathbf{s}}_0)$ for a slab of width $L$. The RTE is written as

$$\hat{\mathbf{s}} \cdot \nabla I(\mathbf{r}, \hat{\mathbf{s}}; \mathbf{r}_0, \hat{\mathbf{s}}_0) + \mu_t I(\mathbf{r}, \hat{\mathbf{s}}; \mathbf{r}_0, \hat{\mathbf{s}}_0) = \mu_s \int_{\mathbb{S}^2} p(\hat{\mathbf{s}}, \hat{\mathbf{s}}') I(\mathbf{r}, \hat{\mathbf{s}}'; \mathbf{r}_0, \hat{\mathbf{s}}_0) d\hat{\mathbf{s}}' \ , \quad 0 < z < L \, . \tag{9.68}$$

The boundary conditions are given by

$$I(\mathbf{r}, \hat{\mathbf{s}}; \mathbf{r}_0, \hat{\mathbf{s}}_0) = I_0 \delta(\boldsymbol{\rho} - \boldsymbol{\rho}_0) \delta(\hat{\mathbf{s}} - \hat{\mathbf{s}}_0) \, , \quad z = 0 \, , \quad \hat{\mathbf{s}} \cdot \hat{\mathbf{z}} > 0 \, , \tag{9.69}$$

and

$$I(\mathbf{r}, \hat{\mathbf{s}}; \mathbf{r}_0, \hat{\mathbf{s}}_0) = 0 \, , \quad z = L \, , \quad \hat{\mathbf{s}} \cdot \hat{\mathbf{z}} < 0 \, . \tag{9.70}$$

The solution $I(\mathbf{r}, \hat{\mathbf{s}}; \mathbf{r}_0, \hat{\mathbf{s}}_0)$ can be written as a superposition of the eigenmodes (9.22), (9.23):

$$I(\mathbf{r}, \hat{\mathbf{s}}; \mathbf{r}_0, \hat{\mathbf{s}}_0) = \frac{\mu_t^2}{(2\pi)^2} \sum_{m=-l_{\max}}^{l_{\max}} \sum_{\lambda}$$

$$\int_{\mathbb{R}^2}\left[F_\lambda^{(+)m}(\mathbf{q})I_\lambda^{(+)m}(\mathbf{r},\hat{\mathbf{s}};\mathbf{q})+(-1)^mF_\lambda^{(-)m}(\mathbf{q})I_\lambda^{(-)m}(\mathbf{r},\hat{\mathbf{s}};\mathbf{q})\right]d\mathbf{q}\ , \tag{9.71}$$

where coefficients $F_\nu^{(\pm)m}$ will be determined so that $I(\mathbf{r},\hat{\mathbf{s}};\mathbf{r}_0,\hat{\mathbf{s}}_0)$ satisfies the boundary conditions.

By approximating the delta function for angles using spherical harmonics, we rewrite the boundary conditions as

$$I(\mathbf{r},\hat{\mathbf{s}};\mathbf{r}_0,\hat{\mathbf{s}}_0)=\Theta(\hat{\mathbf{s}}\cdot\hat{\mathbf{z}})I_0\delta(\boldsymbol{\rho}-\boldsymbol{\rho}_0)\sum_{l=0}^{l_{\max}}\sum_{m=-l}^{l}Y_{lm}(\hat{\mathbf{s}})Y_{lm}^*(\hat{\mathbf{s}}_0)$$
$$+\Theta(-\hat{\mathbf{s}}\cdot\hat{\mathbf{z}})I(\mathbf{r},\hat{\mathbf{s}};\mathbf{r}_0,\hat{\mathbf{s}}_0),\quad z=0\ , \tag{9.72}$$
$$I(\mathbf{r},\hat{\mathbf{s}};\mathbf{r}_0,\hat{\mathbf{s}}_0)=\Theta(\hat{\mathbf{s}}\cdot\hat{\mathbf{z}})I(\mathbf{r},\hat{\mathbf{s}};\mathbf{r}_0,\hat{\mathbf{s}}_0),\quad z=L \tag{9.73}$$

for $\hat{\mathbf{s}}\in\mathbb{S}^2$. We multiply $Y_{l'm'}^*(\hat{\mathbf{s}})e^{-i\mu_t\mathbf{q}\cdot\boldsymbol{\rho}}$ on both sides of the boundary conditions and integrate both sides over $\hat{\mathbf{s}}$ and $\boldsymbol{\rho}$. We obtain

$$\sum_{m=-l_{\max}}^{l_{\max}}\sum_{\lambda}\sum_{l=|m|}^{l_{\max}}\mathcal{B}_{l'l}^{m'}\frac{\psi_l^m(\lambda)}{\sqrt{\sigma_l}}e^{-im'(\varphi_{\mathbf{q}}+\pi)}\Big[F_\lambda^{(+)m}(\mathbf{q})d_{m',m}^l[i\tau(\lambda q)]$$
$$+F_\lambda^{(-)m}(\mathbf{q})(-1)^{l+m'}d_{m',-m}^l[i\tau(\lambda q)]\Big]=I_0e^{-i\mu_t\mathbf{q}\cdot\boldsymbol{\rho}_0}\sum_{l=|m|}^{l_{\max}}\mathcal{B}_{l'l}^{m'}Y_{lm'}^*(\hat{\mathbf{s}}_0)\ , \tag{9.74}$$

$$\sum_{m=-l_{\max}}^{l_{\max}}\sum_{\lambda}\sum_{l=|m|}^{l_{\max}}\mathcal{B}_{l'l}^{m'}(-1)^{l+l'}\frac{\psi_l^m(\lambda)}{\sqrt{\sigma_l}}e^{-im'(\varphi_{\mathbf{q}}+\pi)}\Big[F_\lambda^{(+)m}(\mathbf{q})e^{-\mu_t\hat{k}_z(\lambda q)L/\lambda}d_{m',m}^l[i\tau(\lambda q)]$$
$$+F_\lambda^{(-)m}(\mathbf{q})e^{\mu_t\hat{k}_z(\lambda q)L/\lambda}(-1)^{l+m'}d_{m',-m}^l[i\tau(\lambda q)]\Big]=0\ . \tag{9.75}$$

Here, $\mathcal{B}_{ll'}^m$ are defined by

$$\mathcal{B}_{ll'}^m=\frac{1}{2}\sqrt{\frac{(2l+1)(2l'+1)(l-m)!\,(l'-m)!}{(l+m)!\,(l'+m)!}}\int_0^1P_l^m(\mu)P_{l'}^m(\mu)\,d\mu\ . \tag{9.76}$$

Note that

$$\int_0^{2\pi}\int_0^{\pi/2}Y_{lm}(\hat{\mathbf{s}})Y_{l'm'}^*(\hat{\mathbf{s}})\sin\theta\ d\theta d\varphi=\delta_{mm'}\mathcal{B}_{ll'}^m\ , \tag{9.77}$$

$$\int_0^{2\pi}\int_{\pi/2}^{\pi}Y_{lm}(\hat{\mathbf{s}})Y_{l'm'}^*(\hat{\mathbf{s}})\sin\theta\ d\theta d\varphi=\delta_{mm'}(-1)^{l+l'}\mathcal{B}_{ll'}^m\ , \tag{9.78}$$

and $\mathcal{B}_{ll'}^m=\mathcal{B}_{ll'}^{-m}$. Moreover, $\mathcal{B}_{ll'}^m=\frac{1}{2}\delta_{ll'}$ if $l$ and $l'$ have the same parity. The coefficients $F_\nu^{(\pm)m}(\mathbf{q})$ can be obtained using (9.74) and (9.75).

In what follows, we restrict ourselves to the case of $\hat{\mathbf{s}}_0=\hat{\mathbf{z}}$. Let us write $F_\lambda^{(\pm)m}(\mathbf{q})$ as

$$F_\lambda^{(+)m}(\mathbf{q})=I_0f_\lambda^{(+)m}(q)e^{-i\mu_t\mathbf{q}\cdot\boldsymbol{\rho}_0},\quad F_\lambda^{(-)m}(\mathbf{q})=I_0f_\lambda^{(-)m}(q)e^{-i\mu_t\mathbf{q}\cdot\boldsymbol{\rho}_0}e^{-\mu_t\hat{k}_z(\nu q)L/\lambda}\ . \tag{9.79}$$

Then we have

$$\sum_{m=-l_{\max}}^{l_{\max}}\sum_{\lambda}\sum_{l=|m|}^{l_{\max}}\mathcal{B}_{l'l}^{m'}\frac{\psi_l^m(\lambda)}{\sqrt{\sigma_l}}\Big[f_\lambda^{(+)m}(q)d_{m',m}^l[i\tau(\lambda q)]$$

$$+f_v^{(-)m}(q)e^{-\mu_t\hat{k}_z(\lambda q)L/\lambda}(-1)^{l+m'}d^l_{m',-m}[i\tau(\lambda q)]\Big] = \delta_{m'0}\sum_{l=|m|}^{l_{\max}}\mathcal{B}^0_{l'l}\sqrt{\frac{2l+1}{4\pi}}\,, \tag{9.80}$$

$$\sum_{m=-l_{\max}}^{l_{\max}}\sum_{\lambda}\sum_{l=|m|}^{l_{\max}}\mathcal{B}^{m'}_{l'l}\frac{\psi_l^m(\lambda)}{\sqrt{\sigma_l}}(-1)^{l+l'}\Big[f_\lambda^{(+)m}(q)e^{-\mu_t\hat{k}_z(\lambda q)L/\lambda}d^l_{m',m}[i\tau(\lambda q)]$$
$$+f_\lambda^{(-)m}(q)(-1)^{l+m'}d^l_{m',-m}[i\tau(\lambda q)]\Big] = 0\,. \tag{9.81}$$

Coefficients $f_\lambda^{(\pm)m}$ are determined using (9.80) and (9.81). Noting that (9.80) and (9.81) are invariant with respect to the substitution $m' \to -m'$, we have

$$f_\lambda^{(\pm)-m}(q) = (-1)^m f_\lambda^{(\pm)m}(q). \tag{9.82}$$

Hence we can consider only nonnegative $m' \geq 0$ and $m \geq 0$. We have

$$\begin{pmatrix}\mathcal{M}^{(++)}(q) & \mathcal{M}^{(+-)}(q)\\ \mathcal{M}^{(-+)}(q) & \mathcal{M}^{(--)}(q)\end{pmatrix}\begin{pmatrix}f^{(+)}(q)\\ f^{(-)}(q)\end{pmatrix} = \begin{pmatrix}\mathrm{v}^{(+)}\\ \mathrm{v}^{(-)}\end{pmatrix}, \tag{9.83}$$

where

$$\mathcal{M}^{(++)}_{l'm',m\nu}(q) = \sum_{l=|m|}^{l_{\max}} c_{ll'm',m\nu}\,, \tag{9.84}$$

$$\mathcal{M}^{(+-)}_{l'm',m\nu}(q) = e^{-\mu_t\hat{k}_z(\lambda q)L/\lambda}(-1)^{m+m'}\sum_{l=|m|}^{l_{\max}}(-1)^l c_{ll'm',m\nu}\,, \tag{9.85}$$

$$\mathcal{M}^{(-+)}_{l'm',m\nu}(q) = e^{-\mu_t\hat{k}_z(\lambda q)L/\lambda}\sum_{l=|m|}^{l_{\max}}(-1)^{l+l'} c_{ll'm',m\nu}\,, \tag{9.86}$$

$$\mathcal{M}^{(--)}_{l'm',m\nu}(q) = (-1)^{l'+m'+m}\sum_{l=|m|}^{l_{\max}} c_{ll'm',m\nu}\,, \tag{9.87}$$

and

$$\mathrm{v}^{(+)}_{l'm'} = \delta_{m'0}\sum_{l=|m|}^{l_{\max}}\mathcal{B}^0_{l'l}\sqrt{\frac{2l+1}{4\pi}}\,,\quad \mathrm{v}^{(-)}_{l'm'} = 0\,. \tag{9.88}$$

Here,

$$c_{ll'm',m\lambda} = \mathcal{B}^{m'}_{l'l}\frac{\psi_l^m(\lambda)}{\sqrt{\sigma_l}}\Big[d^l_{m',m}[i\tau(\lambda q)] + (1-\delta_{m0})(-1)^m d^l_{m',-m}[i\tau(\lambda q)]\Big]\,. \tag{9.89}$$

The number of rows of $\mathcal{M}^{(++)}_{l'm',m\nu}$, $\mathcal{M}^{(+-)}_{l'm',m\nu}$, $\mathcal{M}^{(-+)}_{l'm',m\nu}$, $\mathcal{M}^{(--)}_{l'm',m\nu}$ is the number of pairs $(l', m')$. The number of columns of these matrices is the number of pairs $(m, \lambda)$. We note that equations with only odd $l'$ (or even $l'$) are linearly independent because of the three-term recurrence relations of associated Legendre polynomials (see (7.10)) (Liemert and Kienle, 2012b). Thus the numbers of equations and unknowns are equal in (9.83).

The specific intensity can be expressed as

$$I(\mathbf{r},\hat{\mathbf{s}};\mathbf{r}_0,\hat{\mathbf{s}}_0) = \frac{{\mu_t}^2 I_0}{2\pi}\sum_{l=|m|}^{l_{\max}}\sum_{m=-l}^{l}\frac{(-i)^m}{\sqrt{\sigma_l}}e^{-im\varphi_{\boldsymbol{\rho}-\boldsymbol{\rho}_0}}Y_{lm}(\hat{\mathbf{s}})K_{lm}(|\boldsymbol{\rho}-\boldsymbol{\rho}_0|,z)\,, \tag{9.90}$$

where $\varphi_{\boldsymbol{\rho}-\boldsymbol{\rho}_0}$ is the polar angle of the two-dimensional vector $\boldsymbol{\rho} - \boldsymbol{\rho}_0$ and

$$K_{lm}(|\boldsymbol{\rho} - \boldsymbol{\rho}_0|, z) = \int_0^\infty qJ_m(\mu_t q|\boldsymbol{\rho} - \boldsymbol{\rho}_0|) \sum_{m\geq 0} \psi_l^m(\lambda) \Big[e^{-\mu_t \hat{k}_z(\lambda q)z/\lambda} f_\lambda^{(+)m}(q) \\ +(-1)^{l+m\prime+m} e^{-\mu_t \hat{k}_z(\lambda q)(L-z)/\lambda} f_\lambda^{(-)m}(q)\Big] \\ \times \Big[d^l_{m',m}[i\tau(\lambda q)] + (1-\delta_{m0})(-1)^m d^l_{m',-m}[i\tau(\lambda q)]\Big]\, dq \qquad (9.91)$$

with $J_m$ the Bessel function of the first kind of order $m$. Hence the energy density $U$ is written as

$$U(\mathbf{r}; \mathbf{r}_0, \hat{\mathbf{s}}_0) = \frac{1}{c}\int_{\mathbb{S}^2} I(\mathbf{r}, \hat{\mathbf{s}}; \mathbf{r}_0, \hat{\mathbf{s}}_0) d\hat{\mathbf{s}} = \frac{I_0}{c\sqrt{\pi\sigma_0}} K_{00}(|\boldsymbol{\rho} - \boldsymbol{\rho}_0|, z)\,. \qquad (9.92)$$

Numerical results are shown below in Fig. 9.3. By setting $\mathbf{r}_0 = \mathbf{0}, \hat{\mathbf{s}}_0 = \hat{\mathbf{z}}$, we can write $U(\mathbf{r}; \mathbf{r}_0, \hat{\mathbf{s}}_0) = u(\rho, z)$ with $\rho = \sqrt{x^2 + y^2}$. In Fig. 9.3, the energy density $u(\rho, z)$ is compared to Monte Carlo simulation.

The specific intensity for the half-space geometry is obtained in the limit of $L \to \infty$. In the half space, the ballistic subtraction can be performed as follows. First we decompose the specific intensity into the ballistic and scattering terms as given in (7.21): $I(\mathbf{r}, \hat{\mathbf{s}}) = I_b(\mathbf{r}, \hat{\mathbf{s}}) + I_s(\mathbf{r}, \hat{\mathbf{s}})$. Then we further split $I_s(\mathbf{r}, \hat{\mathbf{s}})$ as

$$I_s(\mathbf{r}, \hat{\mathbf{s}}) = I_p(\mathbf{r}, \hat{\mathbf{s}}) + \psi(\mathbf{r}, \hat{\mathbf{s}})\,, \qquad (9.95)$$

where the particular solution $I_p(\mathbf{r}, \hat{\mathbf{s}})$ is defined for the whole space as

$$\hat{\mathbf{s}} \cdot \nabla I_p(\mathbf{r}, \hat{\mathbf{s}}) + \mu_t I_p(\mathbf{r}, \hat{\mathbf{s}}) = \mu_s \int_{\mathbb{S}^2} p(\hat{\mathbf{s}}, \hat{\mathbf{s}}') I_p(\mathbf{r}, \hat{\mathbf{s}}') d\hat{\mathbf{s}}' + \Theta(z) S[I_b](\mathbf{r}, \hat{\mathbf{s}}) \qquad (9.96)$$

for $(\mathbf{r}, \hat{\mathbf{s}}) \in \mathbb{R}^3 \times \mathbb{S}^2$. Here, $S[I_b](\mathbf{r}, \hat{\mathbf{s}})$ was defined in (7.24). It is possible to find $I_p(\mathbf{r}, \hat{\mathbf{s}})$ depending on the source term in the RTE for $I_b(\mathbf{r}, \hat{\mathbf{s}})$ (Liemert and Kienle, 2013b; Machida, Hoshi, Kagawa, and Takada, 2020). Since the RTE for the general solution $\psi(\mathbf{r}, \hat{\mathbf{s}})$ does not have the internal source term, $\psi(\mathbf{r}, \hat{\mathbf{s}})$ can be given as a superposition of eigenmodes. The coefficients of the superposition are determined using the boundary conditions. Thus, the specific intensity is obtained as $I(\mathbf{r}, \hat{\mathbf{s}}) = I_b(\mathbf{r}, \hat{\mathbf{s}}) + I_p(\mathbf{r}, \hat{\mathbf{s}}) + \psi(\mathbf{r}, \hat{\mathbf{s}})$.

The specific intensity was explicitly calculated in the case of $l_{\max} = 3$ (Liemert and Kienle, 2014).

### 9.6. Time-dependent case

Based on the spherical-harmonic expansion, solutions to the time-dependent three-dimensional RTE were obtained by Liemert and Kienle, 2012d and Liemert and Kienle, 2012e. In this section we solve the time-dependent RTE following Liemert and Kienle, 2012e. Let us consider the Green's function $G(\mathbf{r}, \hat{\mathbf{s}}, t; \mathbf{r}_0, \hat{\mathbf{s}}_0)$ in the whole space, which obeys

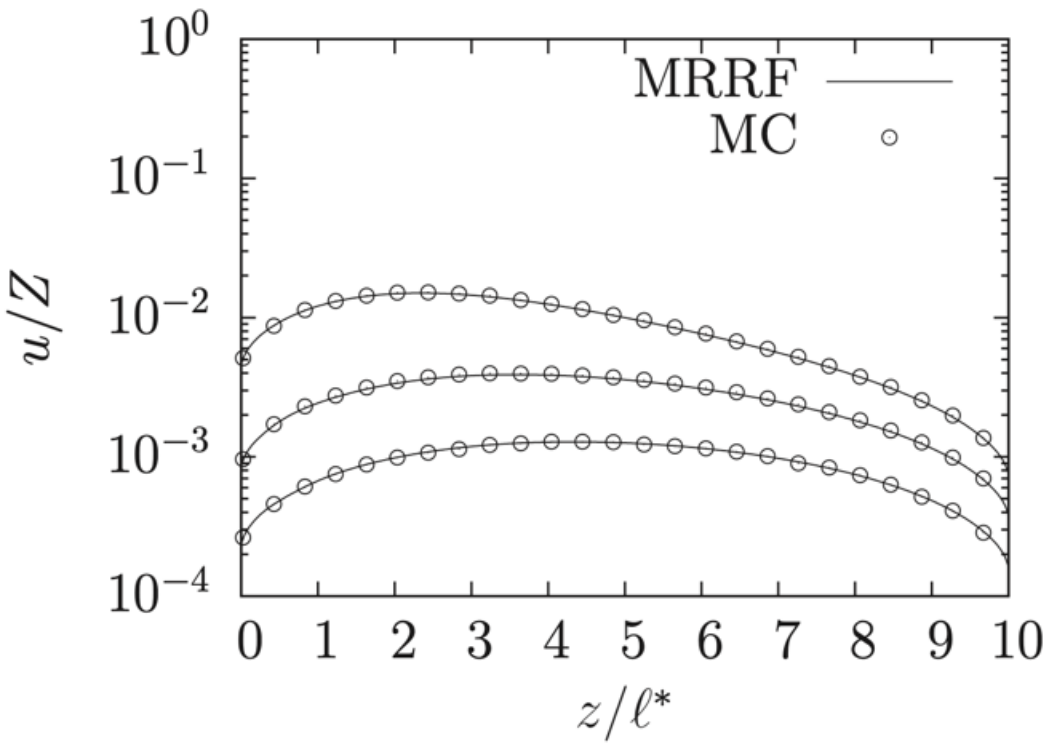

Figure 9.3. (Corrigendum: Machida, Panasyuk, Schotland, and Markel, 2010). The energy density $u(\rho, z)$ (see (9.92)) is plotted with open circles from Monte Carlo simulation. The width of the slab is $L = 10\ell^*$, where $\ell^* = (\mu_a + (1-\mathrm{g})\mu_s)^{-1}$. From the top to bottom, $\rho = 4\ell^*, 7\ell^*$, and $10\ell^*$. Normalization factor: $Z = (\ell^*)^{-2}$.

$$\left(\frac{1}{c}\frac{\partial}{\partial t} + \hat{\mathbf{s}}\cdot\nabla + \mu_t\right) G(\mathbf{r}, \hat{\mathbf{s}}, t; \mathbf{r}_0, \hat{\mathbf{s}}_0) = \mu_s \int_{\mathbb{S}^2} p(\hat{\mathbf{s}}, \hat{\mathbf{s}}') G(\mathbf{r}, \hat{\mathbf{s}}', t; \mathbf{r}_0, \hat{\mathbf{s}}_0)\, d\hat{\mathbf{s}}' + \delta(\mathbf{r}-\mathbf{r}_0)\delta(\hat{\mathbf{s}}-\hat{\mathbf{s}}_0)\delta(t). \tag{9.97}$$

The Fourier transform is given by

$$\tilde{G}(\mathbf{k}, \hat{\mathbf{s}}, t) = \int_{\mathbb{R}^3} e^{-i\mathbf{k}\cdot\mathbf{r}} G(\mathbf{r}, \hat{\mathbf{s}}, t; \mathbf{r}_0, \hat{\mathbf{s}}_0) d\mathbf{r}\ , \mathbf{k} \in \mathbb{R}^3. \tag{9.98}$$

Then we have

$$\left(\frac{1}{c}\frac{\partial}{\partial t} + i\mathbf{k}\cdot\hat{\mathbf{s}} + \mu_t\right) \tilde{G}(\mathbf{k}, \hat{\mathbf{s}}, t) = \mu_s \int_{\mathbb{S}^2} p(\hat{\mathbf{s}}, \hat{\mathbf{s}}') \tilde{G}(\mathbf{k}, \hat{\mathbf{s}}', t) d\hat{\mathbf{s}}' + e^{-i\mathbf{k}\cdot\mathbf{r}_0}\delta(\hat{\mathbf{s}}-\hat{\mathbf{s}}_0)\delta(t). \tag{9.99}$$

We express the Fourier transform of the Green's function as

$$\tilde{G}(\mathbf{k}, \hat{\mathbf{s}}, t) \approx \sum_{l=0}^{l_{\max}} \sum_{m=-l}^{l} f_{lm}(\mathbf{k}, t) \mathcal{R}_{\hat{\mathbf{k}}} Y_{lm}(\hat{\mathbf{s}})\ , \tag{9.100}$$

where

$$f_{lm}(\mathbf{k}, t) = \int_{\mathbb{S}^2} \tilde{G}(\mathbf{k}, \hat{\mathbf{s}}, t) \mathcal{R}_{\hat{\mathbf{k}}} Y_{lm}^*(\hat{\mathbf{s}})\, d\hat{\mathbf{s}}\ . \tag{9.101}$$

We note that $\hat{\mathbf{k}} = \mathbf{k}/|\mathbf{k}| \in \mathbb{R}^3$ in (9.100) and (9.101) is different from $\hat{\mathbf{k}}(\nu, \mathbf{q})$ in (4.32). Let us write the coefficients $f_{lm}(\mathbf{k}, t)$ in a vector form as $\mathbf{f}^m(\mathbf{k}, t) = \big(f_{lm}(\mathbf{k}, t)\big)$. Then the RTE can be expressed as

$$\frac{1}{c}\frac{d}{dt}\mathbf{f}^m(\mathbf{k}, t) + Q^m(k)\mathbf{f}^m(\mathbf{k}, t) = \mathbf{q}^m(\mathbf{k}, t)\,, \tag{9.102}$$

where matrix $Q^m(k)$ and vector $\mathbf{q}^m(\mathbf{k}, t)$ are given by

$$Q^m(k) = \begin{pmatrix} \sigma_l & ikb_l(m) & 0 & 0 & & \cdots & 0 \\ ikb_l(m) & \sigma_{l+1} & ikb_{l+1}(m) & 0 & & \cdots & \vdots \\ 0 & ikb_{l+1}(m) & \ddots & \ddots & & & 0 \\ 0 & 0 & \ddots & \ddots & & & 0 \\ \vdots & & & & & \sigma_{l_{\max}-1} & b_{l_{\max}}(m) \\ 0 & \cdots & 0 & 0 & & b_{l_{\max}}(m) & \sigma_{l_{\max}} \end{pmatrix}, \tag{9.103}$$

where $b_l(m)$ was given in (9.31) and

$$\mathbf{q}^m(\mathbf{k}, t) = \begin{pmatrix} e^{-i\mathbf{k}\cdot\mathbf{r}_0}\mathcal{R}_{\hat{\mathbf{k}}}Y^*_{|m|,m}(\hat{\mathbf{s}}_0)\delta(t) \\ \vdots \\ e^{-i\mathbf{k}\cdot\mathbf{r}_0}\mathcal{R}_{\hat{\mathbf{k}}}Y^*_{lm}(\hat{\mathbf{s}}_0)\delta(t) \\ \vdots \\ e^{-i\mathbf{k}\cdot\mathbf{r}_0}\mathcal{R}_{\hat{\mathbf{k}}}Y^*_{l_{\max},m}(\hat{\mathbf{s}}_0)\delta(t) \end{pmatrix}. \tag{9.104}$$

The solution of the above matrix-vector equation is obtained as

$$\mathbf{f}^m(\mathbf{k}, t) = ce^{-Q^m(k)ct}\mathbf{q}^m(\mathbf{k}, t) . \tag{9.105}$$

Hence,

$$\tilde{G}(\mathbf{k}, \hat{\mathbf{s}}, t) = ce^{-i\mathbf{k}\cdot\mathbf{r}_0} \sum_{m=-l_{\max}}^{l_{\max}} \sum_{l'=|m|}^{l_{\max}} \left(\mathcal{R}_{\hat{\mathbf{k}}}Y^*_{l'm}(\hat{\mathbf{s}}_0)\right) \left\{e^{-Q^m(k)ct}\right\}_{ll'} \left(\mathcal{R}_{\hat{\mathbf{k}}}Y_{lm}(\hat{\mathbf{s}})\right) . \tag{9.106}$$

We obtain

$$G(\mathbf{r}, \hat{\mathbf{s}}, t; \mathbf{r}_0, \hat{\mathbf{s}}_0) = \sum_{m=-l_{\max}}^{l_{\max}} \sum_{l'=|m|}^{l_{\max}} \left(\mathcal{R}_{\hat{\mathbf{R}}}Y^*_{l'm}(\hat{\mathbf{s}}_0)\right) \psi^m_{ll'}(R, t) \left(\mathcal{R}_{\hat{\mathbf{R}}}Y_{lm}(\hat{\mathbf{s}})\right) . \tag{9.107}$$

Here,

$$\psi^m_{ll'}(R, t) = c(-1)^m \sum_{m'=-\min(l,l')}^{\min(l,l')} (-1)^{m'} \sum_{L=|l-l'|}^{l+l'} C^{L,0}_{l,m,l',-m} C^{L,0}_{l,m',l',-m'} \times \frac{i^L}{2\pi^2} \int_0^\infty \left\{e^{-Q^m(k)ct}\right\}_{ll'} j_L(kR)k^2 \, dk , \tag{9.108}$$

where $j_L$ is the spherical Bessel function of the first kind. The ballistic subtraction is possible making use of this expression (Liemert and Kienle, 2012e).

Let us consider light propagation from a time-harmonic source:

$$\left(\frac{1}{c}\frac{\partial}{\partial t} + \hat{\mathbf{s}}\cdot\nabla + \mu_t\right) G(\mathbf{r}, \hat{\mathbf{s}}, t; \mathbf{r}_0, \hat{\mathbf{s}}_0) = \mu_s \int_{\mathbb{S}^2} p(\hat{\mathbf{s}}, \hat{\mathbf{s}}') G(\mathbf{r}, \hat{\mathbf{s}}', t; \mathbf{r}_0, \hat{\mathbf{s}}_0) \, d\hat{\mathbf{s}}' + \delta(\mathbf{r} - \mathbf{r}_0)\delta(\hat{\mathbf{s}} - \hat{\mathbf{s}}_0)e^{i\omega t} . \tag{9.109}$$

The solution can be written as

$$G(\mathbf{r}, \hat{\mathbf{s}}, t; \mathbf{r}_0, \hat{\mathbf{s}}_0) = I_\omega(\mathbf{r}, \hat{\mathbf{s}}; \mathbf{r}_0, \hat{\mathbf{s}}_0)e^{i\omega t} . \tag{9.110}$$

This $I_\omega$ satisfies

$$\left(\hat{\mathbf{s}}\cdot\nabla + \mu_t + i\frac{\omega}{c}\right) I_\omega(\mathbf{r}, \hat{\mathbf{s}}; \mathbf{r}_0, \hat{\mathbf{s}}_0) = \mu_s \int_{\mathbb{S}^2} p(\hat{\mathbf{s}}, \hat{\mathbf{s}}') I_\omega(\mathbf{r}, \hat{\mathbf{s}}'; \mathbf{r}_0, \hat{\mathbf{s}}_0) \, d\hat{\mathbf{s}}' + \delta(\mathbf{r} - \mathbf{r}_0)\delta(\hat{\mathbf{s}} - \hat{\mathbf{s}}_0). \tag{9.111}$$

Thus the problem reduces to a time-independent problem. Let us introduce

$$f_{lm}(\mathbf{k}; \omega) = \int_{\mathbb{S}^2} I_\omega(\mathbf{r}, \hat{\mathbf{s}}; \mathbf{r}_0, \hat{\mathbf{s}}_0)\mathcal{R}_{\hat{\mathbf{k}}}Y^*_{lm}(\hat{\mathbf{s}}) \, d\hat{\mathbf{s}} . \tag{9.112}$$

Let us define

$$\sigma_l(\omega) = \mu_a + \left(1 - \frac{\beta_l}{2l+1}\right)\mu_s + i\frac{\omega}{c}\,. \tag{9.113}$$

We obtain

$$\left(A^m(k;\omega) + i\frac{\omega}{c}\mathbb{I}_m\right)\mathbf{f}^m(\mathbf{k};\omega) = \mathbf{q}^m(\mathbf{k})\,, \tag{9.114}$$

where matrix $Q^m(k;\omega)$ is given in (9.103) but now $\sigma_l(\omega)$ are complex numbers. Vectors $\mathbf{q}^m(\mathbf{k})$ are given by

$$\mathbf{q}^m(\mathbf{k}) = \begin{pmatrix} e^{-i\mathbf{k}\cdot\mathbf{r}_0}\mathcal{R}_{\hat{\mathbf{k}}}Y^*_{|m|,m}(\hat{\mathbf{s}}_0) \\ \vdots \\ e^{-i\mathbf{k}\cdot\mathbf{r}_0}\mathcal{R}_{\hat{\mathbf{k}}}Y^*_{lm}(\hat{\mathbf{s}}_0) \\ \vdots \\ e^{-i\mathbf{k}\cdot\mathbf{r}_0}\mathcal{R}_{\hat{\mathbf{k}}}Y^*_{l_{\max},m}(\hat{\mathbf{s}}_0) \end{pmatrix}. \tag{9.115}$$

Let us consider the homogeneous equation. We have

$$\sum_{l'=|m|}^{l_{\max}} B_{ll'}(m;\omega)\psi^m_{l'}(\lambda) = \lambda\psi^m_l(\lambda). \tag{9.116}$$

Here, the matrix $B(m;\omega)$ was introduced in (9.38) but now it depends on $\omega$ because $\sigma_l(\omega)$ in $\beta_l(m)$ in (9.39) depends on $\omega$. We obtain

$$f^m_l(\mathbf{k};\omega) = e^{-i\mathbf{k}\cdot\mathbf{r}_0}\sum_{l'=|m|}^{l_{\max}}\sum_{\lambda}\frac{\psi^m_l(\lambda)\psi^m_{l'}(\lambda)^*}{(1+ik\lambda)\sqrt{\sigma_l(\omega)\sigma_{l'}(\omega)}}\mathcal{R}_{\hat{\mathbf{k}}}Y^*_{l'm}(\hat{\mathbf{s}}_0)\,. \tag{9.117}$$

Similar to the calculation by Panasyuk, Schotland, and Markel, 2006 for $\omega = 0$, the specific intensity ($\omega \neq 0$) can be expressed as

$$I_\omega(\mathbf{r},\hat{\mathbf{s}};\mathbf{r}_0,\hat{\mathbf{s}}_0) = \sum_{m=-l_{\max}}^{l_{\max}}\sum_{l'=|m|}^{l_{\max}}\left(\mathcal{R}_{\hat{\mathbf{R}}}Y^*_{l'm}(\hat{\mathbf{s}}_0)\right)\psi^m_{ll'}(R,\omega)\left(\mathcal{R}_{\hat{\mathbf{R}}}Y_{lm}(\hat{\mathbf{s}})\right)\,, \tag{9.118}$$

where

$$\psi^m_{ll'}(R,\omega) = \frac{c(-1)^m}{2\pi\sigma_l(\omega)\sigma_{l'}(\omega)}\sum_{M=-\min(l,l')}^{\min(l,l')}(-1)^M\sum_{L=|l-l'|}^{l+l'}C^{L,0}_{l,m,l',-m}C^{L,0}_{l,M,l',-M}$$
$$\times\sum_{\lambda>0}\frac{\psi^M_l(\lambda)\psi^M_{l'}(\lambda)^*}{\lambda^3}k_L\left(\frac{\mu_t R}{\lambda}\right)\,. \tag{9.119}$$

In Fig. 9.4, the Green's function $G(\mathbf{r},\hat{\mathbf{s}},t;\mathbf{r}_0,\hat{\mathbf{s}}_0)$ in (9.107) is plotted. Results were compared to Monte Carlo simulation.

Although the time-dependent RTE has not been solved in the presence of boundaries, Liemert and Kienle, 2013d approximately solved the time-dependent RTE in the half space using the Green's function for an infinite medium.

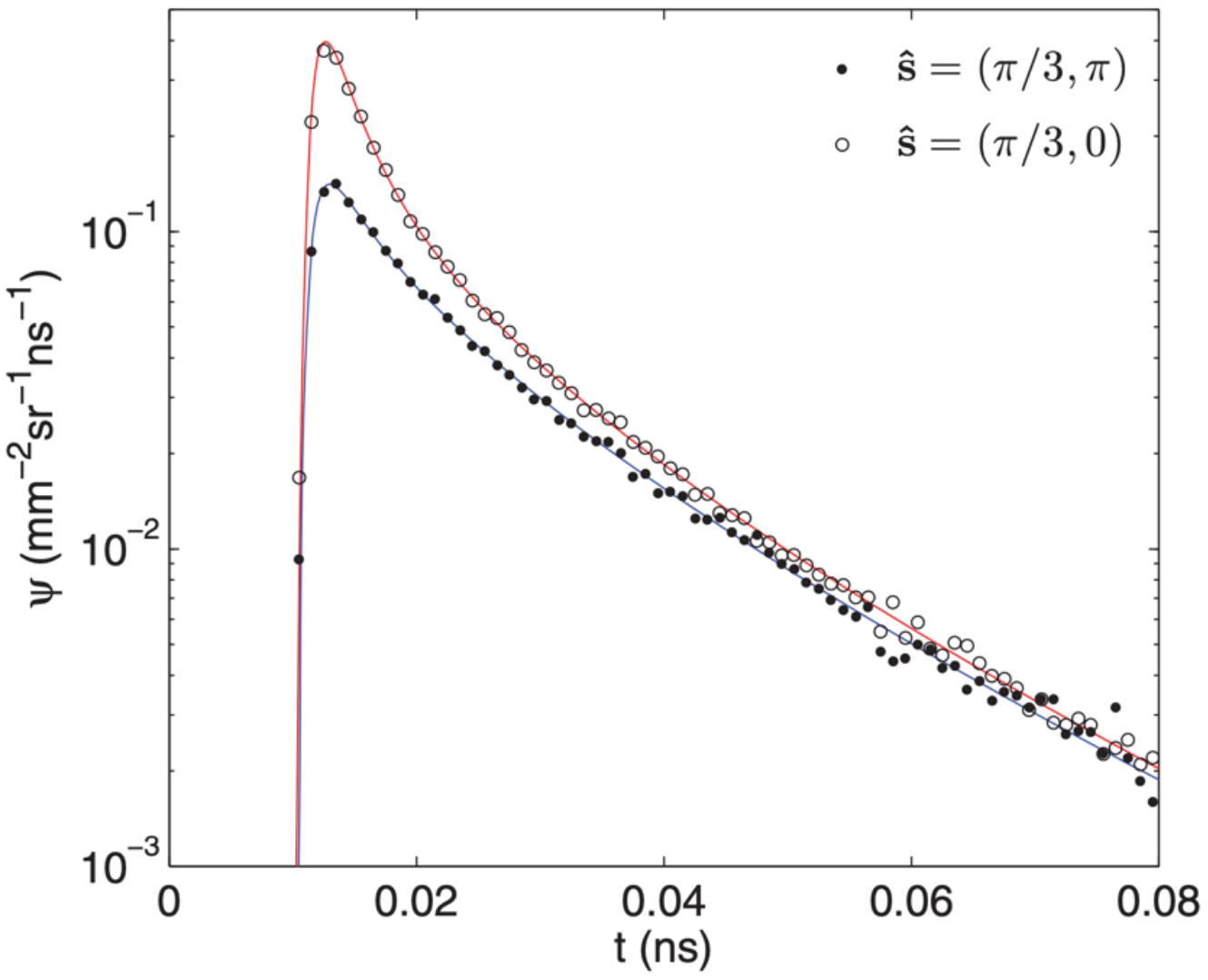

Figure 9.4. (Fig. 3 in Liemert and Kienle, 2012e). Time-resolved Green's function $\Psi = G(\mathbf{r}, \hat{\mathbf{s}}, t; \mathbf{r}_0, \hat{\mathbf{s}}_0)$ is shown for $(\theta_0, \varphi_0) = (\pi/4, 0)$. The observed position was chosen at $r = 3$ mm and $\hat{\mathbf{r}} = \hat{\mathbf{z}}$. Optical parameters were set to $\mu_a = 0.1$ mm$^{-1}$, $\mu_s = 10$ mm$^{-1}$, and g $= 0.9$. Furthermore, $c = 3 \times 10^8\ m/s$.

# 10. Applications

## 10.1. Optical tomography

Light propagation in biological tissue is governed by the RTE. On a large scale, say the propagation distance is 1cm when the transport mean free path is 1mm, the solution to the diffusion equation becomes almost identical to the solution of the RTE. Optical tomography is formulated as inverse problems which determine coefficients of the RTE or diffusion equation. The diffusion equation is usually used for optical tomography because the equation is easier to solve. Hence the tomography is often called diffuse optical tomography (Arridge, 1999; Boas, et al., 2001). The use of the RTE is difficult in general because the forward problem must be solved repeatedly when the inverse problem is solved in any ways (Arridge, 1999; Arridge and Schotland, 2009). Although the numerical cost of algorithms such as finite-difference and finite-element methods, and Monte Carlo simulation is high, one way to achieve optical tomography based on the RTE is to use rotated reference frames.

Here, we consider the reconstruction of $\mu_a$ in the RTE from boundary measurements (Machida, et al., 2016). When the reconstruction is performed, direct reconstruction methods offer an alternative approach to iterative algorithms (Arridge and Schotland, 2009). By direct reconstruction, we mean the use of inversion formulas and associated fast algorithms. In particular, we will employ the Rytov approximation.

In the experiment, the target (a lemon slice or a lotus root slice of thickness approximately 4mm) was placed in the middle of a sample chamber of width $L$ =1cm. The background medium consists of 0.7% Intralipid plus 0.07% India ink. Near-infrared light in this box obeys the RTE for the slab geometry. We assume that $\mu_s$ is constant everywhere in the medium but $\mu_a$ varies with position. Considering the Fresnel reflection on the boundaries, we write the RTE as

$$\hat{\mathbf{s}}\cdot\nabla I(\mathbf{r},\hat{\mathbf{s}})+(\mu_a(\mathbf{r})+\mu_s)I(\mathbf{r},\hat{\mathbf{s}})=\mu_s\int_{\mathbb{S}^2}p(\hat{\mathbf{s}},\hat{\mathbf{s}}')I(\mathbf{r},\hat{\mathbf{s}}')d\hat{\mathbf{s}}'+S(\mathbf{r},\hat{\mathbf{s}}) \qquad (10.1)$$

with the boundary condition

$$I(\mathbf{r},\hat{\mathbf{s}})=R(|\hat{\mathbf{s}}\cdot\hat{\mathbf{s}}'|)I(\mathbf{r},\hat{\mathbf{s}}'),\quad(\mathbf{r},\hat{\mathbf{s}})\in\Gamma_-, \qquad (10.2)$$

where $R(|\hat{\mathbf{s}}\cdot\hat{\mathbf{s}}'|)$ is the reflection coefficient and $\hat{\mathbf{s}}'$ is the incident direction corresponding to the reflected direction $\hat{\mathbf{s}}$. The notation $\Gamma_-$ means the set of incoming light:

$$\Gamma_-=\{\mathbf{r}\in\mathbb{R}^3,\hat{\mathbf{s}}\in\mathbb{S}^2;\ z=0\text{ or }L,\hat{\mathbf{s}}\cdot\hat{\boldsymbol{\nu}}<0\}, \qquad (10.3)$$

where $\hat{\boldsymbol{\nu}}$ is the outer unit normal vector on the boundary of the slab. We note that the laser beam which illuminates the sample chamber is described by the source term $S(\mathbf{r},\hat{\mathbf{s}})$. Let $\mathfrak{n}$ be the refractive index of the slab and we assume that there is air outside the slab. Assuming unpolarized light, the reflection coefficient $R(x)$ is given by

$$R(x)=\begin{cases}\dfrac{1}{2}\left[\left(\dfrac{x-\mathfrak{n}x_0}{x+\mathfrak{n}x_0}\right)^2+\left(\dfrac{x_0-\mathfrak{n}x}{x_0+\mathfrak{n}x}\right)^2\right], & x\ge x_c,\\ 1, & x<x_c,\end{cases} \qquad (10.4)$$

where

$$x_0=\sqrt{1-\mathfrak{n}^2(1-x^2)},\quad x_c=\frac{\sqrt{\mathfrak{n}^2-1}}{\mathfrak{n}}. \qquad (10.5)$$

Let us write

$$\mu_a(\mathbf{r})=\bar{\mu}_a+\delta\mu_a(\mathbf{r}). \qquad (10.6)$$

We introduce the Green's function $G(\mathbf{r},\hat{\mathbf{s}};\mathbf{r}',\hat{\mathbf{s}}')$ as

$$\hat{\mathbf{s}}\cdot\nabla G(\mathbf{r},\hat{\mathbf{s}};\mathbf{r}',\hat{\mathbf{s}}')+\bar{\mu}_tG(\mathbf{r},\hat{\mathbf{s}};\mathbf{r}',\hat{\mathbf{s}}')=\mu_s\int_{\mathbb{S}^2}p(\hat{\mathbf{s}},\hat{\mathbf{s}}'')G(\mathbf{r},\hat{\mathbf{s}}'';\mathbf{r}',\hat{\mathbf{s}}')\,d\hat{\mathbf{s}}''$$
$$+\delta(\mathbf{r}-\mathbf{r}')\delta(\hat{\mathbf{s}}-\hat{\mathbf{s}}'), \qquad (10.7)$$

where $\bar{\mu}_t=\bar{\mu}_a+\mu_s$. The boundary condition is given by

$$G(\mathbf{r},\hat{\mathbf{s}};\mathbf{r}',\hat{\mathbf{s}}')=R(|\hat{\mathbf{s}}\cdot\hat{\mathbf{s}}''|)G(\mathbf{r},\hat{\mathbf{s}}'';\mathbf{r}',\hat{\mathbf{s}}'),\quad(\mathbf{r},\hat{\mathbf{s}})\in\Gamma_-, \qquad (10.8)$$

where $\hat{\mathbf{s}}''$ is the incident direction corresponding to the reflected direction $\hat{\mathbf{s}}$. By the spherical-harmonic expansion in rotated reference frames, the Green's function is given by (see (9.21) and (9.90))

$$G(\mathbf{r},\hat{\mathbf{s}};\mathbf{r}',\hat{\mathbf{s}}')=\frac{1}{(2\pi)^2}e^{-i\mathbf{q}\cdot(\boldsymbol{\rho}-\boldsymbol{\rho}')}\sum_{l=0}^{l_{\max}}\sum_{m=-l}^{l}k_{lm}(\mathbf{q},z)Y_{lm}(\hat{\mathbf{s}})Y_{lm}^*(\hat{\mathbf{s}}'), \qquad (10.9)$$

where $\mathbf{r}=(\boldsymbol{\rho},z)$ and

$$k_{lm}(\mathbf{q},z)=(-1)^me^{-im\varphi_q}\sum_{M=0}^{l_{\max}}\sum_{n}\frac{\psi_l^M(\lambda_n^M)}{\sqrt{\sigma_l}}\left[e^{-\bar{\mu}_t\hat{k}_z(\lambda_n^Mq)z/\lambda_n^M}f^{(+)}(\lambda_n^M,q)\right.$$
$$\left.+(-1)^{l+m+M}e^{-\bar{\mu}_t\hat{k}_z(\lambda_n^Mq)(L-z)/\lambda_n^M}f^{(-)}(\lambda_n^M,q)\right]$$
$$\times\left[d_{mM}^l[i\tau(\lambda_n^Mq)]+(1-\delta_{M0})(-1)^Md_{m,-M}^l[i\tau(\lambda_n^Mq)]\right]. \qquad (10.10)$$

Let $I_0(\mathbf{r},\hat{\mathbf{s}})$ be the solution to (10.1) in which $\mu_a(\mathbf{r})$ is replaced by $\bar{\mu}_a$. We note the identity,

$$I(\mathbf{r},\hat{\mathbf{s}}) = I_0(\mathbf{r},\hat{\mathbf{s}}) - \int_{\mathbb{S}^2}\int_{\Omega} G(\mathbf{r},\hat{\mathbf{s}};\mathbf{r}',\hat{\mathbf{s}}')\delta\mu_a(\mathbf{r}')I(\mathbf{r}',\hat{\mathbf{s}}')d\mathbf{r}'d\hat{\mathbf{s}}' \ . \tag{10.11}$$

Then the first-order Rytov approximation is given by

$$-\ln\left(\frac{I(\mathbf{r},\hat{\mathbf{s}})}{I_0(\mathbf{r},\hat{\mathbf{s}})}\right) = \frac{1}{I_0(\mathbf{r},\hat{\mathbf{s}})}\int_{\mathbb{S}^2}\int_{\Omega} G(\mathbf{r},\hat{\mathbf{s}};\mathbf{r}',\hat{\mathbf{s}}')\delta\mu_a(\mathbf{r}')I_0(\mathbf{r}',\hat{\mathbf{s}}')d\mathbf{r}'d\hat{\mathbf{s}}' \ . \tag{10.12}$$

Hereafter we assume a point source oriented in the inward normal direction located on the $z = 0^+$ plane. Thus,

$$I_0(\mathbf{r},\hat{\mathbf{s}}) = G(\mathbf{r},\hat{\mathbf{s}};\mathbf{r}_s,\hat{\mathbf{z}}), \ \ \mathbf{r}_s = \begin{pmatrix}\boldsymbol{\rho}_s\\ 0^+\end{pmatrix} . \tag{10.13}$$

Light exiting the slab in the outward normal direction is collected by a point detector that is located on the plane $z = L$. Suppose that the outgoing light is detected at $\mathbf{r}_d$, where

$$\mathbf{r}_d = \begin{pmatrix}\boldsymbol{\rho}_d\\ L\end{pmatrix} . \tag{10.14}$$

It will prove convenient to introduce the data function $\mathcal{D}(\boldsymbol{\rho}_s,\boldsymbol{\rho}_d)$, which is given by

$$\mathcal{D}(\boldsymbol{\rho}_s,\boldsymbol{\rho}_d) = -I_0(\mathbf{r},\hat{\mathbf{s}})\ln\left(\frac{I(\mathbf{r},\hat{\mathbf{s}})}{I_0(\mathbf{r},\hat{\mathbf{s}})}\right) . \tag{10.15}$$

We have

$$\mathcal{D}(\boldsymbol{\rho}_s,\boldsymbol{\rho}_d) = \int_{\mathbb{S}^2}\int_{\Omega} G(\mathbf{r}_d,\hat{\mathbf{z}};\mathbf{r}',\hat{\mathbf{s}}')G(\mathbf{r}',\hat{\mathbf{s}}';\mathbf{r}_d,-\hat{\mathbf{z}})\delta\mu_a(\mathbf{r}')\,d\mathbf{r}'d\hat{\mathbf{s}}' \ . \tag{10.16}$$

The data function is written as

$$\begin{aligned}\mathcal{D}(\boldsymbol{\rho}_s,\boldsymbol{\rho}_d) = \frac{1}{(2\pi)^4}\int_{\Omega}\int_{\mathbb{R}^2}\int_{\mathbb{R}^2} & e^{i(\mathbf{q}_1-\mathbf{q}_2)\cdot\boldsymbol{\rho}}e^{-i(\mathbf{q}_1\cdot\boldsymbol{\rho}_s-\mathbf{q}_2\cdot\boldsymbol{\rho}_d)}\\ & \times\kappa(\mathbf{q}_1,\mathbf{q}_2,z)\delta\mu_a(\mathbf{r})\,d\mathbf{q}_1d\mathbf{q}_2d\mathbf{r} .\end{aligned} \tag{10.17}$$

where

$$\kappa(\mathbf{q}_1,\mathbf{q}_2,z) = \sum_{l=0}^{l_{\max}}\sum_{m=-l}^{l}(-1)^m k_{lm}(\mathbf{q}_1,z)k_{lm}^*(\mathbf{q}_2,L-z) \ . \tag{10.18}$$

Now we have a linear inverse problem to reconstruct $\delta\mu_a(\mathbf{r})$ using measured data $\mathcal{D}(\boldsymbol{\rho}_s,\boldsymbol{\rho}_d)$. After moving to the Fourier space for $(\boldsymbol{\rho}_s,\boldsymbol{\rho}_d)$, the inversion can be done with singular value decomposition (SVD). When the inverse problem is solved, regularization is necessary to obtain stable solutions $\delta\mu_a(\mathbf{r})$. We employ the step function as the regularizer (truncated SVD).

In Figure 10.1, reconstructed $\delta\mu_a(\mathbf{r})$ is shown for (upper panels) a lemon slice and (lower panels) a lotus root slice. The field of view is 9 cm $\times$ 9 cm. The number above each panel shows the distance from the source plane. The parameters were set to

$$\bar{\mu}_a = 0.4\ \mathrm{cm}^{-1}\ , \ \ \mu_s = 20\ \mathrm{cm}^{-1}\ , \ \ \mathrm{g} = 0.65\ . \tag{10.19}$$

Datasets were acquired by raster scanning the beam over a 29 × 29 square lattice with a lattice spacing 3.21 mm (841 source positions). For each source, a 397 × 397 pixel region of interest is read out from the CCD (157609 detectors). In total, a dataset of $1.3 \times 10^8$ source-detector pairs was acquired. In the numerical calculation, we set $l_{\max} = 9$.

The transport-based optical tomography was also formulated with the three-dimensional $F_N$ method (Machida, 2019).

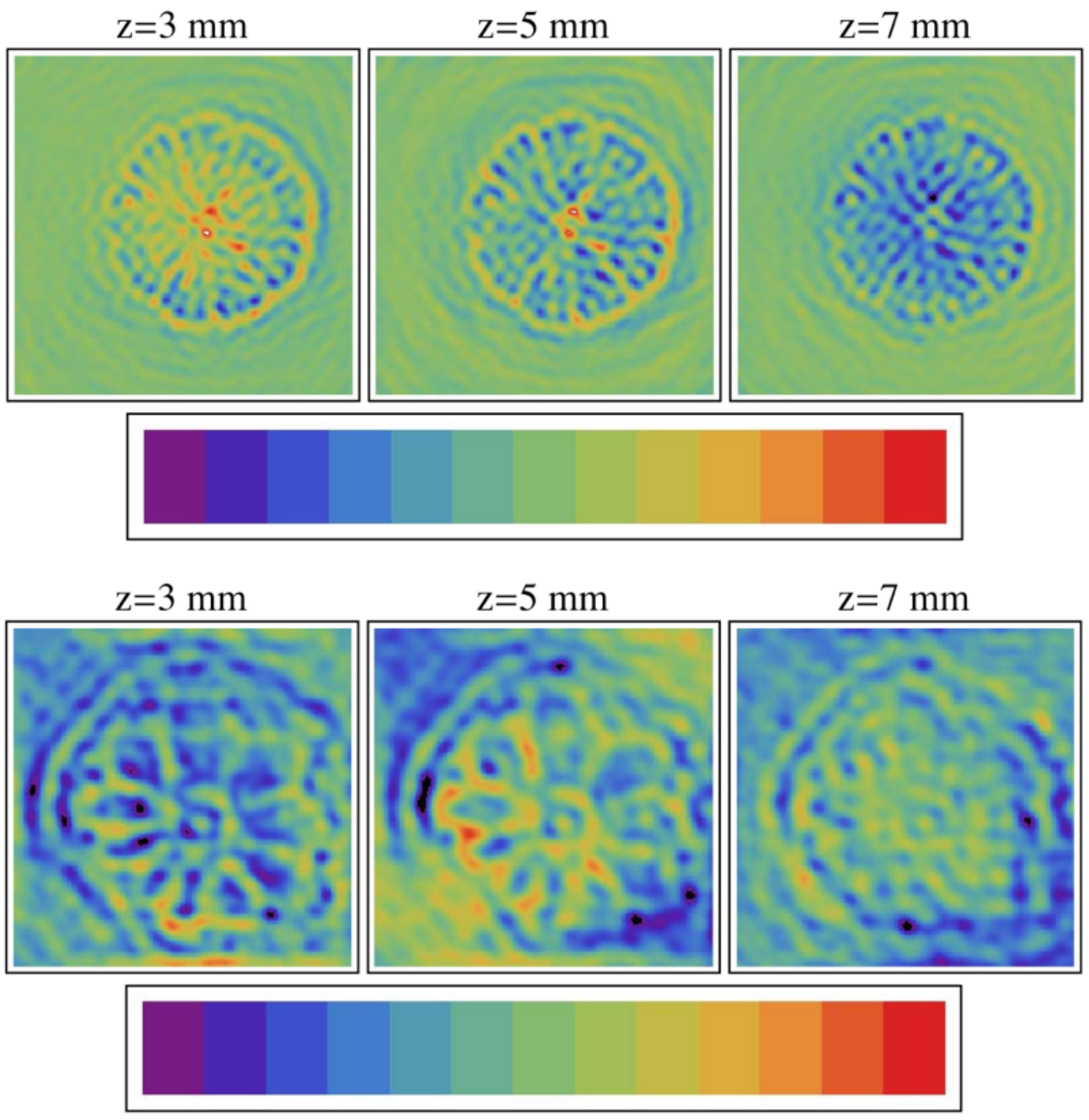

Figure 10.1 (Figs. 7 and 8 in Machida, et al., 2016). Reconstructions of a lemon slice (upper panels) and a lotus root slice (lower panels). The number above each panel shows the distance from the source plane.

## 10.2. Parameter identification

We here consider the spatial-frequency domain imaging (SFDI), for which light illumination in the spatial-frequency domain is used. We first explain that the penetration distance becomes shorter as the spatial frequency $q_0$ gets larger.

Let us begin with the following expression of the specific intensity for the half space (see (7.18), (7.19)).

$$
\begin{aligned}
I(\mathbf{r}, \hat{\mathbf{s}}) = e^{i\mu_t \mathbf{q}_0 \cdot \boldsymbol{\rho}} \sum_{M=-l_{\max}}^{l_{\max}} \Bigg[ \sum_{j=0}^{J^M-1} a_j^M \Psi^M_{\nu_j(M)}(\hat{\mathbf{s}}, \mathbf{q}_0) e^{-\mu_t \hat{k}_z(\nu_j(M) q_0) z/\nu_j(M)} \\
+ \int_0^1 a^M(\nu) \Psi^M_\nu(\hat{\mathbf{s}}, \mathbf{q}_0) e^{-\mu_t \hat{k}_z(\nu q_0) z/\nu}\, d\nu \Bigg],
\end{aligned}
\tag{10.20}
$$

where $a_j^M$ and $a^M(\nu)$ are coefficients which depend on the boundary conditions. Here, $\Psi^M_{\nu_j(M)}(\hat{\mathbf{s}}, \mathbf{q}_0)$ and $\Psi^M_\nu(\hat{\mathbf{s}}, \mathbf{q}_0)$ are three-dimensional singular eigenfunctions:

$$
\Psi^M_\nu(\hat{\mathbf{s}}, \mathbf{q}_0) = \mathcal{R}_{\hat{\mathbf{k}}(\nu, \mathbf{q}_0)} \Phi^M_\nu(\hat{\mathbf{s}}). \tag{10.21}
$$

We note that $\nu_0 = \nu_0^0 > 1$ is the largest eigenvalue. Let us define

$$I_{\rm as}(\mathbf{r},\hat{\mathbf{s}}) = e^{-i\mu_t\mathbf{q}_0\cdot\boldsymbol{\rho}}a_0^0\Psi_{\nu_0}^M(\hat{\mathbf{s}},\mathbf{q}_0)e^{-\mu_t\hat{k}_z(\nu_0 q_0)z/\nu_0}\,. \tag{10.22}$$

For large $z$, the behavior of $I(\mathbf{r},\hat{\mathbf{s}})$ can be investigated using $I_{\rm as}(\mathbf{r},\hat{\mathbf{s}})$ because as $z\to\infty$,

$$\int_{\mathbb{R}^2}\int_{\mathbb{S}^2}|I(\mathbf{r},\hat{\mathbf{s}})-I_{\rm as}(\mathbf{r},\hat{\mathbf{s}})|\,d\hat{\mathbf{s}}d\boldsymbol{\rho} = o\left(\exp\left(-\mu_t z\sqrt{\nu_0^{-2}+q_0^2}\right)\right)\,. \tag{10.23}$$

We can estimate $\nu_0$ as (Machida, et al., 2009)

$$\frac{\sqrt{1+\eta}}{\sqrt{3\frac{\mu_a}{\mu_t}\left(1-\mathrm{g}\frac{\mu_s}{\mu_t}\right)}} \le \nu_0 \le \frac{1+\sqrt{\eta}}{\sqrt{3\frac{\mu_a}{\mu_t}\left(1+\mathrm{g}\frac{\mu_s}{\mu_t}\right)}}\,, \tag{10.24}$$

where

$$\eta = \frac{4}{5}\frac{\mu_a}{\mu_a+\mu_s(1-\mathrm{g}^2)}\,. \tag{10.25}$$

For $\mu_a \ll \mu_s$, which is typical in biological tissue, we have

$$\frac{1}{\nu_0} \approx \sqrt{\frac{3\mu_a}{\mu_t}(1-\varpi\mathrm{g})} \approx \frac{1}{\mu_t}\sqrt{3\mu_a\mu_s'}\,, \tag{10.26}$$

where

$$\mu_s' = (1-\mathrm{g})\mu_s\,. \tag{10.27}$$

Thus, if $\mu_a$ is small, the specific intensity decays as

$$I \sim e^{-z\sqrt{3\mu_a\mu_s'+(\mu_t q_0)^2}}\,. \tag{10.28}$$

That is, the penetration depth becomes shallow for nonzero $q_0$. This facts motivates us to investigate optical properties of the top layer of a layered medium by SFDI.

In Machida, et al., 2020, $\mu_a$ and $\mu_s$ of the top layer of a four-layer solid phantom was experimentally determined using the specific intensity (10.20). In Nothelfer, et al., 2019, the correction for surface scattering in SFDI was considered based on (10.20) for precise retrieval of the bulk optical properties. The method was validated using phantoms with different surface roughness.

# 11. Concluding remarks

When the three-dimensional RTE is solved by separation of variables, separated solutions have a connection to separated solutions for the one-dimensional RTE via rotated reference frames. By making use of this connection, the three-dimensional RTE can be solved using the knowledge of the one-dimensional transport theory.

It is certainly a limitation that coefficients $\mu_a, \mu_s$ must be constant. One of the reasons that the RTE with constant coefficients is still worth being studied is benchmark (Ganapol, 2008). Since the analytical and numerical solutions presented in this review are quite accurate, they can be used for benchmark. When absorption and scattering coefficients spatially vary, it is possible to treat the spatially varying part as perturbation. In 10.1, $\mu_a(\mathbf{r})$ is reconstructed using solutions of the RTE obtained with rotated reference frames.

## Acknowledgements

The author started the research on the RTE and rotated reference frames when he was a postdoc of John C. Schotland and Vadim A. Markel, which he greatly appreciates. To deepen the theory on the RTE with rotated reference frames, fruitful discussion on the occasion of the International Conference on Transport Theory (ICTT), which have been held every other year, was crucial. This work was supported by JSPS KAKENHI Grant Number JP17K05572 and JST PRESTO Grant Number JPMJPR2027.

# A. Singular eigenfunctions

In this Appendix, we take $1/\mu_t$ as the unit of length. Then the one-dimensional homogeneous RTE is written as

$$\mu\frac{\partial}{\partial z}I(z,\hat{\mathbf{s}})+I(z,\hat{\mathbf{s}})=\varpi\int_{\mathbb{S}^2}p(\hat{\mathbf{s}},\hat{\mathbf{s}}')I(z,\hat{\mathbf{s}}')\,d\hat{\mathbf{s}}'\ ,\qquad (A.1)$$

where $z\in\mathbb{R}$, $\hat{\mathbf{s}}\in\mathbb{S}^2$. Separated solutions to (A.1) are given by (Case, 1960; McCormick and Kuščer, 1966; Mika, 1961)

$$I(z,\hat{\mathbf{s}})=\Phi_\nu^m(\hat{\mathbf{s}})e^{-z/\nu}\ ,\qquad (A.2)$$

where $\nu\in\mathbb{R}$ is a separation constant and $\Phi_\nu^m(\hat{\mathbf{s}})$ is introduced (4.14). We note $-l_{\max}\le m\le l_{\max}$. For $\Phi_\nu^m(\hat{\mathbf{s}})$, (4.17) is derived. We can rewrite (4.17) as

$$\mathcal{K}\Phi_\nu^m(\hat{\mathbf{s}})=\nu\Phi_\nu^m(\hat{\mathbf{s}})\ ,\qquad (A.3)$$

where operator $\mathcal{K}$ was defined by

$$(\mathcal{K}^{-1}\psi)(\hat{\mathbf{s}})=\frac{1}{\mu}\left[\psi(\hat{\mathbf{s}})-\varpi\int_{\mathbb{S}^2}p(\hat{\mathbf{s}},\hat{\mathbf{s}}')\psi(\hat{\mathbf{s}}')\,d\hat{\mathbf{s}}'\right],\quad 0<|\mu|\le 1\ .\qquad (A.4)$$

for a function $\psi(\hat{\mathbf{s}})$. The operator $\mathcal{K}$ was introduced in Larsen, 1974. As described in Sec. 4.4, the separation constant $\nu$ takes discrete values $\pm\nu_j^m$ ($\nu_j^m>1, j=0,1,\cdots,M^m$) and real numbers in the interval between $-1$ and 1.

Singular eigenfunctions are given by (4.25). The function $\lambda^m(\nu)$ is given by

$$\lambda^m(\nu)=1-\frac{\varpi\nu}{2}\mathcal{P}\int_{-1}^{1}\frac{g^m(\nu,\mu)}{\nu-\mu}(1-\mu^2)^{|m|}\,d\mu\ .\qquad (A.5)$$

For $\nu_j^m$, we have (4.47).

The normalization factor (4.8) can be expressed as follows (McCormick and Kuščer, 1966): for $\nu=\nu_j^m$,

$$\mathcal{N}^m(\nu)=\frac{1}{2}\nu^2g^m(\nu,\nu)\left.\frac{d\Lambda^m(w)}{dw}\right|_{w=\nu}\ ,\qquad (A.6)$$

and for $\nu\in(-1,1)$,

$$\mathcal{N}^m(\nu)=\nu\Lambda^{m+}(\nu)\Lambda^{m-}(\nu)(1-\nu^2)^{-|m|}\ ,\qquad (A.7)$$

where

$$\Lambda^{m\pm}(\nu)=\lim_{\epsilon\to 0^+}\Lambda^m(\nu\pm i\epsilon)\ .\qquad (A.8)$$

## B. Eq. (6.4)

We have

$$\left(\mathcal{R}_{\hat{\mathbf{k}}} Y_{l_1 m_1}(\hat{\mathbf{s}})\right)\left(\mathcal{R}_{\hat{\mathbf{k}}} Y_{l_2 m_2}(\hat{\mathbf{s}})\right) = \sum_{m_1'=-l_1}^{l_1} D^{l_1}_{m_1' m_1}(\varphi_{\hat{\mathbf{k}}}, \theta_{\hat{\mathbf{k}}}, 0) Y_{l_1 m_1'}(\hat{\mathbf{s}}) \sum_{m_2'=-l_2}^{l_2} D^{l_2}_{m_2' m_2}(\varphi_{\hat{\mathbf{k}}}, \theta_{\hat{\mathbf{k}}}, 0) Y_{l_2 m_2'}(\hat{\mathbf{s}})$$
$$= \sum_{m_1'=-l_1}^{l_1} \sum_{m_2'=-l_2}^{l_2} D^{l_1}_{m_1' m_1}(\varphi_{\hat{\mathbf{k}}}, \theta_{\hat{\mathbf{k}}}, 0) D^{l_2}_{m_2' m_2}(\varphi_{\hat{\mathbf{k}}}, \theta_{\hat{\mathbf{k}}}, 0)$$
$$\times \sum_{l=|l_1-l_2|}^{l_1+l_2} \sum_{m=-l}^{l} \sqrt{\frac{(2l_1+1)(2l_2+1)}{4\pi(2l+1)}} C^{l0}_{l_1 0 l_2 0} C^{lm}_{l_1 m_1' l_2 m_2'} Y_{lm}(\hat{\mathbf{s}}) \,, \qquad (B.1)$$

where we used (Sec. 5.6 of Varshalovich, Moskalev, and Khersonskii, 1988)

$$Y_{l_1 m_1'}(\hat{\mathbf{s}}) Y_{l_2 m_2'}(\hat{\mathbf{s}}) = \sum_{l=|l_1-l_2|}^{l_1+l_2} \sum_{m=-l}^{l} \sqrt{\frac{(2l_1+1)(2l_2+1)}{4\pi(2l+1)}} C^{l0}_{l_1 0 l_2 0} C^{lm}_{l_1 m_1' l_2 m_2'} Y_{lm}(\hat{\mathbf{s}}) \,. \qquad (B.2)$$

Thus we can write

$$\left(\mathcal{R}_{\hat{\mathbf{k}}} Y_{l_1 m_1}(\hat{\mathbf{s}})\right)\left(\mathcal{R}_{\hat{\mathbf{k}}} Y_{l_2 m_2}(\hat{\mathbf{s}})\right) = \sum_{l=|l_1-l_2|}^{l_1+l_2} \sum_{m=-l}^{l} f^{lm}_{l_1 m_1, l_2 m_2} Y_{lm}(\hat{\mathbf{s}}) \,, \qquad (B.3)$$

where

$$f^{lm}_{l_1 m_1, l_2 m_2} = \sqrt{\frac{(2l_1+1)(2l_2+1)}{4\pi(2l+1)}}$$
$$\times \sum_{m_1'=-l_1}^{l_1} \sum_{m_2'=-l_2}^{l_2} D^{l_1}_{m_1' m_1}(\varphi_{\hat{\mathbf{k}}}, \theta_{\hat{\mathbf{k}}}, 0) D^{l_2}_{m_2' m_2}(\varphi_{\hat{\mathbf{k}}}, \theta_{\hat{\mathbf{k}}}, 0) C^{l0}_{l_1 0 l_2 0} C^{lm}_{l_1 m_1' l_2 m_2'} \,. \qquad (B.4)$$

We note that (see Sec. 4.6 of Varshalovich, Moskalev, and Khersonskii, 1988)

$$D^{l_1}_{m_1' m_1}(\varphi_{\hat{\mathbf{k}}}, \theta_{\hat{\mathbf{k}}}, 0) D^{l_2}_{m_2' m_2}(\varphi_{\hat{\mathbf{k}}}, \theta_{\hat{\mathbf{k}}}, 0) = \sum_{L=|l_1-l_2|}^{l_1+l_2} \sum_{M=-L}^{L} \sum_{M'=-L}^{L} C^{LM}_{l_1 m_1' l_2 m_2'}$$
$$\times D^{L}_{MM'}(\varphi_{\hat{\mathbf{k}}}, \theta_{\hat{\mathbf{k}}}, 0) C^{LM'}_{l_1 m_1 l_2 m_2} \,. \qquad (B.5)$$

We can write

$$f^{lm}_{l_1 m_1, l_2 m_2} = \sqrt{\frac{(2l_1+1)(2l_2+1)}{4\pi(2l+1)}} C^{l0}_{l_1 0 l_2 0} \sum_{L=|l_1-l_2|}^{l_1+l_2} \sum_{M=-L}^{L} \sum_{M'=-L}^{L} D^{L}_{MM'}(\varphi_{\hat{\mathbf{k}}}, \theta_{\hat{\mathbf{k}}}, 0) C^{LM'}_{l_1 m_1 l_2 m_2}$$
$$\times \sum_{m_1'=-l_1}^{l_1} \sum_{m_2'=-l_2}^{l_2} C^{LM}_{l_1 m_1' l_2 m_2'} C^{lm}_{l_1 m_1' l_2 m_2'} \,. \qquad (B.6)$$

We note that $C_{l_1 0 l_2 0}^{l0}$ vanishes unless $|l_1 - l_2| \leq l \leq l_1 + l_2$ and moreover (Sec. 8.1 of Varshalovich, Moskalev, and Khersonskii, 1988),

$$\sum_{m_1'=-l_1}^{l_1} \sum_{m_2'=-l_2}^{l_2} C_{l_1 m_1' l_2 m_2'}^{LM} C_{l_1 m_1' l_2 m_2'}^{lm} = \delta_{Ll} \delta_{Mm} \, . \qquad (B.7)$$

Hence we obtain

$$f_{l_1 m_1, l_2 m_2}^{lm} = \sqrt{\frac{(2l_1 + 1)(2l_2 + 1)}{4\pi(2l + 1)}} C_{l_1 0 l_2 0}^{l0} \sum_{M'=-l}^{l} D_{mM'}^{l}(\varphi_{\hat{\mathbf{k}}}, \theta_{\hat{\mathbf{k}}}, 0) C_{l_1 m_1 l_2 m_2}^{LM'} \, . \qquad (B.8)$$

Thus we arrive at

$$\left(\mathcal{R}_{\hat{\mathbf{k}}} Y_{l_1 m_1}(\hat{\mathbf{s}})\right)\left(\mathcal{R}_{\hat{\mathbf{k}}} Y_{l_2 m_2}(\hat{\mathbf{s}})\right) = \sum_{l=|l_1-l_2|}^{l_1+l_2} \sum_{m'=-l}^{l} \sum_{m=-l}^{l} \sqrt{\frac{(2l_1 + 1)(2l_2 + 1)}{4\pi(2l + 1)}}$$

$$\times C_{l_1 0 l_2 0}^{l0} C_{l_1 m_1 l_2 m_2}^{Lm} D_{m'm}^{l}(\varphi_{\hat{\mathbf{k}}}, \theta_{\hat{\mathbf{k}}}, 0) Y_{lm'}(\hat{\mathbf{s}})$$

$$= \sum_{l=|l_1-l_2|}^{l_1+l_2} \sum_{m=-l}^{l} \sqrt{\frac{(2l_1 + 1)(2l_2 + 1)}{4\pi(2l + 1)}} C_{l_1 0 l_2 0}^{l0} C_{l_1 m_1 l_2 m_2}^{Lm} \mathcal{R}_{\hat{\mathbf{k}}} Y_{lm}(\hat{\mathbf{s}}) \, . \qquad (B.9)$$

Using (B.2),

$$\left(\mathcal{R}_{\hat{\mathbf{k}}} Y_{l_1 m_1}(\hat{\mathbf{s}})\right)\left(\mathcal{R}_{\hat{\mathbf{k}}} Y_{l_2 m_2}(\hat{\mathbf{s}})\right) = \mathcal{R}_{\hat{\mathbf{k}}}\left(Y_{l_1 m_1}(\hat{\mathbf{s}}) Y_{l_2 m_2}(\hat{\mathbf{s}})\right).$$